\documentclass[10pt,reqno]{amsart}
\usepackage{amsmath} 
\usepackage{amssymb}
\usepackage{dsfont}
\usepackage[dvips,draft,final]{graphics}
\usepackage[T1]{fontenc}
\usepackage{fancyhdr}
\usepackage{url}
\usepackage[colorlinks,linktocpage,linkcolor=blue]{hyperref}

\usepackage{color}
\usepackage{graphicx}
\newtheorem{thm}{Theorem}[section]
\newtheorem{lem}[thm]{Lemma}
\newtheorem{cor}[thm]{Corollary}
\newtheorem{prop}[thm]{Proposition}

\newtheorem{rmk}[thm]{Remark}
\newtheorem{exa}[thm]{Example}
\newtheorem{assumption}[thm]{Assumption}
\numberwithin{equation}{section}

\newcommand{\R}{\mathbb{R}}

\allowdisplaybreaks

\def\epsilon{\varepsilon}
\def\phi {\varphi}

\definecolor{gray}{RGB}{168,168,168}

\providecommand{\abs}[1]{\left\lvert#1\right\rvert}
\providecommand{\norm}[1]{\left\lVert#1\right\rVert}

\renewcommand{\leq}{\leqslant}
\renewcommand{\geq}{\geqslant}
\renewcommand{\d}{{\rm d}}
\providecommand{\abs}[1]{\left\lvert#1\right\rvert}
\providecommand{\norm}[1]{\left\lVert#1\right\rVert}
\def\thefootnote{{}}

\title[ ]{}

\title[Reconstruction and well-posedness of an inverse coefficient problem]{Reconstruction and well-posedness of an inverse coefficient problem for a diffusion equation$^*$}

\author[Siyu Cen, Bangti Jin, Yavar KIAN and Zhi Zhou]{ Siyu Cen$^1$, Bangti Jin$^2$, Yavar KIAN$^1$ and Zhi Zhou$^3$}

\date{}

\begin{document}
\begin{abstract}  In this work, we investigate the inverse coefficient problem of determining an unknown coefficient in a diffusion equation from a single boundary flux measurement. More precisely, assuming that the domain is cylindrical, we develop a unified theoretical and numerical framework for the stable reconstruction of a space-time-dependent coefficient, by reformulating the inverse problem as a nonlinear equation. The inherent nonlinearity of the problem is transferred to an auxiliary nonlinear equation, enabling the application of techniques and concepts from nonlinear analysis. By combining representation formulas with regularity properties of solutions to diffusion equations, we extend to nonlinear inverse problems an approach that has so far been restricted to linear ones. The theoretical reconstruction procedure is complemented by a numerical analysis of the proposed method, including convergence and error estimates, as well as numerical experiments. The proposed reconstruction algorithm is directly inspired by the theoretical framework. \vskip5pt

  \vskip5pt


\end{abstract}
\maketitle

\renewcommand{\thefootnote}{\fnsymbol{footnote}}
\footnotetext{\hspace*{-5mm} 
\begin{tabular}{@{}r@{}p{16cm}@{}}
$^*$
&The work of S. C. and Y. K. is supported by the French National Research Agency ANR and Hong Kong RGC Joint Research Scheme for the project IdiAnoDiff (grant ANR-24-CE40-7039). The work of B. J. is supported by Hong Kong RGC General Research Fund (Project 14306824), ANR / Hong Kong RGC Joint Research Scheme (A-CUHK402/24). The work of Z. Z. is supported by National Natural Science Foundation of China (Project 12422117 and Project 12426312) and Hong Kong Research Grants Council (15302323), and an internal grant of Hong Kong Polytechnic University (Project ID: P0053938, Work Programme: 4-ZZVA).\\
$^1$ 
& Univ Rouen Normandie, CNRS, Normandie Univ, LMRS UMR 6085, F-76000 Rouen, France  (\texttt{siyu.cen@univ-rouen.fr}, \texttt{yavar.kian@univ-rouen.fr})\\
$^2$
& Department of Mathematics, The Chinese University of Hong Kong, Shatin, New Territories, Hong Kong, P.R. China (\texttt{b.jin@cuhk.edu.hk, bangti.jin@gmail.com})\\
$^3$
& Department of Applied Mathematics, The Hong Kong Polytechnic University, Kowloon, Hong Kong, China. (\texttt{zhizhou@polyu.edu.hk})
\end{tabular}}

\section{Introduction}\label{sec:intro}

This work is devoted to the reconstruction
of a coefficient in a diffusion equation from a single boundary measurement. 
More precisely, let $\Omega=\omega\times (0,\ell)$ with $\omega$ being a   bounded domain of $\R^{d-1}$, $d=2,3$, with a $C^4$ boundary, and consider the following initial-boundary value problem (IBVP)
\begin{equation}\label{eqn:forward}
\left\{\begin{aligned}
\partial_t u-\Delta u+q(t,x')u&=0,\quad &&(t,x',x_d)=(t,x)\in (0,T)\times \Omega,\\  
u(t,x', \ell)&=g_0(t,x'),\quad &&(t,x')\in(0,T)\times \omega,\\
u(t,x)&=g_1(t,x) ,\quad &&(t,x)\in(0,T)\times \partial\omega\times(0,\ell),\\ 
\partial_{x_d}u(t,x', 0)&=0,\quad &&(t,x')\in(0,T)\times \omega,\\
u(0,x)&=u_0(x),\quad &&x\in\Omega.
\end{aligned}\right.
\end{equation}
We investigate the following inverse coefficient problem:
\begin{itemize}
\item[{\bf(ICP)}]{\em Reconstruct the  coefficient $q$ from the boundary  measurement of the flux on the top face $\omega\times\{\ell\}$ of the cylinder $\Omega$ associated with problem \eqref{eqn:forward}.}
\end{itemize}

Problem \eqref{eqn:forward} governs the evolution of the diffusion concentration of a particle located at position $x$ at time $t$, and arises in a wide range of applications, including physics, biology, and finance. For instance, $q$ may represent the radiative coefficient in heat conduction \cite{Yamamoto:2001} or the perfusion coefficient in human physiology \cite{Pennes:1948,Trucu:2010}. Problem \eqref{eqn:forward} can also be interpreted as a linearized reaction--diffusion equation, where $q$ corresponds to the linear reaction term. Such models appear in numerous fields, including chemistry, biology, geology, physics, and ecology, with applications ranging from population dynamics \cite{FIS} and Rayleigh--B\'enard convection \cite{NW} to combustion theory \cite{YBZFK38}. In this context, \textbf{(ICP)} consists in reconstructing the coefficient $q$, which characterizes the underlying physical law governing diffusion in the medium, from measurements of the flux generated by the diffusion process \eqref{eqn:forward} on a portion of the boundary.

Motivated by its intrinsic nonlinearity, ill-posedness, and numerous applications, inverse coefficient problems similar to \textbf{(ICP)} have attracted considerable attention within the mathematical community. Significant progress has been made on fundamental questions such as uniqueness \cite{Fei,I1,Isakov:1991,KLLY,SU} and on related problems \cite{I2,KiUh}. However, other important issues, including stability and reconstruction, remain only partially understood. Nevertheless, a number of studies have provided partial answers under various measurement settings, including boundary observations, final-time measurements, and combinations thereof.  The stable reconstruction of space-dependent coefficients has been investigated using both final-time observations \cite{ZhangZhangZhou:2022} and various types of boundary and internal measurements \cite{Yamamoto:2001}.  To the best of our knowledge, the most general stability result for the recovery of a  coefficient of a parabolic equation depending on both time and space variables, using excitations and measurements restricted to the boundary, was established in \cite{CK2}. We mention also the work \cite[Corollary 1.5.]{KianYamamoto:2019}, where a Lipschitz stability estimate was derived.

In contrast to the  literature on stability, reconstruction methods for inverse coefficient problems such as \textbf{(ICP)} have mainly been developed from a numerical perspective and are generally based on available stability estimates (see, e.g., \cite{JinShinZhou:2023,ZhangZhangZhou:2022}). Overall, despite being a central topic in inverse problem theory (see, e.g., \cite{CGR,FKOU,Ga,Na}), theoretically justified reconstruction procedures for inverse coefficient problems remain poorly understood and are often associated with unstable and difficult-to-implement algorithms. These limitations are closely related to the nonlinear nature of \textbf{(ICP)} and to the instability inherent in many inverse coefficient problems \cite{KRS}. In contrast, for linear inverse problems, such as inverse source problems, several stable and theoretically justified reconstruction procedures have been developed. Most of these approaches rely on reformulating the inverse problem as an operator equation whose well-posedness can be established in a suitable functional framework \cite{JKZ,KianYamamoto:2019}. Such theoretically grounded reconstruction methods can substantially improve numerical performance by providing additional analytical properties, including rigorous error estimates \cite{CJKZ:2026}.

The main objective of this work is to provide a comprehensive treatment of \textbf{(ICP)} by developing both theoretical and numerical reconstruction methods endowed with stability guarantees. In particular, one of our primary goals is to fill a gap in the existing literature by establishing, to the best of our knowledge, the first theoretically justified reconstruction procedure for recovering coefficients from boundary measurements of a partial differential equation through the well-posedness analysis of a nonlinear equation.

Our analysis relies on a detailed study of estimates, regularity properties, and representation formulas for solutions to \eqref{eqn:forward} corresponding to coefficients $q$ belonging either to $C([0,T];H^s(\omega))$, with $s\in(0,1/2]$, or to $C([0,T];L^2(\omega))\cap L^\infty(0,T;L^\infty(\omega))$. These results are obtained by combining several analytical ingredients, including specific properties of cylindrical domains and product estimates in Sobolev spaces. By exploiting these tools together with suitable regularity assumptions, we reformulate the two variants of \textbf{(ICP)} in Section \ref{main}, in terms of solving the nonlinear equations \eqref{t2a} and \eqref{t3a}, respectively. We establish local well-posedness under a suitable smallness assumption on either the measurements or the final time $T$, and we further prove global well-posedness under additional structural assumptions. Our approach follows the general philosophy underlying  methods for nonlinear problems, e.g., nonlinear partial differential equations (see, e.g., \cite{CH}). To the best of our knowledge, this is the first work to employ such a methodology for the analysis and reconstruction of \textbf{(ICP)}.
 
In addition to the reconstruction procedures derived from the well-posedness analysis of \eqref{t2a} and \eqref{t3a}, we establish two Lipschitz stability estimates for problem \textbf{(ICP)} under weak \textit{a priori} assumptions, with explicit dependence of the stability constants on the problem parameters. Furthermore, our theoretical results are complemented by a numerical analysis based on the Galerkin finite element method (FEM) and the backward Euler scheme for the discretization of the nonlinear operators arising in the reconstruction process. We prove linear convergence of the resulting iterative reconstruction algorithm in the discrete setting and derive \textit{a priori} error estimates for the reconstructed coefficient.

The remainder of this article is organized as follows. In Section \ref{main}, we present the main results, and in Section \ref{sec:prelim}, we collect several preliminary results on the direct problem and investigate the properties of solutions to \eqref{eqn:forward}. In Section \ref{sec:recon_local}, we present the local well-posedness analysis, stated in Theorem \ref{t2}, of the nonlinear integral equations associated with \textbf{(ICP)}. In Section \ref{sec:stab}, we prove the stability results stated in Theorem \ref{t4} and Corollary \ref{c3}, and in Section \ref{sec:recon_global}, prove global well-posedness under suitable \textit{a priori} assumptions on the coefficient $q$ stated in Theorem \ref{t3}. We derive explicit reconstruction formulas and the error analysis in Section \ref{sec:error}, and in Section \ref{sec:num} we present numerical experiments to illustrate the theoretical findings.

\section{Main results}\label{main}

We present two formulations of  \textbf{(ICP)}, each leading to a reconstruction procedure based on solving a nonlinear equation, and investigate their well-posedness. The first formulation yields a general reconstruction procedure for coefficients $q\in C([0,T];H^{s_1}(\omega))$, with $s_1\in(0,\frac12]$, and we also establish Lipschitz stability estimates. The second formulation of \textbf{(ICP)} is for coefficients
$q\in C([0,T];L^2(\omega))\cap L^\infty(0,T;L^\infty(\omega))$. Under additional structural assumptions, we prove global well-posedness of the associated nonlinear equation and derive the reconstruction procedure.

\subsection{First reconstruction result}
Fix $0<s_0<s_1\leq \frac{1}{2}$ and  $q_1\in C([0,T];H^{s_1}(\omega))$ a known reference coefficient. Fix also $q_2=q_1+q$, $q\in C([0,T];H^{s_1}(\omega))$. For $j=1,2$, consider the problem
\begin{equation}\label{eq3}
\left\{\begin{aligned}\partial_t u_j-\Delta u_j+q_j(t,x')u_j=0,\quad &(t,x)\in 
(0,T)\times \Omega,\\  
u_j(t,x', \ell)=1\quad &(t,x')\in(0,T)\times \omega,\\
u_j(t,x)=g(t,x) ,\quad &(t,x)\in(0,T)\times \partial\omega\times(0,\ell),\\ 
\partial_{x_d}u_j(t,x', 0)=0\quad &(t,x')\in(0,T)\times \omega,\\
u_j(0,x)=u_0(x),\quad &x\in\Omega,\end{aligned}\right.
\end{equation}
with $g$, $u_0$ such that there exists $G\in   C^1([0,T];H^{s_0}(\Omega))\cap  C([0,T];H^{2+s_0}(\Omega)) $ satisfying
\begin{equation}\label{trace}
\left\{\begin{aligned}
G(t,x', \ell)=1\quad &(t,x')\in(0,T)\times \omega,\\
G(t,x)=g(t,x) ,\quad &(t,x)\in(0,T)\times \partial\omega\times(0,\ell),\\ 
\partial_{x_d}G(t,x', 0)=0\quad &(t,x')\in(0,T)\times \omega,\\
G(0,x)=u_0(x),\quad &x\in\Omega.\end{aligned}\right.\end{equation}

We will prove in Theorem \ref{t1} that, for $j=1,2$, problem \eqref{eq3} admits a unique solution $u_j\in C([0,T];H^2(\Omega))\cap C^1([0,T];L^2(\Omega))$. We study the following inverse problem:
\begin{itemize}
\item[{\bf(ICP1)}]{\em Reconstruct  the unknown coefficient $q$ from the difference of the flux measurement between the system \eqref{eq3}, with $j=1,2$, given by $\partial_{x_d}u_1(t,x',\ell)-\partial_{x_d}u_2(t,x',\ell)$, $(t,x')\in(0,T)\times\omega$.}
\end{itemize}

To solve \textbf{(ICP1)}, we formulate it as the resolution of a nonlinear integral equation, and then prove its well-posedness under the following regularity assumption.

\begin{assumption}\label{assum:recon_loc}
    Fix $m:(t,x')\mapsto \partial_{x_d}u_1(t,x',\ell)-\partial_{x_d}u_2(t,x',\ell)$. There exists $H_m\in  C^1([0,T];H^{s_0}(\Omega))\cap  C([0,T];H^{2+s_0}(\Omega)) $ satisfying
\begin{equation}\label{trace2}
\left\{\begin{aligned}
H_m(t,x', \ell)=m(t,x')\quad &(t,x')\in(0,T)\times \omega,\\
H_m(t,x)=0 ,\quad &(t,x)\in(0,T)\times \partial\omega\times(0,\ell),\\ 
H_m(t,x', 0)=0\quad &(t,x')\in(0,T)\times \omega,\\
H_m(0,x)=0,\quad &x\in\Omega.\end{aligned}\right.\end{equation}
\end{assumption}
Our first main result can be stated as follows.

\begin{thm}\label{t2} Fix $s_1\in(0,\frac12]$, and let Assumption \ref{assum:recon_loc}  be fulfilled. Then, there exists a  map $\mathcal M_{q_1}$, taking values in $C([0,T];H^{s_1}(\omega))$, and a map  $\mathcal K_{q_1}:C([0,T];H^{s_1}(\omega))\to C([0,T];H^{s_1}(\omega))$, both depending only on $\Omega$, $u_0$, $g$ and $q_1$,  such that \emph{\textbf{(ICP1)}} for $q\in C([0,T];H^{s_1}(\omega))$ can be reformulated  as solving  in  $C([0,T];H^{s_1}(\omega))$  the following nonlinear equation
\begin{equation}\label{t2a}q=\mathcal M_{q_1}(m)+\mathcal K_{q_1}(q).\end{equation}
Moreover,  there exists $T_0>0$, depending on $P=\norm{\mathcal M_{q_1}(m)}_{C([0,T];H^{s_1}(\omega))}$, $\Omega$, $g$, $u_0$, $s_1$ and $q_1$,
 such that, for $T\in(0,T_0)$, equation \eqref{t2a} is well-posed in $C([0,T];H^{s_1}(\omega))$ for $m$ equipped with the seminorm $m\mapsto\norm{\mathcal M_{q_1}(m)}_{C([0,T];H^{s_1}(\omega))}$. Similarly, for all $T>0$, there exists $P_0$  depending on $T$, $\Omega$, $g$, $u_0$ and $q_1$,
 such that, for $P\in(0,P_0)$, the 
nonlinear equation \eqref{t2a} is well-posed in $C([0,T];H^{s_1}(\omega))$.
\end{thm}

To the best of our knowledge, Theorem \ref{t2} provides the first positive answer to  \textbf{(ICP1)} based on its reformulation as equation \eqref{t2a} and the subsequent analysis of its well-posedness. The explicit form of equation \eqref{t2a} will be given in Section \ref{sec:recon_local}. The well-posedness analysis of equation \eqref{t2a} is carried out under additional assumptions related either to the length of the time interval $(0,T)$ or to the size of the measurements $P$. This result is consistent with classical nonlinear analysis theory, where these properties are referred to as local well-posedness and arise in situations where phenomena such as finite-time blow-up may occur.
Overall, the analysis of  \textbf{(ICP1)} and its reformulation as the nonlinear equation \eqref{t2a} share similarities with methods from nonlinear analysis (see, e.g., \cite{CH}), and Theorem \ref{t2} appears to be the first result establishing such a connection. We believe that the methodology developed in the proof of Theorem \ref{t2} may have applications beyond  \textbf{(ICP1)}.

The proof of Theorem \ref{t2} relies on a suitable representation of solutions of \eqref{eq1} that essentially exploits the axial-coordinate independence of the coefficient $q$. In this representation, the coefficient is transferred into the nonlinear equation \eqref{t2a}, whose well-posedness is subsequently analyzed. Among other things, we heavily exploit the properties of solutions to parabolic equations with coefficient $q$ lying in the low-regularity Sobolev space $C([0,T];H^s(\omega))$, $s \in (0, \frac12]$, which appears to be completely new in the setting of inverse coefficient problems similar to \textbf{(ICP1)}.

Assumption \ref{assum:recon_loc} is a regularity condition required for reformulating \textbf{(ICP1)} as an integral equation. In Section \ref{sec:stab}, we show that Assumption \ref{assum:recon_loc} is satisfied under suitable regularity assumptions on the parameters $q$ and $q_1$. Thus, the reconstruction method can only be applied to coefficients possessing sufficient regularity, a requirement that is indirectly encoded in Assumption \ref{assum:recon_loc} on the measurements. From a practical viewpoint, the reconstruction procedure, which reformulates  \textbf{(ICP1)} as equation \eqref{t2a}, can only be applied when the measurements exhibit an appropriate level of regularity. The numerical analysis of \textbf{(ICP)} in Section \ref{sec:error} indicates that such limitations are not restrictive and are, to some extent, expected for a reliable reconstruction procedure.

In the spirit of reconstruction methods for inverse source problems \cite{JKZ,KianYamamoto:2019}, the result of Theorem \ref{t2} can, under suitable assumptions, be extended to a more general class of parabolic equations with variable coefficients. In addition, the class of Dirichlet excitations considered in \eqref{eq3} can be substantially enlarged. We have chosen to work with a simplified model in order to avoid unnecessary technical difficulties and to focus primarily on the solution of \textbf{(ICP1)}.

\subsection{Stability estimates}

In addition to the reconstruction result in Theorem \ref{t2}, we also study the stability issue for \textbf{(ICP1)}. We give explicit conditions ensuring the regularity condition of Assumption \ref{assum:recon_loc}. More precisely, consider $s_1\in(0,\frac12]$, $s_0\in (0,s_1)$, $q\in C([0,T];H^{s_1}(\omega))$ and  the problem
\begin{equation}\label{eq10}
\left\{\begin{aligned}\partial_t u-\Delta u+q(t,x')u=0,\quad &(t,x)\in 
(0,T)\times \Omega,\\  
u(t,x', \ell)=1\quad &(t,x')\in(0,T)\times \omega,\\
u(t,x)=g(t,x) ,\quad &(t,x)\in(0,T)\times \partial\omega\times(0,\ell),\\ 
\partial_{x_d}u(t,x', 0)=0\quad &(t,x')\in(0,T)\times \omega,\\
u(0,x)=u_0(x),\quad &x\in\Omega,\end{aligned}\right.
\end{equation}
with $g$, $u_0$ such that there exists $G\in  C^2([0,T];H^{s_0}(\Omega))\cap C^1([0,T];H^{2+s_0}(\Omega))\cap  C([0,T];H^{4}(\Omega))$ satisfying the following compatibility conditions
\begin{equation}\label{trace_stab}
\left\{\begin{aligned}
G(t,x', \ell)=1,\ \quad &(t,x')\in(0,T)\times \omega,\\
G(t,x)=g(t,x) ,\quad &(t,x)\in(0,T)\times \partial\omega\times(0,\ell),\\ 
\partial_{x_d}G(t,x', 0)=0\quad &(t,x')\in(0,T)\times \omega,\\
G(0,x)=u_0(x),\quad &x\in\Omega,\\
\partial_t G(0,x)=\Delta u_0(x)-q(0,x')u_0(x),\quad &x\in\Omega,\\
\partial_{x_d}^2 G(t,x',\ell)=q(t,x'),\quad&(t,x')\in(0,T)\times\omega,\\
 \partial_{x_d}^3G(t,x',0)=0,\quad&(t,x')\in(0,T)\times\omega,\\
\sum_{j=1}^{d-1} \partial_{x_j}^2G(t,x)=\partial_tg(t,x)-\partial_{x_d}^2g(t,x)+q(t,x')g(t,x) ,\quad &(t,x)\in(0,T)\times \partial\omega\times(0,\ell).\end{aligned}\right.\end{equation}
We prove in Proposition \ref{p1} that for $q\in C^1([0,T];H^{s_1}(\omega))\cap C([0,T];W^{2,\infty}(\omega))$, problem \eqref{eq10} admits a unique solution $u\in C^1([0,T];H^2(\Omega))\cap C([0,T];H^4(\Omega))$. We derive the following stability estimate for \textbf{(ICP1)}.

\begin{thm}\label{t4} Let $s_1\in(0,\frac12]$, $s_0\in(0,s_1)$ and, for $j=1,2$, let  $q_j\in C^1([0,T];H^{s_1}(\omega))\cap C([0,T];W^{2,\infty}(\omega))$ and $u_j\in C^1([0,T];H^2(\Omega))\cap C([0,T];H^4(\Omega))$ solve \eqref{eq10}. Assume also that there exists 
 $M>0$ such that
 \begin{equation}\label{t4a}\norm{q_j}_{L^\infty(0,T;H^{s_1}(\omega))}\leq M,\quad j=1,2.\end{equation}
Then, there exist $C,c>0$ depending on $\Omega$, $g$, $u_0$, $s_0$, $s_1$,  $T$, and  $\beta_1,\beta_2\in(0,1)$ depending on $s_0,s_1$, and $C'>0$ depending on $s_0$, $s_1$ and $T$ such that
 \begin{equation}\label{t4b}\begin{aligned}&\norm{q_1-q_2}_{L^\infty(0,T;H^{s_1}(\omega))}\leq \left(\sum_{k=1}^{\infty}\frac{(C\kappa(M))^{k+1}}{\Gamma (k\beta_1)}\right)\left(\norm{m}_{W^{1,\infty}(0,T; H^{s_0}(\omega))}+\norm{m}_{L^\infty(0,T;H^{\frac{3}{2}+s_0}(\omega))}\right),\end{aligned}\end{equation}
with
$$\kappa(M)=(1+M^4)e^{cM^2}\left(1+\sum_{k=1}^{\infty}\frac{(C'M)^{k}}{\Gamma (k\beta_2)}\right).$$
 
 \end{thm}
 
We have also a conditional stability under a priori estimates in $L^\infty(0,T;L^\infty(\omega))$.
\begin{cor}\label{c3} Let $s_1\in(0,\frac12)$, $q_j\in C^1([0,T];H^{s_1}(\omega))\cap C([0,T];W^{2,\infty}(\omega))$ and let $u_j\in C^1([0,T];H^2(\Omega))\cap C([0,T];H^4(\Omega))$ solve \eqref{eq10}. Assume also that there exists 
 $M>0$ such that
 \begin{equation}\label{t5a}\norm{q_j}_{L^\infty(0,T;L^\infty(\omega))}\leq M,\quad j=1,2.\end{equation}
Then there exist $C,c>0$ depending on $\Omega$, $g$, $u_0$, $s_1$ and $T$, and $\beta_1,\beta_2\in(0,1)$ depending on $s_1$ and $C'>0$ depending on $ s_1,T$ such that \eqref{t4b} holds true.\end{cor}

The conditional stability estimates in Theorem \ref{t4} and Corollary \ref{c3} are byproducts of the reconstruction procedure established in Theorem \ref{t2}, where Assumption \ref{assum:recon_loc} is replaced by suitable regularity assumptions on the unknown parameters $q_j$, $j=1,2$. Nevertheless, Theorem \ref{t4} and Corollary \ref{c3} provide several improvements over existing results \cite{CK2,KianYamamoto:2019,Yamamoto:2001,ZhangZhangZhou:2022}. More precisely, they substantially relax both the regularity assumptions and the \textit{a priori} conditions required for the Lipschitz stability estimate \eqref{t4b}. Moreover,  the dependence of the stability constants on the parameter $M$, appearing in the \textit{a priori} conditions \eqref{t4a} and \eqref{t5a}, is made explicit. These properties make Theorem \ref{t4} and Corollary \ref{c3} more precise and potentially more suitable for applications.

\subsection{Second reconstruction result}

Now we present an alternative formulation of \textbf{(ICP)} and consider a reconstruction method based on the global well-posedness of a nonlinear equation similar to \eqref{t2a}, under an additional boundedness constraint. Namely, let $q^\dagger\in L^\infty(0,T;L^\infty(\omega))$ be the exact coefficient and let $u^\dagger=u(q^\dagger)$ denote the solution to problem \eqref{eq10} with $q=q^\dagger$. Consider the following:
\begin{itemize}
\item[{\bf(ICP2)}]{\em Reconstruct the coefficient $q^\dagger$ from the measurement $m(t,x'):=\partial_{x_d}u^\dagger(t,x',\ell)$.}
\end{itemize}

Like in Theorem \ref{t2}, we analyze   \textbf{(ICP2)} with extra regularity assumptions. 
\begin{assumption}\label{assum:alg}
   There exists $G_m\in C([0,T];H^2(\Omega))\cap C^1([0,T];L^2(\Omega))$ satisfying 
        \begin{equation*} 
     \left\{\begin{aligned}
         G_m(0,x)&=\partial_{x_d} u_0(x),&& x\in \Omega\\
        G_m(t,x',x_d)&=\partial_{x_d} g(t,x',x_d), && (t,x',x_d)\in   [0,T]\times\partial\omega\times(0,\ell),\\
        G_m(t,x',0)&=0,&&    (t,x')\in     [0,T] \times\omega,\\
        G_m(t,x',\ell)&=m(t,x'),&&  (t,x')\in  [0,T]\times \omega.\\
      \end{aligned}\right.
    \end{equation*} 
\end{assumption}Note that by the trace theorem, Assumption \ref{assum:alg} implicitly encodes $m\in C([0,T];H^\frac{3}{2}(\omega))$.

For $b>0$, we denote by $\mathcal{Q}_b$ the closed convex subset of $L^\infty(0,T;L^\infty(\omega))$ defined by
\begin{equation}\label{eqn:ad_set}
    \mathcal{Q}_b=\{q\in L^\infty(0,T;L^\infty(\omega))\,:\, \|q\|_{L^\infty(0,T;L^\infty(\omega))}\le b\}.
\end{equation} 
We define also the  projection $P_{\mathcal{Q}_b}$ onto  $\mathcal{Q}_b$  by 
    \begin{equation}\label{eqn:cut-off-op}
        (P_{\mathcal{Q}_b}v)(x',t)=\min(\max(v(x',t),-b),b),\quad \text{for a.e. }(x',t)\in\omega\times(0,T),\quad v\in C([0,T];L^2(\omega)).
    \end{equation}
Our result for   \textbf{(ICP2)} can be stated as follows.
\begin{thm}\label{t3} Let $b>0$ and let Assumption \ref{assum:alg} hold. Then there exist two maps $\mathcal{S}$ and $\mathcal{R}$, taking values in $C([0,T];L^2(\omega))$ and depending only on $\Omega$, $T$, $g$, $u_0$, and $b$, such that problem \textbf{(ICP2)} for $q^\dagger \in \mathcal{Q}_b \cap C([0,T];L^2(\omega))$ can be reformulated as solving the nonlinear equation
\begin{equation}\label{t3a}
    q^\dagger = P_{\mathcal{Q}_b}\big(\mathcal{S}(m) + \mathcal{R}(G_m, q^\dagger)\big)
\end{equation}
within $\mathcal{Q}_b \cap C([0,T];L^2(\omega))$. In addition, equation \eqref{t3a} is globally well-posed in $\mathcal{Q}_b \cap C([0,T];L^2(\omega))$.

\end{thm}
The reconstruction procedure of Theorem \ref{t3} follows the same general philosophy as that of Theorem \ref{t2}:  \textbf{(ICP2)} is first reformulated as a nonlinear equation, and the reconstruction is then obtained through a well-posedness analysis of the equation. However, in contrast to Theorem \ref{t2}, Theorem \ref{t3} provides a reconstruction procedure based on the global well-posedness of  equation \eqref{t3a}. This distinction relies strongly on the formulation of  \textbf{(ICP2)} and on the \textit{a priori} boundedness condition imposed on the class of admissible coefficients. Indeed, unlike \textbf{(ICP1)}, where the objective is to recover the difference between two coefficients from measurements associated with two systems, the goal in \textbf{(ICP2)} is to determine a coefficient from flux measurements corresponding to a single system. Furthermore, the result of Theorem \ref{t3} applies only to the reconstruction of coefficients in $\mathcal{Q}_b$ and satisfying the uniform $L^\infty$ bound $b$. It remains unclear whether either of these restrictions can be removed in the statement of Theorem \ref{t3}.

Despite these limitations compared with the more general reconstruction procedure in Theorem \ref{t2}, the methodology developed for Theorem \ref{t3} appears to be considerably better suited to numerical reconstruction. Indeed, the solution of the nonlinear equation \eqref{t3a} is obtained by means of a contraction mapping argument. This argument is then exploited in Section \ref{sec:error} for the numerical analysis of \textbf{(ICP2)} using the Galerkin finite element method (FEM) together with the backward Euler scheme to discretize the nonlinear operator $q^\dagger\mapsto P_{\mathcal{Q}_b}(\mathcal S (m)+\mathcal R(q^\dagger))$  in \eqref{t3a}. A key feature is that the contraction property is preserved in the fully discrete setting, which guarantees linear convergence of the resulting iterative sequence independently of the discretization parameters. Moreover, we establish an \textit{a priori} error estimate for  \textbf{(ICP2)}. For  the exact coefficient $q^\dagger$ and  the reconstructed coefficient $q^*$, we prove in Theorem \ref{thm:error} that the following estimate holds
\begin{equation*}
\| q^* -q^\dagger  \|_{ \ell^\infty(L^2(\omega))}\le c\big((h^{-1} +\tau^{-1})\delta+|\log h|h^{\frac{1}{2}} + \tau\big),
\end{equation*}
where $\|\cdot\|_{ \ell^\infty(L^2(\omega))}$ denotes the discrete $\ell^\infty$ norm for $L^2(\omega)$-valued sequences and $h$, $\tau$, $\delta$ are the mesh size, time-step size, and noise level, respectively. The four terms in this bound reflect three distinct sources of error: spatial discretization ($| \log h | h^{\frac{1}{2}}$), temporal discretization ($\tau$), and measurement noise ($h^{-1}\delta + \tau^{-1}\delta$). The noise amplification factors $h^{-1}$ and $\tau^{-1}$ reflect the ill-posedness of the inverse problem: refining the discretization improves the approximation of the forward map but simultaneously amplifies the effect of data noise. Balancing these competing effects yields practical guidelines for parameter selection:  given a noise level $ \delta$, the mesh size $h$ and time step $\tau$ should be chosen according to the scaling relations $h \sim \delta^{\frac{2}{3}}$ and $\tau \sim \delta^{\frac{1}{2}}$ to balance the four error terms, yielding an overall convergence rate of $O(\delta^{\frac{1}{3}})$ (up to logarithmic factors) with respect to the noise level $\delta$. This rate is consistent with the stability estimate in Theorem \ref{t3} and provides, to the best of our knowledge, the first rigorous \textit{a priori} error analysis for a nonlinear inverse coefficient problem of this type.

Similar to Theorem \ref{t2}, the regularity conditions on the measurement $m$ imposed in Assumption \ref{assum:alg} can be verified under suitable regularity assumptions of the same type as that in Theorem \ref{t4}.

\section{The direct problem}\label{sec:prelim}
In this section we consider the IBVP \eqref{eqn:forward}
when $q\in C([0,T];H^{s_1}(\omega))$, with $0<s_1\leq \frac{1}{2}$. Using classical lifting argument, this is equivalent to the study of the IBVP with the homogeneous boundary condition
\begin{equation}\label{eq1}
    \left\{\begin{aligned}\partial_t u-\Delta u+q(t,x')u=F(t,x),\quad &(t,x',x_d)=(t,x)\in Q:=
    (0,T)\times\Omega,\\  
    u(t,x)=0,\quad &(t,x)\in(0,T)\times \partial\omega\times(0,\ell),\\ 
    \partial_{x_d}u(t,x', 0)=u(t,x', \ell)=0\quad &(t,x')\in(0,T)\times \omega,\\
    u(0,x)=0,\quad &x\in \Omega.\end{aligned}\right.
\end{equation}
The main result of this section can be stated as follows.
\begin{thm}\label{t1} 
    Let  $q\in C([0,T];H^{s_1}(\omega))$ and $F\in C([0,T];H^{s_0}(\Omega))$, with $0<s_0<s_1\leq \frac{1}{2}$. Then, problem \eqref{eq1} admits a unique solution $u\in C([0,T];H^2(\Omega))\cap C^1([0,T];L^2(\Omega))$. Moreover, if there exists $M>0$ such that
    \begin{equation}\label{t1aa}
        \norm{q}_{L^\infty(0,T;H^{s_1}(\omega))}\leq M,
    \end{equation}
   then there exists $C>0$ depending on $\Omega$, $s_0$, $s_1$ and $T$, and $c>0$ depending only on $\Omega$, $s_0$, $s_1$ such that
    \begin{equation}\label{est1}
        \norm{u}_{C([0,T];H^2(\Omega))}+ \norm{u}_{C^1([0,T];L^2(\Omega))}\leq C\norm{F}_{C([0,T];H^{s_0}(\Omega))}(1+M^2)(1+TM^2)e^{cM^2T}.
    \end{equation} 
\end{thm}

To prove Theorem \ref{t1}, we need the following result about the product of functions in Sobolev space on cylindrical domain of the form $\Omega=\omega\times(0,\ell)$.

\begin{lem}\label{l1} 
    Assume $s_i> s_0>0,\,i=1,2$, and $s_1+s_2-s_0>\frac{d-1}{2}$. Then,  for all $g\in H^{s_1}(\omega)$ and $f\in H^{s_2}(\Omega)$, the map $G:(x',x_d)\mapsto g(x')f(x',x_d)$ lies in $H^{s_0}(\Omega)$ and 
    \begin{equation}\label{l1a}
        \norm{G}_{H^{s_0}(\Omega)}\leq C\norm{g}_{H^{s_1}(\omega)}\norm{f}_{H^{s_2}(\Omega)},
    \end{equation}
    with $C>0$ depending only on $\Omega$ and $s_j$, $j=0,1,2$ .
\end{lem} 
\begin{proof}
We consider only the case $d=3$, and we can deduce the case $d=2$ similarly. Following \cite[Theorem 9.1, Chapter 1]{LM1}, by using the extension map, we may assume without loss of generality that $\omega=\R^2_{x'}$ and $\Omega=\R^2_{x'}\times \R_{x_3}$ with $g\in H^{s_1}(\R^2_{x'})$ and $f\in H^{s_2}(\R^2_{x'}\times \R_{x_3})$ compactly supported. We will use a density argument and  decompose the proof into two steps.\\
\textbf{Step 1} We prove that there exists a constant $C>0$ such that
\begin{equation}\label{l1b}\norm{(x',x_3)\mapsto \phi(x')\psi(x',x_3)}_{H^{s_0}(\R^3)}\leq C\norm{\phi}_{H^{s_1}(\R^2_{x'})}\norm{\psi}_{H^{s_2}(\R^3)},\quad \phi\in C^\infty_0(\R^2_{x'}),\ \psi\in C^\infty_0(\R^3).\end{equation}
Fix $\phi\in C^\infty_0(\R^2_{x'})$, $\psi\in C^\infty_0(\R^3)$ and  denote by $C>0$ a generic constant independent of $\phi$ and $\psi$. Define $\Phi\in C^\infty_0(\R^3)$ by
$$\Phi(x',x_3)=\phi(x')\psi(x',x_3),\quad (x',x_3)\in \R^2_{x'}\times \R_{x_3}.$$
Note that
$$\norm{\Phi}_{H^{s_0}(\R^3)}^2=\int_{\R^3}(1+|\xi|^2)^{s_0}|\mathcal F\Phi(\xi)|^2\d\xi,$$
where, for all $h\in L^1(\R^3)$, $\mathcal F h$ denotes the Fourier transform of $h$ defined by 
$$\mathcal F h(\xi)=(2\pi)^{-\frac 3 2}\int_{\R^3}h(x)e^{-ix\cdot\xi}\d x.$$
Recall that
$$(1+|\xi|^2)^{s_0}\leq ((1+|\xi'|^2)+(1+|\xi_3|^2))^{s_0}\leq 2^{s_0}((1+|\xi'|^2)^{s_0}+(1+|\xi_3|^2)^{s_0}),\quad \xi=(\xi',\xi_3)\in\R^3.$$
Thus, by Fourier-Plancherel formula, we find
$$\norm{\Phi}_{H^{s_0}(\R^3)}^2\leq 2^{s_0}\left(\int_{\R}\int_{\R^2}(1+|\xi'|^2)^{s_0}|\mathcal F_{x'}\Phi(\xi',x_3)|^2\d\xi'\d x_3+\int_{\R}\int_{\R^2}(1+|\xi_3|^2)^{s_0}|\mathcal F_{x_3}\Phi(x',\xi_3)|^2\d x'\d\xi_3\right),$$
where, for all $h\in L^1(\R^3)$, $\mathcal F_{x'} h$ and $\mathcal F_{x_3} h$,  denote  the partial Fourier transforms of $h$ defined by 
$$\mathcal F_{x'} h(\xi',x_3)=(2\pi)^{-1}\int_{\R^2}h(x',x_3)e^{-ix'\cdot\xi'}\d x',\quad \mathcal F_{x_3} h(x',\xi_3)=(2\pi)^{-\frac 1 2}\int_{\R}h(x',x_3)e^{-ix_3\xi_3}\d x_3.$$ 
Therefore, we get
\begin{equation}\label{l1c}
    \norm{\Phi}_{H^{s_0}(\R^3)}\leq C(\norm{\Phi}_{L^2(\R_{x_3};H^{s_0}(\R^2_{x'}))}+\norm{\Phi}_{H^{s_0}(\R_{x_3};L^2(\R^2_{x'}))}).
\end{equation}
Let $x_3\in\R$ and consider the map $\Phi(\cdot,x_3):=x'\mapsto \Phi(x',x_3)$ defined on $\R^2_{x'}$. We have
$$\norm{\Phi(\cdot,x_3)}_{H^{s_0}(\R^2_{x'})}=\norm{\phi\psi(\cdot,x_3)}_{H^{s_0}(\R^2_{x'})}.$$
By the properties of products of functions in Sobolev spaces  \cite[Theorem 5.1]{BH}, we deduce 
$$\norm{\Phi(\cdot,x_3)}_{H^{s_0}(\R^2_{x'})}\leq C\norm{\phi}_{H^{s_1}(\R^2_{x'})}\norm{\psi(\cdot,x_3)}_{H^{s_2}(\R^2_{x'})},$$
with $C>0$ independent of $x_3$. Thus, we get
$$\begin{aligned}
    \norm{\Phi}_{L^2(\R_{x_3};H^{s_0}(\R^2_{x'}))}^2&=\int_\R \norm{\Phi(\cdot,x_3)}_{H^{s_0}(\R^2_{x'})}^2 \d x_3\\
    &\leq C\norm{\phi}_{H^{s_1}(\R^2_{x'})}^2\int_\R\norm{\psi(\cdot,x_3)}_{H^{s_2}(\R^2_{x'})}^2\d x_3\\
    &\leq C\norm{\phi}_{H^{s_1}(\R^2_{x'})}^2\norm{\psi}_{H^{s_2}(\R^3)}^2,
\end{aligned}$$
and it follows that 
\begin{equation}\label{l1d}
    \norm{\Phi}_{L^2(\R_{x_3};H^{s_0}(\R^2_{x'}))}\leq C\norm{\phi}_{H^{s_1}(\R^2_{x'})}\norm{\psi}_{H^{s_2}(\R^3)}.
\end{equation}
Similarly, for all $\xi_3\in\R$, we denote by $\mathcal F_{x_3}\psi(\cdot,\xi_3)$ the map $\R^2\ni x'\mapsto \mathcal F_{x_3}\psi(x',\xi_3)$. Then, we have
$$\norm{\Phi}_{H^{s_0}(\R_{x_3};L^2(\R^2_{x'}))}^2=\int_\R \int_{\R^2}(1+|\xi_3|^2)^{s_0}|\phi(x')|^2|\mathcal F_{x_3}\psi(x',\xi_3)|^2\d x'\d\xi_3.$$
By taking $p>0$ such that $\frac{1}{2}-\frac{s_1}{d-1}<\frac{1}{p}<\frac{s_2}{d-1}-\frac{s_0}{d-1}$,  H\"older's inequality with $\frac{1}{p}+\frac{1}{p'}=\frac{1}{2}$ and the  Sobolev embeddings $H^{s_1}(\R^2_{x'})\hookrightarrow L^p(\R^2_{x'})$ and $H^{s_2-s_0}(\R^2_{x'})\hookrightarrow L^{p'}(\R^2_{x'})$, we get
$$\begin{aligned}
    \norm{\Phi}_{H^{s_0}(\R_{x_3};L^2(\R^2_{x'}))}^2&=\int_\R (1+|\xi_3|^2)^{s_0}\norm{\phi\mathcal F_{x_3}\psi(\cdot,\xi_3)}_{L^2
    (\R^2_{x'})}^2\d\xi_3\\
    &\leq \int_\R (1+|\xi_3|^2)^{s_0}\norm{\phi}_{L^p(\R^2_{x'})}^2\norm{\mathcal F_{x_3}\psi(\cdot,\xi_3)}_{L^{p'}(\R^2_{x'})}^2\d\xi_3\\
    &\leq C\norm{\phi}_{H^{s_1}(\R^2_{x'})}^2\int_\R (1+|\xi_3|^2)^{s_0}\norm{\mathcal F_{x_3}\psi(\cdot,\xi_3)}_{H^{s_2-s_0}(\R^2_{x'})}^2\d \xi_3\\
    &\leq C\norm{\phi}_{H^{s_1}(\R^2_{x'})}^2\int_\R\int_{\R^2} \underbrace{(1+|\xi_3|^2)^{s_0}}_{\leq (1+|\xi|^2)^{s_0}}\underbrace{(1+|\xi'|^2)^{s_2-s_0}}_{\leq (1+|\xi|^2)^{s_2-s_0}}|\mathcal F \psi(\xi',\xi_3)|^2\d\xi' \d\xi_3\\
    &\leq C\norm{\phi}_{H^{s_1}(\R^2_{x'})}^2\int_{\R^3}(1+|\xi|^2)^{s_2}|\mathcal F \psi(\xi)|^2\d\xi\leq C\norm{\phi}_{H^{s_1}(\R^2_{x'})}^2\norm{\psi}_{H^{s_2}(\R^3)}^2.
\end{aligned}$$
By combining this with \eqref{l1c}-\eqref{l1d}, we deduce \eqref{l1b}.\\
\textbf{Step 2} Choose two sequences $(\phi_n)_{n\in\mathbb N}$ and $(\psi_n)_{n\in\mathbb N}$; respectively in $C^\infty_0(\R^2_{x'})$ and $C^\infty_0(\R^3)$ such that
\begin{equation}\label{l1e}\lim_{n\to+\infty}\norm{\phi_n-g}_{H^{s_1}(\R^2_{x'})}=0,\quad \lim_{n\to+\infty}\norm{\psi_n-f}_{H^{s_2}(\R^3)}=0.\end{equation}
Fix also $(\Phi_n)_{n\in\mathbb N}$, $(H_n)_{n\in\mathbb N}$ defined by
$$\Phi_n(x',x_3)=\phi_n(x')\psi_n(x',x_3),\quad H_n(x',x_3)=g(x')\psi_n(x',x_3),\quad n\in\mathbb N,\ (x',x_3)\in\R^2_{x'}\times \R_{x_3}.$$
By H\"older's inequality with  $\frac{1}{p}+\frac{1}{p'}=\frac{1}{2}$ and $ \frac{1}{2}-\frac{s_1}{d-1}<\frac{1}{p}<\frac{s_2}{d-1}$ and Sobolev embedding theorem, we obtain
$$\begin{aligned}\norm{G-H_n}_{L^2(\R^3)}^2&=\int_{\R}\int_{\R^2}|g(x')|^2|\psi_n(x',x_3)-f(x',x_3)|^2\d x'\d x_3\\
&\leq \norm{g}_{L^p(\R^2)}^2\int_{\R} \norm{\psi_n(\cdot,x_3)-f(\cdot,x_3)}_{L^{p'}(\R^2)}^2\d x_3\\
&\leq C\norm{g}_{H^{s_1}(\R^2_{x'})}^2\norm{\psi_n-f}_{H^{s_2}(\R^3)}^2,\end{aligned}$$
with $C>0$ independent of $n$.
In the same way, we can prove  
$$\norm{\Phi_n-H_n}_{L^2(\R^3)}\leq C\norm{\phi_n-g}_{H^{s_1}(\R^2_{x'})}\norm{\psi_n}_{H^{s_2}(\R^3)},$$
and, it follows that 
$$\begin{aligned}
\norm{\Phi_n-G}_{L^2(\R^3)}&\leq \norm{G-H_n}_{L^2(\R^3)}+\norm{\Phi_n-H_n}_{L^2(\R^3)}\\
&\leq C(\norm{g}_{H^{s_1}(\R^2_{x'})}\norm{\psi_n-f}_{H^{s_2}(\R^3)}+\norm{\phi_n-g}_{H^{s_1}(\R^2_{x'})}\norm{\psi_n}_{H^{s_2}(\R^3)}),\end{aligned}$$
with $C>0$ independent of $n$. By combining this with \eqref{l1e}, we deduce that $(\Phi_n)_{n\in\mathbb N}$ converges to $G$ in the sense of $L^2(\R^3)$.
In view of \eqref{l1b}, for all $m,n\in\mathbb N$, we deduce 
$$\norm{\Phi_m-\Phi_n}_{H^{s_0}(\R^3)}\leq C(\norm{\phi_n}_{H^{s_1}(\R^2_{x'})}\norm{\psi_m-\psi_n}_{H^{s_2}(\R^3)}+\norm{\phi_n-\phi_m}_{H^{s_1}(\R^2_{x'})}\norm{\psi_m}_{H^{s_2}(\R^3)})$$
and \eqref{l1e} implies that $(\Phi_n)_{n\in\mathbb N}$ is a Cauchy and therefore convergent sequence of $H^{s_0}(\R^3)$. By the uniqueness of the limit, this implies that $G\in H^{s_0}(\R^3)$ and
$$\norm{G}_{ H^{s_0}(\mathbb{R}^3) }=\lim_{n\to+\infty}\norm{\Phi_n}_{H^{s_0}(\R^3)}\leq C\lim_{n\to+\infty}\norm{\phi_n}_{H^{s_1}(\R^2_{x'})}\norm{\psi_n}_{H^{s_2}(\R^3)}= C\norm{g}_{H^{s_1}(\mathbb{R}^2_{x'})}\norm{f}_{H^{s_2}(\mathbb{R}^3)}.$$
This completes the proof of the lemma.\end{proof}

Using Lemma \ref{l1}, we can complete the proof of Theorem \ref{t1}.

\textbf{Proof of Theorem \ref{t1}.} Consider only the case $d=3$, the case $d=2$ can be similarly deduced. We decompose the proof into three steps.\\
\textbf{Step 1: Unique existence of variational solutions.}
Let $Q=\Omega\times(0,T)$. Let us introduce the unbounded operator $A$ acting on $L^2(\Omega)$ with domain $D(A)=\{f\in H^1(\Omega):\ \Delta f\in L^2(\Omega),\ \partial_{x_d}f|_{x_3=0}=0,\ f|_{x_3=\ell}=0,\ f|_{\partial\omega\times(0,\ell)}=0\}$ defined by $A v=-\Delta v$, $v\in D(A)$. Then $A$ is the unique self-adjoint operator associated with the bilinear form $\mathfrak b$ with domain $D(\mathfrak{b}):=\{f\in H^1(\Omega):\ f|_{x_3=\ell}=0,\ f|_{\partial\omega\times(0,\ell)}=0\}$ defined by
$$\mathfrak {b}(v,w)=\int_{\Omega}\nabla v\cdot\nabla w\d x,\quad v,w\in D(\mathfrak{b}).$$
Following \cite[Theorem 9.1, Chapter 2]{LM1}, we can extend $A$ by density to a bounded linear operator from $D(\mathfrak{b})$ to $D(\mathfrak{b})'$, where $D(\mathfrak{b})'$ denotes the dual space of $D(\mathfrak{b})$, defined by
$$\left\langle Au,v \right\rangle_{D(\mathfrak{b})',D(\mathfrak{b})}=\mathfrak{b}(u,v),\quad u,v\in D(\mathfrak{b}).$$
By Lemma \ref{l1} and the Cauchy-Schwarz inequality, for all $t\in[0,T]$, we obtain 
\begin{align*}\abs{\int_{\Omega}q(t,x')|v(x',x_3)|^2\d x'\d x_3}&\leq \norm{q(t,\cdot)v}_{L^2(\Omega)}\norm{v}_{L^2(\Omega)}\\
&\leq C\norm{q}_{L^\infty(0,T;H^{s_1}(\omega))}\norm{v}_{H^1(\Omega)}\norm{v}_{L^2(\Omega)},\quad v\in D(\mathfrak{b}).
\end{align*}
It follows that
\begin{align*}
&\int_\Omega|\nabla v|^2\d x+\int_{\Omega}q(t,x')|v(x',x_3)|^2\d x'\d x_3\\
\geq& \frac{1}{2}\norm{v}_{H^1(\Omega)}^2-(1+(C\norm{q}_{L^\infty(0,T;H^{s_1}(\omega))})^2)\norm{v}_{L^2(\Omega)}^2,\quad v\in D(\mathfrak{b}).
\end{align*}
Then, using the bilinear form $a(t,\cdot,\cdot)$,  $t\in(0,T)$ defined on $D(\mathfrak{b})$ by
$$a(t,v,w)=\int_\Omega \nabla v\cdot\nabla wdx+\int_{\Omega}q(t,x')v(x',x_3)w(x',x_3)\d x'\d x_3,\quad v,w\in D(\mathfrak{b}),$$
we have
$$a(t,v,w)\leq (1+C\norm{q}_{L^\infty(0,T;H^{s_1}(\omega))})\norm{v}_{H^1(\Omega)}\norm{w}_{H^1(\Omega)},\quad t\in[0,T],\ v,w\in D(\mathfrak{b}),$$ 
$$a(t,v,v)+\left(1+(C\norm{q}_{L^\infty(0,T;H^{s_1}(\omega))})^2\right)\norm{v}_{L^2(\Omega)}^2\geq \frac{1}{2}\norm{v}_{H^1(\Omega)}^2,\quad t\in[0,T],\ v\in D(\mathfrak{b}).$$
Thus, by \cite[Theorem 4.1, Chapter 3]{LM1}, problem \eqref{eq1} admits a unique variational solution $u\in L^2(0,T;D(\mathfrak{b}))\cap H^1(0,T;D(\mathfrak{b})')$. Moreover,  repeating the preceding argument gives  for all $t\in[0,T]$,
\begin{align*}&\int_0^t\int_\Omega F(s,x)u(s,x)\d x\d s\\
&=\int_0^t\left\langle \partial_su(s,\cdot)+A u(s,\cdot),u(s,\cdot)\right\rangle_{D(\mathfrak{b})',D(\mathfrak{b})}\d s+\int_0^t\int_\Omega q(s,x)|u(s,x)|^2\d x\d s\\
&=\int_0^t\frac12\partial_s\left(\int_\Omega|u(s,x)|^2\d x\right)\d s+\int_0^t\int_\Omega(|\nabla u(s,x)|^2+ q(s,x)|u(s,x)|^2)\d x\d s\\
&\geq \frac{1}{2}\int_\Omega |u(t,x)|^2\d x+\int_0^t\left(\int_\Omega|\nabla u(s,x)|^2\d x - C\norm{q}_{L^\infty(0,T;H^{s_1}(\omega))}\norm{u(s,\cdot)}_{H^1(\Omega)}\norm{u(s,\cdot)}_{L^2(\Omega)}\right)\d s\\
&\geq \frac{1}{2}\int_\Omega |u(t,x)|^2\d x+\int_0^t\norm{u(s,\cdot)}_{H^1(\Omega)}^2\d s- (1+C^2\norm{q}_{L^\infty(0,T;H^{s_1}(\omega))}^2)\int_0^t\int_\Omega |u(s,x)|^2\d x\d s,
\end{align*}
where $ C>0$ depends only on $\Omega$. 
Therefore,  under \eqref{t1aa}, we get
$$\int_\Omega |u(t,x)|^2\d x \leq (3+2C^2M^2)\int_0^t\int_\Omega |u(s,x)|^2\d x\d s+\norm{F}_{L^2(Q)}^2,$$
and, by Gronwall's inequality, we get   
$$\int_\Omega |u(t,x)|^2\d x\leq \norm{F}_{L^2(Q)}^2 e^{(3+2C^2M^2)t},\quad t\in[0,T].$$
From this last estimate, we obtain 
\begin{equation}
    \label{e1}\norm{u}_{L^\infty(0,T;L^2(\Omega))}\leq \norm{F}_{L^2(Q)}e^{(\frac32+C^2M^2)T}.\end{equation}
Similarly, we obtain
$$\int_0^T\norm{u(s,\cdot)}_{H^1(\Omega)}^2\d s\leq (2+C^2M^2)\int_0^T\int_\Omega |u(s,x)|^2\d x\d s+\norm{F}_{L^2(Q)}^2,$$
and, by the estimate \eqref{e1}, we get
$$\norm{u}_{L^2(0,T;H^1(\Omega))}^2\leq \norm{F}_{L^2(Q)}^2(T(2+C^2M^2)e^{(3+2C^2M^2)T}+1),$$
and, there exist $C,c>0$ depending on $\Omega$, $s_0$, $s_1$ and $T$  such that
\begin{equation}\label{est2}\norm{u}_{L^2(0,T;H^1(\Omega))}+ \norm{u}_{H^1(0,T;H^{-1}(\Omega))}\leq C(1+TM^2)e^{cM^2T}\norm{F}_{L^2(Q)}.\end{equation}
\ \\
 \textbf{Step 2: Unique existence of a solution in $L^2(0,T;H^2(\Omega))\cap H^1(0,T;L^2(\Omega))$.}
Fix $G$  to be
$$G(t,x',x_3)=-q(t,x')u(t,x',x_3)+F(t,x',x_3),\quad (t,x',x_3)\in(0,T)\times\Omega.$$
Then we can rewrite problem \eqref{eq1} in the following abstract way
\begin{equation}\label{eqa}
\left\{\begin{aligned} u'(t)+Au(t)=G(t),\quad &t\in(0,T),\\  
u(0)=0,\end{aligned}\right.
\end{equation}
where, for a.e. $t\in(0,T)$,  $u(t)$ and $G(t)$ denote respectively the maps $x\mapsto u(t,x)$, $x\mapsto G(t,x)$.
Repeating the  argument of Step 1 yields
$$\norm{G}_{L^2(0,T;L^2(\Omega))}\leq C\norm{q}_{L^\infty(0,T;H^{s_1}(\omega))}\norm{u}_{L^2(0,T;H^1(\Omega))}+\norm{F}_{L^2(0,T;L^2(\Omega))}<\infty.$$
Thus, the map $G\in L^2(0,T;L^2(\Omega))$ and problem \eqref{eqa} admits a unique solution $u\in H^1(0,T;L^2(\Omega))\cap L^2(0,T;D(A))$. Following a separation of variable argument similar to \cite[Lemma 2.4]{GK}, one can check that $D(A)$ embeds continuously into $H^2(\Omega)$ which implies that the unique solution $u$ of \eqref{eq1} lies in $ L^2(0,T;H^2(\Omega))\cap H^1(0,T;L^2(\Omega))$. In addition, there exists $C>0$ depending on $\Omega$ and $T$ such that
\begin{align*} &\norm{u}_{L^2(0,T;H^2(\Omega))}+\norm{u}_{H^1(0,T;L^2(\Omega))}\leq C\norm{G}_{L^2(Q)}\\
\leq& C\left(\norm{q}_{L^\infty(0,T;H^{s_1}(\omega))}\norm{u}_{L^2(0,T;H^1(\Omega))}+\norm{F}_{L^2(0,T;L^2(\Omega))}\right).
\end{align*}
Therefore, under \eqref{t1aa}, by the estimate \eqref{est2}, we can find $C,c>0$ depending  $\Omega$, $T$, $s_0$ and $_1$ such that 
\begin{equation}\label{est3}\norm{u}_{L^2(0,T;H^2(\Omega))}+\norm{u}_{H^1(0,T;L^2(\Omega))}\leq C(1+M)(1+TM^2)e^{cM^2T}\norm{F}_{L^2(Q)}.\end{equation}
\ \\
\textbf{Step 3: Unique existence of solution in $u\in C([0,T];H^2(\Omega))\cap C^1([0,T];L^2(\Omega))$.}
Since $A$ is a self-adjoint operator with a spectrum consisting of a non-decreasing sequence of non-negative eigenvalues $(\lambda_n)_{n\geq1}$ . Let us also introduce 
an orthonormal  basis in the Hilbert space $L^2(\Omega)$ of eigenfunctions $(\phi_n)_{n\geq1}$ of $A$ associated to the non-decreasing sequence of  eigenvalues $(\lambda_n)_{n\geq1}$.  For all $s\geq 0$, we denote by $A^s$ the operator defined by 
\[A^s g=\sum_{n=1}^{+\infty}\left\langle g,\phi_n\right\rangle \lambda_{n}^s\phi_{n},\quad g\in D(A^s)=\left\{g\in L^2(\Omega):\ \sum_{n=1}^{+\infty}\abs{\left\langle g,\phi_{n}\right\rangle}^2 \lambda_{n}^{2s}<\infty\right\}\]
and consider on $D(A^s)$ the norm
\[\|g\|_{D(A^s)}=\left(\sum_{n=1}^{+\infty}\abs{\left\langle g,\phi_{n}\right\rangle}^2 \lambda_{n}^{2s}\right)^{\frac{1}{2}},\quad g\in D(A^s).\]
Consider also the operator valued function $S(t)$  defined by
$$S(t)h=\sum_{n=1}^\infty e^{-\lambda_nt}\left\langle h,\phi_n\right\rangle_{L^2(\Omega)}\phi_n,\quad h\in L^2(\Omega).$$
The solution of problem \eqref{eq1} solves the following integral equation
$$u(t)=-\int_0^tS(t-s)q(s)u(s)\d s+\int_0^tS(t-s)F(s)\d s.$$
Since $D(A)$ embeds continuously into $H^2(\Omega)$, by interpolation, for all $s\in[0,1]$, $D(A^s)$ embeds continuously into $H^{2s}(\Omega)$. Moreover, recalling that $D(A^{\frac{1}{2}})=D(\mathfrak{b})=\{f\in H^1(\Omega):\ f|_{x_3=\ell}=f|_{\partial\omega\times(0,\ell)}=0\}$, $0<s_0<\frac{1}{2}$ and by \cite[Theorem 11.4, Chapter 1]{LM1} and \cite[Theorem 11.6, Chapter 1]{LM1}, we deduce that $D(A^{\frac{s_0}{2}})=H^{s_0}(\Omega)$ with equivalent norms. By \cite[Theorem 3.1, Chapter 1]{LM1}, we get $u\in H^1(0,T;L^2(\Omega))\cap L^2(0,T;D(A))\subset C([0,T];D(A^{\frac{1}{2}}))\subset C([0,T];H^1(\Omega))$. Then, by Lemma \ref{l1}, we find $(t,x',x_3)\mapsto q(t,x')u(t,x',x_3)\in C([0,T];H^{s_0}(\Omega))\subset C([0,T];D(A^{\frac{s_0}{2}}))$. Therefore, we obtain
$$\norm{S(t-s)q(s)u(s)}_{D(A)}=\norm{AS(t-s)q(s)u(s)}_{L^2(\Omega)}=\norm{A^{1-\frac{s_0}{2}}S(t-s)A^{\frac{s_0}{2}}q(s)u(s)}_{L^2(\Omega)},$$
for $0<s<t\leq T$.
Moreover, for all $0<s<t\leq T$, we have
$$\begin{aligned}\norm{A^{1-\frac{s_0}{2}}S(t-s)h}_{L^2(\Omega)}^2&=\sum_{n=1}^\infty \lambda_n^{2-s_0}e^{-2\lambda_n(t-s)}|\left\langle h,\phi_n\right\rangle_{L^2(\Omega)}|^2\\
&\leq \left(\sup_{\tau>0}\tau^{2-s_0}e^{-2\tau}\right)(t-s)^{s_0-2}\norm{h}_{L^2(\Omega)}^2.\end{aligned}$$
It follows that
$$\norm{S(t-s)q(s)u(s)}_{D(A)}\leq C(t-s)^{\frac{s_0}{2}-1}\norm{qu}_{L^\infty(0,T;D(A^{\frac{s_0}{2}}))},$$
and by the dominated convergence theorem, we deduce  
$$t\mapsto -\int_0^tS(t-s)q(s)u(s)\d s\in C([0,T];D(A))\subset C([0,T];H^2(\Omega)).$$
In a similar way, we can prove   
$$t\mapsto \int_0^tS(t-s)F(s)\d s\in C([0,T];D(A))\subset C([0,T];H^2(\Omega)),$$
which implies that $u\in C([0,T];H^2(\Omega))$. Repeating similar arguments gives  $u\in C^1([0,T]; L^2(\Omega))$. Moreover, from the above argumentation and the estimate \eqref{l1a}, we deduce that there exists $C>0$ depending on $\Omega$, $s_0$, $s_1$ and $T$ such that
\begin{align}\label{est4}
&\norm{u}_{C([0,T];H^2(\Omega))}+ \norm{u}_{C^1([0,T];L^2(\Omega))}\\
\leq& C\left(\norm{q}_{L^\infty(0,T;H^{s_1}(\omega))}\norm{u}_{C([0,T];H^1(\Omega))}+\norm{F}_{C([0,T];H^{s_0}(\Omega))}\right).\nonumber\end{align}
By \cite[Theorem 11.6, Chapter 1]{LM1}, there exists $C=C(\Omega, s_0, s_1 ,T)$ such that
\begin{equation}\label{est5}\norm{u}_{C([0,T];H^1(\Omega))}\leq C(\norm{u}_{L^2(0,T;H^2(\Omega))}+\norm{u}_{H^1(0,T;L^2(\Omega))}).\end{equation}
Thus, under condition \eqref{t1aa},  by combining \eqref{est3}-\eqref{est5}, we obtain \eqref{est1}. \qed

Now consider the following problem  connected with \eqref{eq1}
\begin{equation}\label{eq2}
\left\{\begin{aligned}\partial_t v-\Delta v+q(t,x')v=F(t,x),\quad &(t,x',x_3)=(t,x)\in 
Q=(0,T)\times\Omega,\\  
v(t,x)=0,\quad &(t,x)\in(0,T)\times \partial\omega\times(0,\ell),\\ 
v(t,x', 0)=v(t,x', \ell)=0\quad &(t,x')\in(0,T)\times \omega,\\
v(0,x)=0,\quad &x\in \Omega.\end{aligned}\right.
\end{equation}
In a similar way to Theorem \ref{t1} we can prove the following.

\begin{thm}\label{t1b} 
    Let  $q\in C([0,T];H^{s_1}(\omega))$ and $F\in C([0,T];H^{s_0}(\Omega))$, with $0<s_0<s_1\leq\frac{1}{2}$. Then, problem \eqref{eq2} admits a unique solution $v\in C([0,T];H^2(\Omega))\cap C^1([0,T];L^2(\Omega))$. Moreover, the map $\mathcal K:F\mapsto v$ is a continuous linear map from $C([0,T];H^{s_0}(\Omega))$ to $C([0,T];H^2(\Omega))\cap C^1([0,T];L^2(\Omega))$.
\end{thm}

Similarly to Theorem \ref{t1}, we can show the following corollary with bounded coefficients. 
\begin{cor}\label{cor:forward_weak_reg}
    For $q\in L^{\infty}(0,T;L^\infty(\omega))$ and $F\in C([0,T];L^2(\Omega))$  problems \eqref{eq1}, \eqref{eq2} admit respectively a unique strong solution $u,\,v\in L^2(0,T;H^2(\Omega))\cap H^1(0,T;L^2(\Omega))$ satisfying 
    $$u,v\in\underset{s\in(0,1)}{\bigcap}C([0,T];H^{2s}(\Omega)).$$
 Moreover, if there exists $M>0$ such that
 $\norm{q}_{L^{\infty}(0,T;L^\infty(\omega))}\leq M$, then for any $\beta\in[0,1)$, there exist $C,c>0$ depending only on $\beta$, $\Omega$ and $T$   such that
    \begin{equation*}  
        \norm{u}_{C([0,T];H^{2\beta}(\Omega))} \leq C\norm{F}_{C([0,T];L^2(\Omega))}(1+M^4)e^{cM^2} .
    \end{equation*} 
\end{cor}

\section{Reconstruction from local well-posedness of the inverse problem}\label{sec:recon_local}
In this section, we prove Theorem \ref{t2}. This proof is lengthy and is divided into four steps.

\subsection{Derivation of  the nonlinear equation equivalent to \textbf{(ICP1)}}
First, we prove that \textbf{(ICP1)} is equivalent to solving a nonlinear equation of the form \eqref{t2a}. To highlight the dependency on the unknown coefficient $q$, we denote by $u_2(q)$ the solution of \eqref{eq3} with $j=2$ (i.e., $u_2(q)=u_2$). 
By Theorem \ref{t1},  for $j=1,2$, problem \eqref{eq3} admits a unique solution $u_j\in C([0,T];H^2(\Omega))\cap C^1([0,T];L^2(\Omega))$. Then,  $u=u_1-u_2$ solves the problem
\begin{equation}\label{eq4}
\left\{\begin{aligned}\partial_t u-\Delta u+q_1(t,x')u=q(t,x')u_2(t,x',x_d),\quad &(t,x',x_d)=(t,x)\in 
Q:=(0,T)\times \Omega,\\  
\partial_{x_d}u(t,x', 0)=u(t,x', \ell)=0\quad &(t,x')\in (0,T)\times \omega,\\
u(t,x)=0,\quad &(t,x)\in(0,T)\times \partial\omega\times(0,\ell),\\ 
u(0,x)=0,\quad &x\in\Omega.\end{aligned}\right.
\end{equation}
Fix $v=\partial_{x_d}(u_1-u_2)$. Then, we have $v=v_1+v_2$ with $v_1$ solving 
\begin{equation}\label{eq5}
\left\{\begin{aligned}\partial_t v_1-\Delta v_1+q_1(t,x')v_1=q(t,x')\partial_{x_d}u_2(t,x',x_d),\quad &(t,x',x_d)=(t,x)\in 
Q,\\  
v_1(t,x)=0,\quad &(t,x)\in(0,T)\times \partial\omega\times(0,\ell),\\ 
v_1(t,x', 0)=v_1(t,x', \ell)=0\quad &(t,x')\in(0,T)\times \omega,\\
v_1(0,x)=0,\quad &x\in \Omega,\end{aligned}\right.
\end{equation}
and $v_2$ solving 
\begin{equation}\label{eq6}
\left\{\begin{aligned}\partial_t v_2-\Delta v_2+q_1(t,x')v_2=0,\quad &(t,x',x_d)=(t,x)\in 
Q,\\  
v_2(t,x)=0,\quad &(t,x)\in(0,T)\times \partial\omega\times(0,\ell),\\ 
v_2(t,x', 0)=0\quad &(t,x')\in(0,T)\times \omega,\\
v_2(t,x', \ell)=m(t,x'),\quad &(t,x')\in (0,T)\times\omega,\\
v_2(0,x)=0\quad &x\in \Omega.
\end{aligned}\right.
\end{equation}
Since $q\in C([0,T];H^{s_1}(\omega))$, $\partial_{x_d}u_2\in C([0,T];H^{1}(\Omega))$ and by Lemma \ref{l1} and Theorem \ref{t1b},   $v_1\in C([0,T];H^2(\Omega))\cap C^1([0,T];L^2(\Omega))$. In the same way, using \eqref{trace2} and Theorem \ref{t1b},  $v_2\in C([0,T];H^2(\Omega))\cap C^1([0,T];L^2(\Omega))$. It follows that $u, \partial_{x_d}u=v\in C^1([0,T];L^2(\Omega))\cap C([0,T];H^2(\Omega))$ which, combined with Lemma \ref{l1}, imply $\partial_tu,\Delta u, q_1u\in C([0,T];H^1(0,\ell;L^2(\omega)))$. Fixing $t\in[0,T]$ and taking the trace of the equation $\partial_tu-\Delta u +q_1u$ at $x_d=\ell$, gives
\begin{align*}
q(t,x')&=(q_2(t,x')-q_1(t,x'))u_2(t,x',\ell)\\
&=\partial_tu(t,x',\ell)-\Delta u(t,x',\ell)+q_1(t,x')u(t,x',\ell)\\
&=-\partial_{x_d}^2u(t,x',\ell)=-\partial_{x_d}v(t,x',\ell),\quad (t,x')\in[0,T]\times \omega.
\end{align*}
Therefore,  $q\in C([0,T];H^{s_1}(\omega))$ is a solution of the nonlinear  equation
\begin{equation}\label{int}q=\mathcal M_{q_1}(m)+\mathcal N_{q_1}(q\partial_{x_d}u_2(q)),\end{equation}
where $\mathcal M_{q_1}(m)(t,x')=-\partial_{x_d}v_2(t,x',\ell)$ with $v_2$ solving problem \eqref{eq6} and $\mathcal N_{q_1}(q\partial_{x_d}u_2(q))(t,x')=-\partial_{x_d}v_1(t,x',\ell)$ with $v_1$ solving  problem \eqref{eq5}.
Thus \textbf{(ICP1)} for $q\in C([0,T];H^{s_1}(\omega))$ is equivalent to solving in $C([0,T];H^{s_1}(\omega))$ the nonlinear equation \eqref{t2a} with $\mathcal M_{q_1}$ as in \eqref{int} and $\mathcal K_{q_1}(q)=\mathcal N_{q_1}(q\partial_{x_d}u_2(q))$, $q\in C([0,T];H^{s_1}(\omega))$.
This completes the proof of the first statement of Theorem \ref{t2}. Below we study   the well-posedness of the nonlinear equation \eqref{int}. We equip the measurement data $m$ with the induced topology of $\mathcal M_{q_1}(m)$ on $C([0,T];H^{s_1}(\omega))$ and  under suitable regularity assumptions, this would imply the continuity of $m$ in suitable Sobolev spaces (see Section \ref{sec:stab} for more details).

We   consider a local well-posedness of \eqref{int} in $C([0,T];H^{s}(\omega))$, $s\in(0,\frac12]$, under the assumption of either smallness of the final time $T$ or smallness of the measurement $m$. 

\subsection{Existence and continuity for small time.}\label{EC}
Next we prove the second claim of Theorem \ref{t2} and the existence, as well as the continuity with respect to $m$, of solutions  of \eqref{int} for   small final time $T>0$. Fix $T_1>0$, $ s_0\in(0, s_1)$ and consider $T\in(0,T_1]$ a parameter to be determined below.
For $r>0$, we define the set $B_r:=\{q\in C([0,T];H^{s_1}(\omega)): \norm{q}_{C([0,T];H^{s_1}(\omega))}\leq r\}$. We   consider the solution $q$ of the nonlinear equation \eqref{int} as a fixed point of the map $\mathcal G$ on $B_{2P}$ defined by
$$\mathcal G(q)=\mathcal M_{q_1}(m)+\mathcal N_{q_1}(q\partial_{x_d}u_2(q)),\quad q\in B_{2P}.$$
Following \cite[Lemma 2.4]{GK}, we define an unbounded operator $L$ acting on $L^2(\Omega)$ with domain $D(L)=H^1_0(\Omega)\cap H^2(\Omega)$ defined by $Lv=-\Delta v$, $v\in D(L)$. Then, $L$ is a self-adjoint operators with a spectrum consisting of a non-decreasing sequence of non-negative eigenvalues $(\mu_n)_{n\geq1}$ . We define 
an orthonormal  basis in the Hilbert space $L^2(\Omega)$ of eigenfunctions $(\psi_n)_{n\geq1}$ of $L$ associated to the non-decreasing sequence of  eigenvalues $(\mu_n)_{n\geq1}$. For all $t>0$, we define the linear map $H(t)$ from $L^2(\Omega)$ to $H^{s_1}(\omega)$ by  
$$H(t)f:=\partial_{x_d}\left(\sum_{n=1}^\infty e^{-\mu_nt}\left\langle f,\psi_n\right\rangle\psi_n\right)|_{x_d=\ell}=\sum_{n=1}^\infty e^{-\mu_nt}\left\langle f,\psi_n\right\rangle\partial_{x_d}\psi_n|_{x_d=\ell},\quad f\in L^2(\Omega).$$
From now on and in all the remaining parts of this article, we fix $\tau_0:=\frac{3+2(s_1-s_0)}{4}\in(0,1)$. In view of \cite[Theorem 11.4, Chapter 1]{LM1}, \cite[Theorem 11.6, Chapter 1]{LM1}, for $s_0\in(0,1/2)$,  $D(L^{\frac{s_0}{2}})=H^{s_0}(\Omega)$ with equivalent norm, and, for all $t>0$, we have 
\begin{align*}&\norm{H(t)f}_{H^{s_1}(\omega)}^2\leq C\norm{\sum_{n=1}^\infty e^{-\mu_nt}\left\langle f,\psi_n\right\rangle\psi_n}_{H^{\frac{3+2s_1}{2}}(\Omega)}^2 \leq C\norm{\sum_{n=1}^\infty e^{-\mu_nt}\left\langle f,\psi_n\right\rangle\psi_n}_{D(L^{\frac{3+2s_1}{4}})}^2\\
\leq& C\sum_{n=1}^\infty \mu_n^{\frac{3+2s_1}{2}}e^{-2t\mu_n}\abs{\left\langle f,\psi_n\right\rangle}^2 \leq C\left( \sup_{\sigma>0} \sigma^{2\tau_0} e^{-2\sigma}\right)t^{-2\tau_0}\sum_{n=1}^\infty \mu_n^{s_0}\abs{\left\langle f,\psi_n\right\rangle}^2\\
\leq& C t^{-2\tau_0}\norm{f}_{D(L^{\frac{s_0}{2}})}^2 \leq   C t^{-2\tau_0}\norm{f}_{H^{s_0}(\Omega)}^2,\quad f\in H^{s_0}(\Omega),
\end{align*}
with $C>0$ a constant that may change from line to line depending only on $\Omega$. Thus, we obtain
\begin{equation}\label{t2c}\norm{H(t)f}_{H^{s_1}(\omega)}\leq Ct^{-\tau_0}\norm{f}_{H^{s_0}(\Omega)}.\end{equation}
Note that the solution $v_1$ of \eqref{eq5} solves also the following integral equation
\begin{align*} v_1(t,\cdot)=&-\int_0^t\left(\sum_{n=1}^\infty e^{-\mu_n(t-s)}\left\langle q_1v_1(s,\cdot),\psi_n\right\rangle\psi_n\right)\d s\\
&+\int_0^t\left(\sum_{n=1}^\infty e^{-\mu_n(t-s)}\left\langle q\partial_{x_d}u_2(s,\cdot),\psi_n\right\rangle\psi_n\right)\d s,\quad t\in(0,T).
\end{align*}
Therefore, we obtain the following representation of the map $\mathcal N_{q_1}(q\partial_{x_d}u_2(q))$
$$\mathcal N_{q_1}(q\partial_{x_d}u_2(q))(t)=\int_0^tH(t-s)q_1v_1(s,\cdot)\d s-\int_0^tH(t-s)q\partial_{x_d}u_2(s,\cdot)\d s,\quad t\in(0,T)$$
and the map $\mathcal G$ takes the form
$$\mathcal G(q)(t)=\mathcal M_{q_1}(m)(t)+\int_0^tH(t-s)q_1v_1(s,\cdot)\d s-\int_0^tH(t-s)q\partial_{x_d}u_2(s,\cdot)\d s,\quad t\in(0,T),\  q\in B_{2P}.$$
In view of \eqref{est1}, there exist $C,c>0$ depending only on $\Omega$, $g$, $u_0$, $q_1$, $s_1$, $s_0$ and $T_1$ such that 
\begin{equation}\label{tetu}\norm{\partial_{x_d}u_2}_{C([0,T];H^1(\Omega))}\leq \norm{u_2}_{C([0,T];H^2(\Omega))}\leq C(1+P^2)(1+TP^2)e^{cP^2T}
,\quad q\in B_{2P}.\end{equation}
By combining this with Lemma \ref{l1}, we obtain
\begin{equation}\label{t2d}\norm{q\partial_{x_d}u_2}_{C([0,T];H^{s_0}(\Omega))}\leq C(1+P^3)(1+TP^2)e^{cP^2T},\quad q\in B_{2P}.\end{equation}
Similarly, we obtain 
$$\norm{q_1v_1}_{C([0,T];H^{s_0}(\Omega))}\leq C\norm{q\partial_{x_d}u_2}_{C([0,T];H^{s_0}(\Omega))}\leq C(1+P^3)(1+TP^2)e^{cP^2T},\quad q\in B_{2P}.$$
Fix $q,\widetilde{q}\in B_{2P}$ and consider $\widetilde{u}_2$ and $\widetilde{v}_1$, respectively the solution of \eqref{eq3} and \eqref{eq5} with $q=\widetilde{q}$ and $\widetilde{u}_2=u_2(\widetilde{q})$. Then, we get
\begin{align*} 
\mathcal G(q)(t)-\mathcal G(\widetilde{q})(t)=&\int_0^tH(t-s)q_1(v_1(s,\cdot)-\tilde{v_1}(s,\cdot))\d s\\
&-\int_0^tH(t-s)(q\partial_{x_d}u_2(q)-\widetilde{q}\partial_{x_d}u_2(\widetilde{q}))(s,\cdot))\d s,\quad t\in[0,T].
\end{align*}
Note that $u_3=u_2(q)-u_2(\widetilde{q})$ solves the problem
$$\left\{\begin{aligned}\partial_t u_3-\Delta u_3+q_1(t,x')u_3+q(t,x')u_3=(\widetilde{q}-q)(t,x')\widetilde{u}_2(t,x),\quad &(t,x)\in 
(0,T)\times \Omega,\\  
u_3(t,x', \ell)=0\quad &(t,x')\in(0,T)\times \omega,\\
u_3(t,x)=0 ,\quad &(t,x)\in(0,T)\times \partial\omega\times(0,\ell),\\ 
\partial_{x_d}u_3(t,x', 0)=0\quad &(t,x')\in(0,T)\times \omega,\\
u_3(0,x)=0,\quad &x\in\Omega.\end{aligned}\right.$$
In light of Theorem \ref{t1} and Lemma \ref{l1},  there exist $C,c>0$  depending only on $\Omega$, $s_1$, $s_0$, $g$, $u_0$, $q_1$ and $T_1$ such that 
$$\begin{aligned}\norm{u_2(q)-\widetilde{u}_2(\widetilde{q})}_{C([0,T];H^2(\Omega))}&\leq C(1+P^2)(1+TP^2)e^{cP^2T}\norm{(\widetilde{q}-q)\widetilde{u}_2}_{C([0,T];H^{s_0}(\Omega))}\\
&\leq  C(1+P^2)^2(1+TP^2)^2e^{2cP^2T}\norm{\tilde{q}-q}_{C([0,T];H^{s_1}(\omega))}.\end{aligned}$$
Using this last estimate and \eqref{tetu}, we obtain
\begin{align}&\norm{q\partial_{x_d}u_2(q)-\widetilde{q}\partial_{x_d}\widetilde{u}_2(\widetilde{q})}_{C([0,T];H^{s_0}(\Omega))}\nonumber\\
\leq& \norm{(q-\tilde{q})\partial_{x_d}u_2(q)}_{C([0,T];H^{s_0}(\Omega))}+\norm{\widetilde{q}(\partial_{x_d}u_2(q)-\partial_{x_d}u_2(\widetilde{q}))}_{C([0,T];H^{s_0}(\Omega))}\nonumber\\
\leq& C(1+P^2)(1+TP^2)e^{cP^2T}\norm{\widetilde{q}-q}_{C([0,T];H^{s_1}(\omega))}+CP\norm{u_2(q)-\widetilde{u}_2(\widetilde{q})}_{C([0,T];H^2(\Omega))}\nonumber\\
\leq& C(1+P^5)(1+TP^2)^2e^{2cP^2T}\norm{\widetilde{q}-q}_{C([0,T];H^{s_1}(\omega))}. \label{tutu} \end{align}
Similarly, we find 
\begin{align*}\norm{q_1(v_1-\widetilde{v}_1)}_{C([0,T];H^{s_0}(\Omega))}&\leq C\norm{v_1-\widetilde{v}_1}_{C([0,T];H^{1}(\Omega))}\leq C\norm{q\partial_{x_d}u_2-\widetilde{q}\partial_{x_d}\widetilde{u}_2}_{C([0,T];H^{s_0}(\Omega))}\\
&\leq C(1+P^5)(1+TP^2)^2e^{2cP^2T}\norm{\widetilde{q}-q}_{C([0,T];H^{s_1}(\omega))}.
\end{align*}
By combining \eqref{t2c}-\eqref{tutu}, we get 
\begin{align*}\norm{\mathcal Gq-\mathcal G\widetilde{q}}_{C([0,T];H^{s_1}(\omega))}&\leq C(1+P^5)(1+TP^2)^2e^{2cP^2T}\left(\int_0^Tt^{-\tau_0}\d t\right) \norm{q-\widetilde{q}}_{C([0,T];H^{s_1}(\omega))}\\
&\leq C(1+P^5)(1+TP^2)^2e^{2cP^2T} T^{1-\tau_0}\norm{q-\widetilde{q}}_{C([0,T];H^{s_1}(\omega))},\ q,\widetilde{q}\in B_{2P},
\end{align*}
where $C>0$ depends only on $\Omega$, $s_1$, $s_0$, $g$, $u_0$, $q_1$ and $T_1$.
Thus, by fixing $$T_0=\min((8C(1+P^5)e^{2c})^{-\frac{1}{1-\tau_0}},P^{-2}),$$
for $T\leq T_0$, we obtain
\begin{align*}
\norm{\mathcal Gq-\mathcal G\widetilde{q}}_{C([0,T];H^{s_1}(\omega))}&\leq \tfrac12 \norm{q-\widetilde{q}}_{C([0,T];H^{s_1}(\omega))},\quad q,\widetilde{q}\in B_{2P},\\
\norm{\mathcal Gq}_{C([0,T];H^{s_1}(\omega))}&\leq \norm{\mathcal M_{q_1}(m)}_{C([0,T];H^{s_1}(\omega))}+\tfrac12\norm{q}_{C([0,T];H^{s_1}(\omega))}\leq2P,\quad  q\in B_{2P}.
\end{align*}
Thus, $\mathcal G$ is a contraction of $B_{2P}$ which admits a unique fixed point $q\in B_{2P}$. This proves the existence of a solution of \eqref{int} in $C([0,T];H^{s_1}(\omega))$. Moreover, the continuity with respect to $m$ of this solution can be easily deduced from the above arguments.

\subsection{Existence and continuity for small data} In this step, we consider the third claim of  Theorem \ref{t2} and the existence, as well as the continuity, of solutions  of \eqref{int} for sufficiently small data $P=\norm{\mathcal M_{q_1}(m)}_{C([0,T];H^{s_1}(\omega))}$. For this purpose,  for $\lambda>0$, we define a weighted norm in $C([0,T])$ by 
\begin{equation*} 
\|v\|_{C_\lambda([0,T])}=
         \max_{t\in [0,T]} | e^{-\lambda t}v(t)|.
\end{equation*}
This norm is equivalent to the standard norm of $C([0,T])$. Similarly, for a Banach space $X$, we denote $\|v\|_{C_\lambda([0,T];X)} $ as the weighted norm for vector-valued space $C_\lambda([0,T];X)$. For $r>0$, we define the set $B_r^\lambda:=\{q\in C([0,T];H^{s_1}(\omega)): \norm{q}_{C_\lambda([0,T];H^{s_1}(\omega))}\leq r\}$. Then $B_r^\lambda \subset B_{e^{\lambda T}r}$. Fix $q,\tilde{q}\in B_{2P}^\lambda$ and let $\widetilde{u_2}$ and $\widetilde{v_1}$, be the solution of \eqref{eq3} and \eqref{eq5} with $q=\widetilde{q}$ and $u_2=\widetilde{u_2}$. In all this step $P$ and $\lambda$ will be two positive constants to be determined and $C>0$ will be a constant independent of $P$ and $\lambda$. For all $t\in[0,T]$, we get
\begin{align*} e^{-\lambda t}\mathcal G(q)(t)-e^{-\lambda t}\mathcal G(\widetilde{q})(t) 
 =&\int_0^te^{-\lambda (t-s)}H(t-s)e^{-\lambda s}q_1(v_1(s,\cdot)-\widetilde{v}_1(s,\cdot))\d s\\
&-\int_0^te^{-\lambda (t-s)}H(t-s)e^{-\lambda s}(q\partial_{x_d}u_2(s,\cdot)-\widetilde{q}\partial_{x_d}\widetilde{u}_2(s,\cdot))\d s.
\end{align*}
Then, similarly to Step 1, since $q,\widetilde{q}\in B_{2e^{\lambda T}P}$, we find
\begin{align} &\norm{\mathcal G(q)-\mathcal G(\widetilde{q})}_{C_\lambda([0,T];H^{s_1}(\omega))}\nonumber\\
&\leq C\left(\norm{q_1(v_1-\widetilde{v}_1)}_{C_\lambda([0,T];H^{s_0}(\Omega)}+\norm{q\partial_{x_d}u_2(q)-\widetilde{q}\partial_{x_d}\widetilde{u}_2(\widetilde{q})}_{C_\lambda([0,T];H^{s_0}(\Omega)}\right)\int_0^T e^{-\lambda t}t^{-\tau_0}\d t\nonumber\\
&\leq C\left(\norm{v_1-\widetilde{v}_1}_{C_\lambda([0,T];H^{1}(\Omega))}+\norm{q\partial_{x_d}u_2(q)-\widetilde{q}\partial_{x_d}\widetilde{u}_2(\widetilde{q})}_{C_\lambda([0,T];H^{s_0}(\Omega))}\right)\int_0^\infty e^{-\lambda t}t^{-\tau_0}\d t\nonumber\\
&\leq C\lambda^{\tau_0-1}\left(\norm{v_1-\widetilde{v}_1}_{C_\lambda([0,T];H^{1}(\Omega))}+\norm{q\partial_{x_d}u_2(q)-\widetilde{q}\partial_{x_d}\widetilde{u}_2(\widetilde{q})}_{C_\lambda([0,T];H^{s_0}(\Omega))}\right).\label{t2aa} \end{align}
In the same way, by fixing $u_3=u_2(q)-u_2(\tilde{q})$, for all $t\in [0,T]$, we obtain the representation
\begin{align*}
e^{-\lambda t}u_3(\cdot,t)
=&-\int_0^te^{-\lambda (t-s)}S(t-s)e^{-\lambda s}q_2u_3(s,\cdot)\d s\\
&-\int_0^te^{-\lambda (t-s)}S(t-s)e^{-\lambda s}(q-\widetilde{q})\widetilde{u}_2(s,\cdot)\d s.
\end{align*}
By combining this with \eqref{t2d} and repeating the arguments of Step 1, for all $t\in [0,T]$, we find
$$\begin{aligned}&\norm{e^{-\lambda t}u_3(\cdot,t)}_{H^2(\Omega)}\\
&\leq C\norm{q_2}_{C([0,T];H^{s_1}(\omega))}\int_0^te^{-\lambda (t-\tau)}(t-\tau)^{-\frac{(2-s_0)}{2}}\norm{e^{-\lambda \tau}u_3(\cdot,\tau)}_{H^{1}(\Omega)}\d\tau\\
&\ \ \ +C\norm{q-\widetilde{q}}_{C_\lambda([0,T];H^{s_1}(\omega))}\norm{\widetilde{u}_2}_{C([0,T];H^1(\Omega))}\int_0^Te^{-\lambda \tau}\tau^{-\frac{(2-s_0)}{2}}\d\tau\\
&\leq C\lambda^{-\frac{s_0}{2}}\left((1+Pe^{\lambda T})\norm{u_3}_{C_\lambda([0,T];H^2(\Omega))}+\norm{q-\widetilde{q}}_{C_\lambda([0,T];H^{s_1}(\omega))}(1+(Pe^{\lambda T})^4)e^{c(Pe^{\lambda T})^2T}\right).\end{aligned}$$
By choosing $P\leq P_0:=e^{-\lambda T}$ and $\lambda\geq \lambda_0=(4C)^{\frac{2}{s_0}}$, gives
$$\norm{u_2(q)-\widetilde{u}_2(\widetilde{q})}_{C_\lambda([0,T];H^2(\Omega))}\leq C\norm{q-\widetilde{q}}_{C_\lambda([0,T];H^{s_1}(\omega))}.$$
Thus, similarly to \eqref{tutu}, we find
$$\begin{aligned}&\norm{q\partial_{x_d}u_2(q)-\widetilde{q}\partial_{x_d}\widetilde{u}_2(\widetilde{q})}_{C_\lambda([0,T];H^{s_0}(\Omega))}\\
\leq& \norm{(q-\widetilde{q})\partial_{x_d}u_2(q)}_{C_\lambda([0,T];H^{s_0}(\Omega))}+\norm{\widetilde{q}(\partial_{x_d}u_2(q)-\partial_{x_d}\widetilde{u}_2(\widetilde{q}))}_{C_\lambda([0,T];H^{s_0}(\Omega))}\\
\leq& C\norm{q-\widetilde{q}}_{C_\lambda([0,T];H^{s_1}(\omega))}.\end{aligned}$$
By repeating the above argumentation, we obtain 
$$\norm{v_1-\widetilde{v}_1}_{C_\lambda([0,T];H^{2}(\Omega))}\leq C\norm{q-\widetilde{q}}_{C_\lambda([0,T];H^{s_1}(\omega))}$$
and, combining this with \eqref{t2aa} yields
$$\norm{\mathcal G(q)-\mathcal G(\widetilde{q})}_{C_\lambda([0,T];H^{s_1}(\omega))}\leq C\lambda^{\tau_0-1}\norm{q-\widetilde{q}}_{C_\lambda([0,T];H^{s_1}(\omega))}.$$
By fixing $\lambda_1:=(2C)^{\frac{1}{1-\tau_0}}$, $\lambda=\max(\lambda_0,\lambda_1)$ and   $\norm{\mathcal M_{q_1}(m)}_{C_\lambda([0,T];H^{s_1}(\omega))}=P= e^{-\lambda T}$,    $\mathcal G$ is a contraction of $B_{2P}^\lambda$. Thus, by assuming that $\norm{\mathcal M_{q_1}(m)}_{C([0,T];H^{s_1}(\omega))}=P\leq P_0$, we obtain  the existence of a solution of \eqref{int} in $C([0,T];H^{s_1}(\omega))$, and one can check the continuity of  such a solution with respect to $m$.

\subsection{Uniqueness of solutions}
We complete the proof of Theorem \ref{t2} by proving the uniqueness of solutions of the nonlinear equation \eqref{t2a} or equivalently \eqref{int}. Let $q,\widetilde{q}\in C([0,T];H^{s_1}(\omega))$ be two solutions of \eqref{int} and fix $P_1=\max(\norm{q}_{C([0,T];H^{s_1}(\omega))}, \norm{\widetilde{q}}_{C([0,T];H^{s_1}(\omega))})$. Note that 
\begin{align*}
&q(t)-\widetilde{q}(t)=\mathcal G(q)(t)-\mathcal G(\tilde{q})(t)\\
=&\int_0^tH(t-s)q_1(v_1-\widetilde{v}_1)(s,\cdot)\d s-\int_0^tH(t-s)(q\partial_{x_d}u_2(s,\cdot)-\widetilde{q}\partial_{x_d}\widetilde{u}_2(s,\cdot))\d s,\quad t\in(0,T).
\end{align*}
Therefore, by repeating the above argumentation, we 
obtain 
\begin{align*}
&\norm{q(t)-\widetilde{q}(t)}_{H^{s_1}(\omega)}=\norm{\mathcal Gq(t)-\mathcal G\widetilde{q}(t)}_{H^{s_1}(\omega)}\\
\leq& C(P_1^9+1)e^{cP_1^2} \int_0^t(t-s)^{-\tau_0}\norm{q(s)-\widetilde{q}(s)}_{H^{s_1}(\omega)}\d s,\quad t\in(0,T],
\end{align*}
with $C,c>0$ depending on $\Omega$, $g$, $u_0$, $q_1$ and $T$.
Let us introduce the following Gronwall's inequality which can be deduced from \cite[Lemma 7.1.1]{H} (see also \cite[Theorem 1]{YGD}).
\begin{lem}\label{l2}
  Let $M,r>0$ and $h,d\in L^1(0,T)$ be nonnegative functions satisfying
\[
	h(t)\le d(t)+M\int_0^t(t-s)^{r-1}h(s)\d s,\quad t\in(0,T).
\]
  Then there exist $C_1,C_2>0$ depending only on $T$ and $r$ such that 
\[
	h(t)\le d(t)+C_1\left(\sum_{k=1}^\infty \frac{(C_2M)^k}{\Gamma(kr)}\right)\int_0^t(t-s)^{r-1}d(s)\d s,\quad t\in(0,T).
\]

\end{lem}

From Lemma \ref{l2}, we deduce  that $q=\tilde{q}$. This proves the uniqueness of solutions of \eqref{int} and   completes the proof of Theorem \ref{t2}.

\section{Stability estimate}\label{sec:stab}

This section is devoted to the proof of  Theorem \ref{t4} and Corollary \ref{c3}.
We start with the following improved regularity result for the IBVP \eqref{eq10}.
\begin{prop}\label{p1} Let $q\in C^1([0,T];H^{s_1}(\omega))\cap C([0,T];W^{2,\infty}(\omega))$, $s_1\in(0,\frac12]$. Then, problem \eqref{eq10} admits a unique solution $u\in C^1([0,T];H^2(\Omega))\cap C([0,T];H^4(\Omega))$.
\end{prop}
\begin{proof} Note first that $u$ can be decomposed into $u=G+v$, with the function $G\in C^2([0,T];H^{s_0}(\Omega))\cap C^1([0,T];H^{2+s_0}(\Omega))\cap  C([0,T];H^{4}(\Omega))$ satisfies \eqref{trace_stab}
and $v$ solves \eqref{eq1} with $F=-(\partial_t-\Delta+q)G$. By Lemma \ref{l1} and \eqref{trace_stab}, we deduce that $F\in C^1([0,T];H^{s_0}(\Omega))\cap C([0,T];D(A))$, and $F(0,\cdot)\equiv0$. By Theorem \ref{t1},   $v\in C^1([0,T];L^2(\Omega))\cap C([0,T];H^2(\Omega))$. Moreover, $v_1=\partial_tv$ solves 
$$
\left\{\begin{aligned}\partial_t v_1-\Delta v_1+q(t,x')v_1=\widetilde{F}(t,x),\quad &(t,x',x_d)=(t,x)\in 
Q,\\  
v_1(t,x)=0,\quad &(t,x)\in(0,T)\times \partial\omega\times(0,\ell),\\ 
-\partial_{x_d}v_1(t,x', 0)=v_1(t,x', \ell)=0\quad &(t,x')\in(0,T)\times \omega,\\
v_1(0,x)=0,\quad &x\in \Omega,\end{aligned}\right.
$$
where 
$$  \widetilde{F}(t,x',x_d)=-\partial_tq(t,x')v(t,x)+\partial_tF(t,x),\quad (t,x',x_d)=(t,x)\in 
Q.$$
By Lemma \ref{l1} again, we deduce that $\widetilde{F}\in C([0,T];H^{s_0}(\Omega))$ and by Theorem \ref{t1},  $v_1\in C^1([0,T];L^2(\Omega))\cap C([0,T];H^2(\Omega))$. It follows that $v\in C^1([0,T];H^2(\Omega))$. In addition, we obtain
$$
\left\{\begin{aligned}-\Delta v=-v_1(t,x)-q(t,x')v(t,x)+F(t,x),\quad &(t,x',x_d)=(t,x)\in 
Q,\\  
v(t,x)=0,\quad &(t,x)\in(0,T)\times \partial\omega\times(0,\ell),\\ 
-\partial_{x_d}v(t,x', 0)=v(t,x', \ell)=0\quad &(t,x')\in(0,T)\times \omega.\end{aligned}\right.
$$
Since $(t,x)\mapsto -v_1(t,x)-q(t,x')v(t,x)+F(t,x)\in C([0,T];H^2(\Omega))$
and by \cite[Lemma 2.4]{GK}, we deduce   $v\in C([0,T];H^4(\Omega))$. This proves   $u\in C^1([0,T];H^{2}(\Omega))\cap  C([0,T];H^{4}(\Omega))$.\end{proof}

We can now prove Theorem \ref{t4} and Corollary \ref{c3}. 

\textbf{Proof of Theorem \ref{t4}.} Below $C>0$ denotes a constant depending on $\Omega$, $s_1$, $s_0$, $g$, $u_0$ and $T$ that might change from line to line.
For $j=1,2$, fix $q_j\in C^1([0,T];H^{s_1}(\omega))\cap C([0,T];W^{2,\infty}(\omega))$, and, by Proposition \ref{p1}, let $u_j\in C^1([0,T];H^2(\Omega))\cap C([0,T];H^4(\Omega))$ be the solution of \eqref{eq10} with $q=q_j$. Then, for
$m:(t,x')\mapsto \partial_{x_d}u_1(t,x',\ell)-\partial_{x_d}u_2(t,x',\ell)$, we have $m\in C^1([0,T];H^{\frac 1 2}(\omega))\cap C([0,T]; H^{\frac 5 2}(\omega))$. Moreover, for $u:=u_1-u_2$, we get 
$$
\left\{\begin{aligned}
\partial_{x_d}u(t,x)=0 ,\quad &(t,x)\in(0,T)\times \partial\omega\times(0,\ell),\\
\partial_{x_d}u(0,x)=0,\quad &x\in\Omega,\end{aligned}\right.$$
which implies that
$$
\left\{\begin{aligned}
m(t,x')=\partial_{x_d}u(t,x',\ell)=0 ,\quad &(t,x')\in(0,T)\times \partial\omega,\\ 
m(0,x')=\partial_{x_d}u(0,x',\ell)=0,\quad &x'\in\omega.\end{aligned}\right.$$
Thus, in view of \cite[Theorem 2.3, Chapter 4]{LM2}, there exists $G_m\in C^1([0,T];H^{s_0+\frac12}(\Omega))\cap C([0,T]; H^{2+s_0}(\Omega))$ satisfying the following conditions
\begin{equation}\label{trace3}
\left\{\begin{aligned}
G_m(t,x', \ell)=m(t,x')\quad &(t,x')\in(0,T)\times \omega,\\
G_m(t,x)=0 ,\quad &(t,x)\in(0,T)\times \partial\omega\times(0,\ell),\\ 
G_m(t,x', 0)=0\quad &(t,x')\in(0,T)\times \omega,\\
G_m(0,x)=0,\quad &x\in\Omega,\end{aligned}\right.\end{equation}
\begin{align}\label{trace4}
&\norm{G_m}_{W^{1,\infty}(0,T;H^{s_0+\frac{1}{2}}(\Omega))}+\norm{G_m}_{L^\infty(0,T;H^{2+s_0}(\Omega))}\\\nonumber
\leq& C\left(\norm{m}_{W^{1,\infty}(0,T; H^{s_0}(\omega))}+\norm{m}_{L^\infty(0,T;H^{\frac{3}{2}+s_0}(\omega))}\right).\end{align}
Following the discussion of Section \ref{sec:recon_local}, by fixing $u=u_1-u_2$ and letting $v=\partial_{x_d}u$, $v$ can be split into two terms $v=v_1+v_2$ with $v_1$ solving \eqref{eq5} and $v_2$ solving \eqref{eq6}. In view of \eqref{trace3}, we can split
$v_2$ into two terms $v_2=G_m+w$, where $w$ solves 
$$\left\{\begin{aligned}\partial_t w-\Delta w+q_1(t,x')w=-(\partial_t-\Delta +q_1)G_m(t,x),\quad &(t,x',x_d)=(t,x)\in 
Q,\\  
w(t,x)=0,\quad &(t,x)\in(0,T)\times \partial\omega\times(0,\ell),\\ 
w(t,x', 0)=w(t,x', \ell)=0\quad &(t,x')\in(0,T)\times \omega,\\
w(0,x)=0,\quad &x\in \Omega.\end{aligned}\right.$$
Thus, by Theorem \ref{t1} and Lemma \ref{l1}, we deduce   $w\in C([0,T];H^2(\Omega))\cap C^1([0,T];L^2(\Omega))$ and, in light of \eqref{trace4} and \eqref{t4a}, we get
\begin{align*}
&\norm{w}_{C([0,T];H^2(\Omega))}+ \norm{w}_{C^1([0,T];L^2(\Omega))}\\
\leq& C(1+M^4)e^{cM^2}\norm{(\partial_t-\Delta +q_1)G_m}_{C([0,T];H^{s_0}(\Omega))}\\
\leq&     C(1+M^5)e^{cM^2}(\norm{G_m}_{W^{1,\infty}(0,T;H^{s_0}(\Omega))}+\norm{G_m}_{L^\infty(0,T;H^{2+s_0}(\Omega))})\\
\leq&  C(1+M^5)e^{cM^2}(\norm{m}_{W^{1,\infty}(0,T; H^{s_0}(\omega))}+\norm{m}_{L^\infty(0,T;H^{\frac{3}{2}+s_0}(\omega))}).
\end{align*}
Using this estimate and \eqref{trace4}, we obtain
$$\norm{v_2}_{C([0,T];H^2(\Omega))}+ \norm{v_2}_{C^1([0,T];L^2(\Omega))}\leq C(1+M^5)e^{cM^2}(\norm{m}_{W^{1,\infty}(0,T; H^{s_0}(\omega))}+\norm{m}_{L^\infty(0,T;H^{\frac{3}{2}+s_0}(\omega))})$$
and then it follows that
\begin{equation}\label{t4d}\norm{\partial_{x_d}v_2|_{x_d=\ell}}_{C([0,T];H^{s_1}(\omega))}\leq C(1+M^5)e^{cM^2}(\norm{m}_{W^{1,\infty}(0,T; H^{s_0}(\omega))}+\norm{m}_{L^\infty(0,T;H^{\frac{3}{2}+s_0}(\omega))}).\end{equation}
By fixing $q=q_2-q_1$,   $q$ satisfies the following condition
\begin{equation}\label{t4e}q(t,x')=-\partial_{x_d}v_1(t,x',\ell)-\partial_{x_d}v_2(t,x',\ell),\quad (t,x')\in[0,T]\times\omega.\end{equation}
Using the notation of Section \ref{EC}, we define the operator valued map $Y$ by
$$Y(t)g=\sum_{n=1}^\infty e^{-\mu_n t}\left\langle g,\psi_n\right\rangle\psi_n,\quad t>0,\ g\in L^2(\Omega).$$
Then, we observe that $v_1\in   C([0,T];H^2(\Omega))$ is given by
$$v_1(t,\cdot)=-\int_0^tY(t-s)q_1v_1(s,\cdot)\d s+\int_0^tY(t-s)q\partial_{x_d}u_2(s,\cdot)\d s,\quad t\in[0,T].$$
Therefore, in view of Lemma \ref{l1} and the arguments used in the last step of Theorem \ref{t1}, fixing $\theta=\frac{1+2s_0-2s_1}{4}$,  for all $t\in[0,T]$, we obtain
\begin{align*}
&\norm{v_1(t,\cdot)}_{H^{\frac32 +s_1}(\Omega)}\\
\leq &CM\int_0^t(t-\tau)^{\theta-1}\norm{v_1(\tau,\cdot)}_{H^1(\Omega)}\d\tau+C\int_0^t(t-\tau)^{\theta-1}\norm{q(\tau,\cdot)}_{H^{s_1}(\omega)}\norm{u_2(\tau,\cdot)}_{H^2(\Omega)}\d\tau\\
\leq &CM\int_0^t(t-\tau)^{\theta-1}\norm{v_1(\tau,\cdot)}_{H^{\frac32 +s_1}(\Omega)}\d\tau+C(1+M^4)e^{cM^2}\int_0^t(t-\tau)^{\theta-1}\norm{q(\tau,\cdot)}_{H^{s_1}(\omega)}\d\tau.\end{align*}
By Lemma \ref{l2}, for all $t\in[0,T]$, we get
\begin{align*}
&\norm{\partial_{x_d}v_1(t,\cdot,\ell)}_{H^{s_1}(\omega)}\leq C\norm{v_1(t,\cdot)}_{H^{\frac32 +s_1}(\Omega)}\\
\leq& C(1+M^4)e^{cM^2}\int_0^t(t-\tau)^{\theta-1}\norm{q(\tau,\cdot)}_{H^{s_1}(\omega)}\d\tau\\
&\ \ \ + C(1+M^4)e^{cM^2}\left(\sum_{k=1}^\infty \frac{(C'M)^k}{\Gamma(k\theta)}\right)\int_0^t(t-\tau_1)^{\theta-1}\int_0^{\tau_1}(\tau_1-\tau_2)^{\theta-1}\norm{q(\tau_2,\cdot)}_{H^{s_1}(\omega)}\d\tau_2 \d\tau_1\\
\leq& C(1+M^4)e^{cM^2}\int_0^t(t-\tau)^{\theta-1}\norm{q(\tau,\cdot)}_{H^{s_1}(\omega)}\d\tau\\
&\ \ \ + C(1+M^4)e^{cM^2}\left(\sum_{k=1}^\infty \frac{(C'M)^k}{\Gamma(k\theta)}\right)\int_0^t(t-\tau)^{2\theta-1}\norm{q(\tau,\cdot)}_{H^{s_1}(\omega)}\d\tau\\
\leq& C(1+M^4)e^{cM^2}\left(1+\sum_{k=1}^\infty \frac{(C'M)^k}{\Gamma(k\theta)}\right)\int_0^t(t-\tau)^{\theta-1}\norm{q(\tau,\cdot)}_{H^{s_1}(\omega)}\d\tau,
\end{align*}
where $C'>0$ depends on $s_0$, $s_1$ and $T$. Combining this with  \eqref{t4d}-\eqref{t4e}, for all $t\in[0,T]$, we obtain the following
$$\begin{aligned}\norm{q(t,\cdot)}_{H^{s_1}(\omega)}&\leq C(1+M^5)e^{cM^2}(\norm{m}_{W^{1,\infty}(0,T; H^{s_0}(\omega))}+\norm{m}_{L^\infty(0,T;H^{\frac{3}{2}+s_0}(\omega))})\\
&\ \ \ + C(1+M^4)e^{cM^2}\left(1+\sum_{k=1}^\infty \frac{(C'M)^k}{\Gamma(k\theta)}\right)\int_0^t(t-\tau)^{\theta-1}\norm{q(\tau,\cdot)}_{H^{s_1}(\omega)}\d\tau.\end{aligned}
$$
Thus, applying again Lemma \ref{l2}, we get \eqref{t4b}.\qed

\textbf{Proof of Corollary \ref{c3}.}
Similarly to Theorem \ref{t4},   $q:=q_2-q_1$  satisfies the following condition
$$q(t,x')=-\partial_{x_d}v_1(t,x',\ell)-\partial_{x_d}v_2(t,x',\ell),\quad (t,x')\in[0,T]\times\omega.$$
Then, repeating the arguments of Theorem \ref{t3} with minor modifications, one can derive \eqref{t4b}.
\qed

\section{Reconstruction from global well-posedness under extra constraints}\label{sec:recon_global}
In this section, we prove Theorem \ref{t3}. Consider first the  following two useful intermediate results.

\begin{lem}\label{l4} The map $P_{\mathcal{Q}_b}$ is Lipschitz continuous in $C([0,T];L^2(\omega))$ and  
    \begin{equation}\label{eqn:cut-off_cts}
        \|P_{\mathcal{Q}_b} v(t,\cdot)- P_{\mathcal{Q}_b} w(t,\cdot) \|_{L^2(\omega)}\le   \|  v(t,\cdot)-w(t,\cdot) \|_{L^2(\omega)},\quad v,w\in C([0,T];L^2(\omega)),\ t\in[0,T].
    \end{equation} 
    
\end{lem}
\begin{proof} The inequality \eqref{eqn:cut-off_cts} is a direct consequence of the following well-known inequality
$|\min(\max(a,-b),b)-\min(\max(c,-b),b)|\leq|a-c|$ for $ a,c\in \R.$
\end{proof}

\begin{lem}\label{l5} The set $\mathcal{Q}_b\cap  C([0,T];L^2(\omega))$ is a closed subset of $ C([0,T];L^2(\omega))$.
\end{lem}
\begin{proof} Let $(f_n)_{n\in\mathbb N}$ be a sequence of $\mathcal{Q}_b\cap  C([0,T];L^2(\omega))$ that converges in $ C([0,T];L^2(\omega))$ to $f\in  C([0,T];L^2(\omega))$. Fix $t\in[0,T]$. Then, $(f_n(t,\cdot))_{n\in\mathbb N}$ converges to $f(t,\cdot)$ in $L^2(\omega)$ and, by eventually considering a subsequence, we may assume without loss of generality that $(f_n(t,\cdot))_{n\in\mathbb N}$ converges almost everywhere in $\omega$ to $f(t,\cdot)$. Meanwhile, since  $\{f_n:\ n\in\mathbb N\}\subset\mathcal{Q}_b$,  we have
$$-b\leq f_n(t,x)\leq b,\quad a.e.\  x\in\omega.$$
Thus, by taking the limit, we deduce  
$$-b\leq f(t,x)\leq b,\quad a.e.\  x\in\omega,$$
which implies that $\|f(t,\cdot)\|_{L^\infty(\omega)}\leq b$. Taking the essential supremum with respect to $t\in(0,T)$, we deduce that $f\in \mathcal{Q}_b$ which implies that  $\mathcal{Q}_b\cap  C([0,T];L^2(\omega))$ is a closed subset of $ C([0,T];L^2(\omega))$.
\end{proof}
Now we can prove Theorem \ref{t3}.\\
\ \\
 \textbf{Proof of Theorem \ref{t3}.} We divide the proof of the theorem into two steps.
  
\textbf{Step 1.} In this step, we prove the first claim of the theorem. Specifically, we establish the equivalence between  \textbf{(ICP2)} and equation \eqref{t3a}. Let $q^\dagger\in \mathcal Q_b\cap C([0,T];L^2(\omega))$ be the unknown coefficient of \textbf{(ICP2)} and let $u^\dagger=u(q^\dagger)$ be the solution to problem \eqref{eq3}, with $q_j=q^\dagger$. By  Corollary \ref{cor:forward_weak_reg},  we deduce that $u^\dagger=u(q^\dagger)\in C([0,T];H^{\frac{7}{4}}(\Omega))$ and   the flux measurement $m(t,x')\in C([0,T];H^{\frac{1}{4}}(\omega))$. In light of Assumption \ref{assum:alg},  we can split $\partial_{x_d} u^\dagger$ into two terms $\partial_{x_d}u^\dagger =G_m+v(m,q^\dagger)$, where $v=v(m,q^\dagger)$ solves 
$$\left\{\begin{aligned}\partial_t v-\Delta v+q^\dagger(t,x')v=-(\partial_t-\Delta +q^\dagger)G_m,\quad &(t,x',x_d)=(t,x)\in 
Q,\\  
v(t,x)=0,\quad &(t,x)\in(0,T)\times \partial\omega\times(0,\ell),\\ 
v(t,x', 0)=v(t,x', \ell)=0\quad &(t,x')\in(0,T)\times \omega,\\
v(0,x)=0,\quad &x\in \Omega.\end{aligned}\right.$$
Since  $G_m\in C([0,T];H^2(\Omega))\cap C^1([0,T];L^2(\Omega))$, $q^\dagger\in L^\infty(0,T;L^\infty(\omega))$ and by Corollary \ref{cor:forward_weak_reg}, we deduce   $\partial_{x_d}u^\dagger\in C([0,T];H^{\frac{7}{4}}(\Omega))$. 
Similarly to the proof of Theorem \ref{t2}, taking the trace on $x_d=\ell$ of the equation \eqref{eq3}, we obtain
$$\begin{aligned}q^\dagger(t,x')=q^\dagger(t,x')u^\dagger(t,x',\ell)&=-(\partial_tu^\dagger(t,x',\ell)-\Delta u^\dagger(t,x',\ell))\\
&=-(\underbrace{\partial_tu^\dagger(t,x',\ell)-\Delta' u^\dagger(t,x',\ell)}_{=0})+\underbrace{\partial_{x_d}^2u^\dagger(t,x',\ell)}_{=\partial_{x_d}(G_m+v(m,q^\dagger))(t,x',\ell)}\\
&=\partial_{x_d}G_m(t,x',\ell)+\partial_{x_d}v(m,q^\dagger)(t,x',\ell),\quad (t,x')\in[0,T]\times \omega,\end{aligned}$$
where $\Delta':=\sum_{j=1}^{d-1}\partial_{x_j}^2$. 
This proves that \textbf{(ICP2)}, for  $q\in \mathcal{Q}_b\cap  C([0,T];L^2(\omega))$, is equivalent to solving in $\mathcal{Q}_b\cap  C([0,T];L^2(\omega))$  the following nonlinear equation
\begin{equation}\label{eqn:int_P}
q^\dagger=P_{\mathcal{Q}_b}q^\dagger=P_{\mathcal{Q}_b}\left[\partial_{x_d}G_m|_{x_d=\ell}+\partial_{x_d}v(m,q^\dagger)|_{x_d=\ell}\right].
\end{equation}
This completes the proof of the first claim of the theorem by fixing
$$\mathcal S(m)=\partial_{x_d}G_m|_{x_d=\ell},\quad \mathcal R(G_m,q)=\partial_{x_d}v(m,q)|_{x_d=\ell},\quad q\in \mathcal{Q}_b\cap  C([0,T];L^2(\omega)).$$ 
\textbf{Step 2.} In this step, we prove the well-posedness of \eqref{t3a}, or equivalently
\eqref{eqn:int_P}, in $\mathcal{Q}_b\cap  C([0,T];L^2(\omega))$. Consider the map $\mathcal{G}: \mathcal{Q}_b\cap  C([0,T];L^2(\omega))\mapsto \mathcal{Q}_b\cap  C([0,T];L^2(\omega))$ defined by 
$$\mathcal{G} q=P_{\mathcal{Q}_b}\left[\partial_{x_d}G_m|_{x_d=\ell}+\partial_{x_d}v(m,q)|_{x_d=\ell}\right],\quad q\in  \mathcal{Q}_b\cap C([0,T];L^2(\omega)).$$
By Lemma \ref{l5},   $\mathcal{Q}_b\cap  C([0,T];L^2(\omega))$ is a closed subset of $ C([0,T];L^2(\omega))$. We only need to show that $\mathcal{G} $ is a contraction map with respect to the topology of $C_\lambda([0,T];L^2(\omega))$ for $\lambda>0$ sufficiently large. Fix $q, \tilde{q}\in \mathcal{Q}_b\cap  C([0,T];L^2(\omega))$. In view of Lemma \ref{l4},  for fixed $\lambda>0$, we obtain
\begin{align}\label{t3b}
 &\norm{\mathcal{G} q-\mathcal{G} \widetilde{q}}_{C_\lambda([0,T];L^2(\omega))} =\sup_{t\in[0,T]}e^{-\lambda t}\norm{\mathcal{G} q(t,\cdot)-\mathcal{G} \widetilde{q}(t,\cdot)}_{L^2(\omega)}\notag\\
&=\sup_{t\in[0,T]}e^{-\lambda t}\norm{P_{\mathcal{Q}_b}\left[\partial_{x_d}G_m|_{x_d=\ell}+\partial_{x_d}v(m,q)|_{x_d=\ell}\right](t,\cdot)-P_{\mathcal{Q}_b}\left[\partial_{x_d}G_m|_{x_d=\ell}+\partial_{x_d}v(m,\tilde{q})|_{x_d=\ell}\right](t,\cdot)}_{L^2(\omega)}\notag\\
&\leq \sup_{t\in[0,T]}e^{-\lambda t}\norm{\partial_{x_d}v(m,q)|_{x_d=\ell}(t,\cdot)-\partial_{x_d}v(m,\tilde{q}))|_{x_d=\ell}(t,\cdot)}_{L^2(\omega)}\notag\\
&\leq \norm{\partial_{x_d}v(m,q)|_{x_d=\ell}-\partial_{x_d}v(m,\tilde{q})|_{x_d=\ell}}_{C_\lambda([0,T];L^2(\omega))}. \end{align}
Then, $w=v(m,q)-v(m,\widetilde{q})$ solves 
$$\left\{\begin{aligned}\partial_t w-\Delta w+q(t,x')w=-(q-\widetilde{q})(G_m+v(m,\widetilde{q})),\quad &(t,x',x_d)=(t,x)\in 
Q,\\  
w(t,x)=0,\quad &(t,x)\in(0,T)\times \partial\omega\times(0,\ell),\\ 
w(t,x', 0)=w(t,x', \ell)=0\quad &(t,x')\in(0,T)\times \omega,\\
w(0,x)=0,\quad &x\in \Omega.\end{aligned}\right.$$
Since $S(t)f:=\sum_{n=1}^\infty e^{-\mu_n t} \left\langle f,\psi_n\right\rangle\psi_n$, $t\geq0$ and $f\in L^2(\Omega)$, we get 
\begin{equation*}
    w(t)=\int_0^t S(t-s) [(\widetilde{q}-q)(G_m+v(m,\widetilde{q}))](s) \mathrm{d} s-\int_0^t S(t-s) [qw](s) \mathrm{d} s,\quad t\in(0,T).
\end{equation*}
Since $q,\tilde{q}\in \mathcal{Q}_b$, we have $\max(\norm{q}_{L^\infty((0,T)\times\omega)},\norm{\tilde{q}}_{L^\infty((0,T)\times\omega)})\leq b $. Thus, by Corollary \ref{cor:forward_weak_reg} and   Sobolev embedding theorem, we can   find $C>0$ depending on $b$, $\Omega$, $T$, $m$, $g$ and $u_0$ such that 
$$\norm{v(m,\widetilde{q})}_{L^\infty([0,T];L^\infty(\Omega))}\leq C\norm{v(m,\widetilde{q})}_{C([0,T];H^\frac{7}{4}(\Omega))}\leq C.$$
By combining this estimate with the  regularity properties of $G_m$, the estimate \eqref{t2c} and  the  Sobolev embedding theorem, we can find $\epsilon\in(0,1)$ and  $C>0$ depending on $\epsilon$, $b$, $\Omega$, $T$, $m$, $g$ and $u_0$ such that, for all $t\in[0,T]$, we have 
\begin{align*}
    \|w(t)\|_{H^{\frac{3+\epsilon}{2}}(\Omega)}
    \le& C\int_0^t (t-s)^{-\frac{3+\epsilon}{4}}\|(\widetilde{q}-q)(s)\|_{L^2(\omega)}\mathrm{d} s + C \int_0^t (t-s)^{-\frac{(3+\epsilon)}{4}}\|w(s)\|_{L^2(\Omega)} \mathrm{d}s\\
     \le& C\int_0^t (t-s)^{-\frac{3+\epsilon}{4}}\|(\widetilde{q}-q)(s)\|_{L^2(\omega)}\mathrm{d} s + C \int_0^t (t-s)^{-\frac{3+\epsilon}{4}}\|w(s)\|_{H^{\frac{3+\epsilon}{2}}(\Omega)} \mathrm{d}s.
\end{align*}
Therefore, by multiplying $e^{-\lambda t}$ on both sides, we can find  $C>0$ depending on $\epsilon$, $b$, $\Omega$, $T$, $m$, $g$ and $u_0$, such that 
\begin{align*}
 \|w \|_{C_\lambda([0,T];H^{\frac{3+\epsilon}{2}}(\Omega))}
    \le& C\sup_{t\in[0,T]}\int_0^t e^{-\lambda(t-s)}(t-s)^{-\frac{3+\epsilon}{4}}e^{-\lambda s}\|(\widetilde{q}-q)(s)\|_{L^2(\omega)}\mathrm{d} s\\
    &+C\sup_{t\in[0,T]}\int_0^t e^{-\lambda(t-s)}(t-s)^{-\frac{3+\epsilon}{4}}e^{-\lambda s}\|w(s)\|_{H^{\frac{3+\epsilon}{2}}(\Omega)}\mathrm{d} s\\
    \le& C\lambda^{-\frac{1-\epsilon}{4}}\|q-\widetilde{q} \|_{C_\lambda([0,T];L^2(\omega))}+ C\lambda^{-\frac{1-\epsilon}{4}}\|w\|_{C_\lambda([0,T];H^{\frac{3+\epsilon}{2}}(\Omega))}.
\end{align*} 
Thus, there exists $\lambda_0>0$  sufficiently large, depending only on $\epsilon$, $b$, $\Omega$, $T$, $m$, $g$ and $u_0$, such that for $\lambda\geq\lambda_0$, we have $C\lambda^{-\frac{1-\epsilon}{4}}\leq\frac12$ and 
\begin{equation}\label{t3d}\|w \|_{C_\lambda([0,T];H^{\frac{3+\epsilon}{2}}(\Omega))}\leq 2C\lambda^{-\frac{1-\epsilon}{4}}\|q-\widetilde{q} \|_{C_\lambda([0,T];L^2(\omega))}.\end{equation}
Meanwhile, there exist $C_1,C_2>0$ depending only on $\Omega$ such that
$$\|\partial_{x_d}w|_{x_d=\ell} \|_{C_\lambda([0,T];L^2(\omega))}\leq C_1\|\partial_{x_d}w|_{x_d=\ell} \|_{C_\lambda([0,T];H^{\frac\epsilon2}(\omega))}
\leq C_2\|w \|_{C_\lambda([0,T];H^{\frac{3+\epsilon}{2}}(\Omega))}.$$
By combining this with \eqref{t3b}-\eqref{t3d}, we deduce  
$$\norm{\mathcal{G} q-\mathcal{G} \widetilde{q}}_{C_\lambda([0,T];L^2(\omega))}\leq \|\partial_{x_d}w|_{x_d=\ell} \|_{C_\lambda([0,T];L^2(\omega))}\leq C\lambda^{-\frac{1-\epsilon}{4}}\|q-\widetilde{q} \|_{C_\lambda([0,T];L^2(\omega))},$$
with $C>0$ depending on $\epsilon$, $b$, $\Omega$, $T$, $m$, $g$ and $u_0$. Then, we can find $\lambda_1\geq\lambda_0$  sufficiently large, depending only on $\epsilon$, $b$, $\Omega$, $T$, $m$, $g$ and $u_0$, such that for $\lambda=\lambda_1$, we obtain
\begin{equation}\label{t3c}\norm{\mathcal{G} q-\mathcal{G} \widetilde{q}}_{C_\lambda([0,T];L^2(\omega))}\leq \tfrac{1}{2}\norm{q-\widetilde{q}}_{C_\lambda([0,T];L^2(\omega))}.\end{equation}
Therefore, Banach fixed point theorem implies the existence of a unique fixed point $q^\star\in\mathcal{Q}_b\cap  C([0,T];L^2(\omega))$ of $\mathcal{G}$. Using this property, similar to Theorem \ref{t2}, we obtain the well-posedness of \eqref{t3a}. \qed

\begin{rmk}\label{rmk:recon}
    Theorem \ref{t3} provides an explicit numerical algorithm for reconstructing the coefficient $q$.  For $q^0\in \mathcal{Q}_b \cap C([0,T];L^2(\omega))$, consider 
    \begin{equation}\label{eqn:recon}
        q^{k+1}(x^\prime,t)= P_{\mathcal{Q}_b}\left[  \partial_{x_d}v(q^k) (t,x^\prime,\ell)\right],\quad k=0,1,\dots,
    \end{equation} 
    where $v=v(q)$ is the solution to 
    \begin{equation}\label{eqn:dudx}
     \left\{\begin{aligned}
         \partial_t v -\Delta v+q(t,x')v  &=0, && \text{in }  (0,T)\times\Omega,\\
         v(0,x)&=\partial_{x_d} u_0(x),&& x\in  \Omega\\
         v(t,x',x_d)&=\partial_{x_d} g(t,x',x_d) && (t,x',x_d)\in  (0,T)\times \partial\omega\times(0,\ell)\,\\
        v(t,x',0)&=0,&&    (t,x')\in     (0,T)\times\omega,\\
        v(t,x',\ell)&=m(t,x'),&&  (t,x')\in   (0,T)\times\omega.\\
      \end{aligned}\right.
    \end{equation} 
    Theorem \ref{t3} implies that for sufficiently large $\lambda$, the sequence $\{q^k\}_{k\ge0}$ generated by \eqref{eqn:recon} converges to $q^\dagger$ in $\mathcal{Q}_b\cap C([0,T];L^2(\omega))$ and, for all $\epsilon\in(0,1)$, the following holds
\begin{equation*}
        \|q^k-q^\dagger\|_{C_\lambda([0,T];L^2(\omega))}\le \big(c\lambda^{ -\frac{1-\epsilon}{4}}\big)^k  \|q^0-q^\dagger\|_{C_\lambda([0,T];L^2(\omega))},
    \end{equation*}
    where $c>0$ depends only on $\epsilon,u_0,g,b, m,\Omega,T$.
\end{rmk}

\section{Numerical scheme and error estimate}\label{sec:error}
In this section, we study the numerical  analysis for \textbf{(ICP2)}. We develop a discrete reconstruction scheme based on the formulas \eqref{eqn:recon}-\eqref{eqn:dudx}, and analyze the convergence of the   scheme and   the error between the reconstructed and exact coefficients.

\subsection{Numerical scheme for (ICP2)}\label{subsec:recon_scheme} 
We first present the fully discrete scheme for problem \eqref{eqn:dudx}. For time discretization,  we divide the interval $[0, T]$ into $N$ uniform subintervals with a step size $\tau = T/N$ and set the time grids $\{t_n=n\tau\}_{n=0}^{N}$. We assume that $\tau b<1$, which guarantees the coercivity of
the bilinear form at each time step for coefficients satisfying
$|q|\le b$. We employ the backward Euler scheme \cite{Thomee:2006} to approximate the time derivative $\partial_t  v(t_n)$ (with $v^n:=v(t_n)$):
\begin{equation*}
    \bar{\partial}_\tau  v^n:= \tau^{-1}(v^n-v^{n-1}).
\end{equation*} 

For spatial discretization, we employ the Galerkin finite element method (FEM) \cite{Thomee:2006}. Let $\mathcal{T}_h$ be a quasi-uniform simplicial triangulation of domain $\Omega$  with a mesh size $h$. The mesh $\mathcal{T}_h$ is assumed to conform exactly to the domain \(\Omega\). Otherwise, the estimate can be extended using curved elements when \(d=2\) \cite{Zlamal:1973,Zlamal:1974}, or isoparametric elements when \(d\ge2\) \cite{Ciarlet:1972,Lenoir:1986}.
Over $\mathcal{T}_h$, we define a continuous piecewise linear finite element space $V_h$ by $V_h=\{\varphi_h\in H^1(\Omega):\varphi_h|_{T} \text{ is a linear function }\forall \,T\in\mathcal{T}_h\}$ and set $V_{h}^{0}=\{\varphi_h\in V_h: \varphi_h=0 \text{ on } \partial\Omega  \}$.  
The following inverse inequality holds on the finite element space $V_h$ \cite[Lemma 4.5.3]{Brenner:2008}: for $0\le r\le s \le 1$ and $1\le q\le p\le \infty$ we have
\begin{align}\label{eqn:inverse_ineq}
	\|\varphi_h\|_{W^{s,p}(\Omega)}\le Ch^{r-s+\frac{d}{p}-\frac{d}{q}}\| \varphi_h\|_{W^{r,q}(\Omega)},\quad \forall\varphi_h\in V_h.
\end{align}
For all $T\in \mathcal{T}_h$, the following trace inequalities hold \cite[Lemmas 1.46 and 1.49]{Pietro:2012book}: 
\begin{align}
\|\varphi_h\|_{L^2(\partial T)}\le& ch^{-\frac{1}{2}}\|\varphi_h\|_{L^2(T)},\quad \forall \varphi_h\in V_h, \label{eqn:dis_trace}\\
\|\varphi\|_{L^2(\partial T)}\le& c( h^{\frac{1}{2}}\|\nabla\varphi\|_{L^2(T)}+h^{-\frac{1}{2}}\|\varphi\|_{L^2(T)} ),\quad \forall \varphi\in H^1(T),\label{eqn:cts_trace_ineq}
\end{align}
with $c$ independent of $h$ and $T$. 
Let $\mathcal T_h^\partial$ be the boundary mesh induced by
$\mathcal T_h$ on $\partial\Omega$, and let
$V_h^\partial=\{v_h|_{\partial\Omega}:v_h\in V_h\}$.

We denote by $\mathcal I_h$ the Scott--Zhang interpolation operator \cite{ScottZhang:1990,Brenner:2008}. For each
boundary node, we choose the associated averaging simplex to be a boundary face. With this choice, the trace of $\mathcal I_hv$ depends only on the trace of $v$, and $\mathcal I_h$ maps $H_0^1(\Omega)$ into $V_h^0$. The following interpolation error estimate holds
\cite[Theorem 4.8.12]{Brenner:2008}: for $s\in[1,2]$
and $p\in[1,\infty]$,
\begin{align}\label{eqn:error_I_h}
\|v-\mathcal I_hv\|_{L^p(\Omega)}
+h\|\nabla(v-\mathcal I_hv)\|_{L^p(\Omega)}
\le Ch^s\|v\|_{W^{s,p}(\Omega)},
\quad \forall v\in W^{s,p}(\Omega).
\end{align}
We define the $L^2(\Omega)$-projection $P_h\colon L^2(\Omega)\rightarrow  V_h$ by 
\begin{equation*}
    (P_h v,\varphi_h)=( v,\varphi_h),\quad\forall\varphi_h\in V_h,
\end{equation*}
and Ritz-projection $R_h\colon H^1(\Omega)\rightarrow  V_h$ by
\begin{align*}
	 (\nabla R_h v,\nabla \varphi_h)=( \nabla v,\nabla \varphi_h),\quad\forall\varphi_h\in V_h^0\quad \text{and} \quad R_hv=\mathcal{I}_h v \text{ on }\partial\Omega.
\end{align*}
The following error estimates hold \cite[p. 32]{Thomee:2006}: for every $s\in [1,2]$, we have
\begin{align}\label{eqn:error_P_h}
	\| v-P_h v\|_{L^2(\Omega)}+h\|\nabla(v-P_h v)\|_{L^2(\Omega)}\le& Ch^s\| v\|_{H^s(\Omega)},\quad\forall v \in H^s(\Omega) ,\\
    \| v-R_h v\|_{L^2(\Omega)}+h\|\nabla(v-R_h v)\|_{L^2(\Omega)}\le& Ch^s\| v\|_{H^s(\Omega)},\quad\forall v \in H^s(\Omega)\cap H_0^1(\Omega).\label{eqn:error_R_h}
\end{align}
We define the discrete Dirichlet Laplacian  $-\Delta_h\colon V_h^{0}\to V_h^{0}$   by
\begin{align*}
    (-\Delta_h v_h,\varphi_h)&= (\nabla v_h, \nabla \varphi_h),\quad \forall v_h,\varphi_h\in V_h^{0}.
\end{align*} 
We define a boundary interpolation operator $\mathcal I_h^\partial:H^{\frac12}(\partial\Omega)\to V_h^\partial$ as follows. For any $v\in H^{\frac12}(\partial\Omega)$, let $\widetilde v\in H^1(\Omega)$ be an extension satisfying $\widetilde v|_{\partial\Omega}=v$, and set $\mathcal I_h^\partial v:=(\mathcal I_h\widetilde v)|_{\partial\Omega}$. Since the averaging simplices associated with boundary nodes lie on $\partial\Omega$, this definition is independent of the choice of extension. The standard scaling argument gives
\[
\|\mathcal I_h^\partial v\|_{L^2(\partial\Omega)}\le c\|v\|_{L^2(\partial\Omega)},\quad \forall v\in H^{\frac12}(\partial\Omega).
\]
By density, $\mathcal I_h^\partial$ extends uniquely to a bounded linear operator from $L^2(\partial\Omega)$ to $V_h^\partial$, still denoted by $\mathcal{I}_h^\partial$. This extension is also given by the boundary-face averaging formulas and is stable in $L^\infty(\partial\Omega)$. Thus,
\begin{equation}\label{eqn:Ih_bdry_stab}
\|\mathcal I_h^\partial v\|_{L^p(\partial\Omega)}\le c\|v\|_{L^p(\partial\Omega)},\quad \forall v\in L^p(\partial\Omega),\quad p\in [2,\infty] ,
\end{equation}
where $c$ is independent of $h$.

For the error analysis, we strengthen Assumption \ref{assum:alg} by imposing the following additional regularity and compatibility conditions. The regularity requirements on the exact solution can be ensured under sufficient regularity of the exact coefficient $q^\dagger$ and the problem data, together with suitable compatibility conditions, by an argument similar to that used in Proposition \ref{p1}.

\begin{assumption}\label{assum:reg_err}  
The exact coefficient $q^\dagger\in W^{1,\infty}(0,T;L^\infty(\omega))$. The exact   data $m^\dagger:=\partial_{x_d}u(q^\dagger)|_{x_d=\ell}$ satisfy $ m^\dagger\in   C([0,T];H^{\frac{3}{2}}(\omega))\cap C^1([0,T];H^1(\omega))\cap C^2([0,T]; L^2(\omega))$.   
Moreover, we define the full boundary datum by
\[
b^\dagger(t,x)=
\begin{cases}
\partial_{x_d}g(t,x),
&x\in\partial\omega\times(0,\ell),\\
0,
&x\in\omega\times\{0\},\\
m^\dagger(t,x'),
&x=(x',\ell)\in\omega\times\{\ell\}.
\end{cases}
\]
Suppose that $b^\dagger\in C^2(0,T;L^2(\partial\Omega))$ and that there exists $G\in  C([0,T]; H^2(\Omega))\cap C^1([0,T];H^\frac{3}{2}(\Omega))\cap C^2([0,T];H^{\frac{1}{2}}(\Omega))\cap W^{1,1}(0,T;L^\infty(\Omega))$ satisfying $\partial_{x_d}u_0-G(0)\in D(-\Delta)$, where $-\Delta$  denotes the Dirichlet Laplacian on  $L^2(\Omega)$, and
\begin{equation} \label{eqn:bdry_lift}
     \left\{\begin{aligned}
     -\Delta G(t,x)  &=0, &&    (t,x)\in [0,T]\times\Omega ,\\
         G(t,x)&=b^\dagger(t,x), && (t,x)\in [0,T]\times\partial\Omega.
      \end{aligned}\right.
    \end{equation}  
     The noisy  data  $z^\delta(t,x^\prime)\in C([0,T];H^{\frac{1}{2}}(\omega))\cap W^{1,1}(0,T;L^\infty(\omega)  ) $ satisfy
\begin{equation*}
        \|z^\delta-m^\dagger\|_{C([0,T];L^2(\omega))}\le \delta\quad \text{and}\quad {\|z^\delta\|_{C([0,T];H^{\frac{1}{2}}(\omega))}+\|z^\delta\|_{W^{1,1}(0,T;L^\infty(\omega)  )}\le c. }
    \end{equation*} 
Moreover, the full noisy boundary datum $b^\delta$ defined by
\[
b^\delta(t,x)=
\begin{cases}
\partial_{x_d}g(t,x),
&x\in\partial\omega\times(0,\ell),\\
0,
&x\in\omega\times\{0\},\\
z^\delta(t,x'),
&x=(x',\ell)\in\omega\times\{\ell\}
\end{cases}
\]
belongs to $ C([0,T];H^{\frac12}(\partial\Omega))$, with
$\|b^\delta\|_{C([0,T];H^{\frac12}(\partial\Omega))}\le c,$
and the initial compatibility condition
$z^\delta(0)=m^\dagger(0)
=\left.\partial_{x_d}u_0\right|_{\omega\times\{\ell\}}$ holds in $L^2(\omega)$.
The relevant constants are independent of $\delta$.
\end{assumption}

Now we can establish the discrete reconstruction scheme. For any initial guess $q^0=\{q^{0,n}\}_{n=0}^{N}$  satisfying
$|q^{0,n}|\le b$ a.e. in $\omega$, consider the following iteration scheme, with updates $q^{k+1}=\{q^{k+1,n}\}_{n=0}^{N}$ from $q^k=\{q^{k,n}\}_{n=0}^{N}$ by 
\begin{equation}\label{eqn:recon_dis}
    q^{k+1,n}(x')=P_{\mathcal{Q}_b}\left[   \partial_{x_d}v_h^n(q^{k}) (x^\prime,\ell)\right],
\end{equation}
where $P_{\mathcal{Q}_b}$ is the cut-off operator defined in \eqref{eqn:cut-off-op} and $v_h^n(q)$ is the fully discrete approximation to the solution
of  problem \eqref{eqn:dudx}:  given $v_h^0=R_h (\partial_{x_d}u_0)$, find $v_h^n\in V_h $ for $n=1,\dots,N$ such that  
\begin{equation}\label{eqn:fully-dis-v}
\left\{
\begin{aligned}
(\bar\partial_\tau v_h^n,\varphi_h)
+(\nabla v_h^n,\nabla\varphi_h)
+(q^n v_h^n,\varphi_h)
&=0,
&&\forall\varphi_h\in V_h^0,\\
v_h^n
&=\mathcal I_h^\partial b^\delta(t_n),
&&\text{on }\partial\Omega.
\end{aligned}
\right.
\end{equation}
Here $b^\delta$ is the noisy full boundary datum defined in
Assumption \ref{assum:reg_err}, and $\mathcal I_h^\partial$
is the full-boundary interpolation operator. 

We define the discrete harmonic lifting $G_h^n=G_h(t_n)$
for $n=0,\dots,N$, with $G_h(t)\in V_h$ satisfying
\begin{equation}\label{eqn:bdry_lift_dis}
\left\{
\begin{aligned}
(\nabla G_h(t),\nabla\varphi_h)
&=0,
&&\forall\varphi_h\in V_h^0,\\
G_h(t)
&=\mathcal I_h^\partial b^\delta(t),
&&\text{on }\partial\Omega,
\end{aligned}
\right.
\quad t\in[0,T].
\end{equation}
By Assumption \ref{assum:reg_err},
$b^\delta(0)=G(0)|_{\partial\Omega}$.
Since $G(0)$ is harmonic, the definitions of $R_h$ and $G_h(0)$ give $G_h^0=R_hG(0)$. Consequently, $ v_h^0-G_h^0 = R_h\bigl(\partial_{x_d}u_0-G(0)\bigr)
\in V_h^0$.

Let $w_h^n=v_h^n-G_h^n$. Then $w_h^n\in V_h^0$
for $n=0,\dots,N$. 
Let $P_h^0:L^2(\Omega)\to V_h^0$ denote the
$L^2(\Omega)$-orthogonal projection:
\begin{equation*}
(P_h^0\psi,\varphi_h)=(\psi,\varphi_h),
\quad \forall\psi\in L^2(\Omega),\quad
\forall\varphi_h\in V_h^0.
\end{equation*}
By \eqref{eqn:fully-dis-v} and \eqref{eqn:bdry_lift_dis},
we obtain
\begin{equation*}
\bar\partial_\tau w_h^n-\Delta_h w_h^n
=
-P_h^0\bigl(\bar\partial_\tau G_h^n+q^n v_h^n\bigr),
\quad n=1,\dots,N,
\end{equation*}
with $w_h^0=R_h(\partial_{x_d}u_0-G(0))\in V_h^0$. By the discrete Laplace transform, we obtain the following
solution representation \cite[Chapter 3]{JinZhou:2023book}:
\begin{equation}\label{eqn:dis-solrep-w}
v_h^n
=
G_h^n
+E_{h,\tau}^{n-1}R_h\bigl(\partial_{x_d}u_0-G(0)\bigr)
-\tau\sum_{j=1}^{n}E_{h,\tau}^{n-j}
P_h^0\bigl(\bar\partial_\tau G_h^j+q^jv_h^j\bigr),
\quad n=1,\dots,N,
\end{equation}
where the discrete solution operator
$E_{h,\tau}^n:V_h^0\to V_h^0$ is defined by
\[
E_{h,\tau}^n
=
\frac{1}{2\pi\mathrm{i}}
\int_{\Gamma_{\theta,\sigma}^\tau}
e^{zt_n}
\bigl(\delta_\tau(e^{-z\tau})-\Delta_h\bigr)^{-1}
\,\mathrm{d}z,
\quad n\ge0,
\]
with $\delta_\tau(\xi)=(1-\xi)/\tau$ and
$\Gamma_{\theta,\sigma}^\tau
:=\{z\in\Gamma_{\theta,\sigma}:|\operatorname{Im}z|\le\pi/\tau\}$,
where $\theta\in(\frac{\pi}{2},\pi)$ is chosen sufficiently
close to $\frac{\pi}{2}$ and the contour is oriented
with an increasing imaginary part.

The temporal discretization satisfies the following properties:
for any $\theta\in(\frac{\pi}{2},\pi)$, there exists
$\theta'\in(\frac{\pi}{2},\pi)$ such that, for all
$z\in\Gamma_{\theta,\sigma}^\tau$
\cite[Lemma 3.1]{JinZhou:2023book},
\begin{equation}\label{eqn:delta_tau}
\delta_\tau(e^{-z\tau})\in\Sigma_{\theta'},
\quad
c_1|z|\le|\delta_\tau(e^{-z\tau})|\le c_2|z|,
\quad \text{and}\quad
|\delta_\tau(e^{-z\tau})-z|\le c_3\tau|z|^2,
\end{equation}
where $c_1,c_2,c_3>0$ are independent of $\tau$.
By the positivity and self-adjointness of $-\Delta_h$
and \eqref{eqn:delta_tau}, we obtain the following discrete
resolvent estimate for $\beta\in[0,1]$:
\begin{equation}\label{eqn:dis-resol}
\left\|
(-\Delta_h)^\beta
\bigl(\delta_\tau(e^{-z\tau})-\Delta_h\bigr)^{-1} \varphi_h
\right\|_{L^2(\Omega) }
\le c|z|^{\beta-1} \| \varphi_h \|_{L^2(\Omega)},
\quad \forall ~z\in\Gamma_{\theta,\sigma}^\tau,~~ \varphi_h \in V_h^0,
\end{equation}
where the fractional powers of $-\Delta_h$ are defined by spectral expansion.

Now, we describe smoothing properties for the discrete solution operators $ E_{h,\tau}^n $. The proof of the $L^2(\Omega)$ estimate follows from \cite[Lemma 4.3]{ZhangZhangZhou:2022} and  interpolation theory. The proof of the $L^\infty(\Omega)$ estimate follows from the discrete $L^\infty(\Omega)$ resolvent estimate \cite[Theorem 1.1]{LiSun:2017} and the estimate \eqref{eqn:delta_tau}.
\begin{lem}\label{lem:dis-sol-op}
For the operator  $E_{h,\tau}^n$ and any $\beta\in [0,1]$, there exists   $c>0$   independent of $h,\,\tau$ and $t_n$ such that 
\begin{equation*}
    t_{n+1}^{ \beta }\|(-\Delta_h)^\beta E_{h,\tau}^n \varphi_h\|_{L^2(\Omega)} \le c \|\varphi_h\|_{L^2(\Omega)}\quad \text{and}\quad  t_{n+1}^\beta \|(-\Delta_h)^{\beta}  E_{h,\tau}^n \varphi_h\|_{L^\infty(\Omega)} \le c \|\varphi_h\|_{L^\infty(\Omega)},\quad \forall \varphi_h\in V_h^0.
\end{equation*} 
\end{lem}

Now we state a quantitative discrete trace estimate, which refines the result in \cite[Lemma 3.2]{CJKZ:2026}. 

\begin{lem}\label{lem:dis_embedding}
With  $\partial_{x_d}\varphi_h$ understood element-wise,
there exists $c>0$, independent of $h$ and $\epsilon$, such that
\begin{equation}\label{eqn:quantitative_discrete_trace}
\|\partial_{x_d}\varphi_h|_{x_d=\ell}\|_{L^2(\omega)}
\le
\frac{c}{\epsilon}
\|(-\Delta_h)^{\frac34+\epsilon}\varphi_h\|_{L^2(\Omega)},
\quad
\forall\varphi_h\in V_h^0,\quad 0<\epsilon\le\frac18.
\end{equation}
\end{lem}

\begin{proof} 
For $\varphi_h\in V_h^0$, let $\psi\in H^2(\Omega)\cap H_0^1(\Omega)$ solve $\Delta\psi=\Delta_h\varphi_h$ in $\Omega$.  
By the definitions of $R_h$ and $-\Delta_h$, we have $R_h\psi=\varphi_h$.
Elliptic regularity and the Ritz projection error estimates then imply
\begin{equation}\label{eqn:psi_H2}
\|\psi\|_{H^2(\Omega)}
\le c\|\Delta_h\varphi_h\|_{L^2(\Omega)}\quad\mbox{and}\quad
\|\psi-\varphi_h\|_{L^2(\Omega)}
+h\|\psi-\varphi_h\|_{H^1(\Omega)}
\le ch^2\|\Delta_h\varphi_h\|_{L^2(\Omega)}.
\end{equation}
The inverse inequality \cite[Lemma 2.2]{JinZhou:2023book} implies
$\|\Delta_h\varphi_h\|_{L^2(\Omega)}
\le ch^{-2}\|\varphi_h\|_{L^2(\Omega)}$. Consequently,
\begin{equation*}
\|\psi\|_{L^2(\Omega)}
\le
\|\varphi_h\|_{L^2(\Omega)}
+\|\psi-\varphi_h\|_{L^2(\Omega)}
\le c\|\varphi_h\|_{L^2(\Omega)}.
\end{equation*}
By splitting $\psi-\varphi_h$ through $\mathcal I_h\psi$,
the element-wise trace inequalities and the interpolation estimate \eqref{eqn:error_I_h},
we obtain
\begin{align*}
\|\partial_{x_d}(\psi-\varphi_h)|_{x_d=\ell}\|_{L^2(\omega)}
&\le
\|\partial_{x_d}(\psi-\mathcal I_h\psi)|_{x_d=\ell}\|_{L^2(\omega)}
+\|\partial_{x_d}(\mathcal I_h\psi-\varphi_h)|_{x_d=\ell}\|_{L^2(\omega)}\\
&\le
ch^{\frac12}\|\psi\|_{H^2(\Omega)}
+ch^{-\frac12}\|\mathcal I_h\psi-\varphi_h\|_{H^1(\Omega)}\le ch^{\frac12}\|\Delta_h\varphi_h\|_{L^2(\Omega)}.
\end{align*}
Moreover, the multiplicative trace inequality, interpolation, and the estimate \eqref{eqn:psi_H2} yield
\begin{align*}
\|\partial_{x_d}\psi|_{x_d=\ell}\|_{L^2(\omega)}\le
c\|\psi\|_{H^1(\Omega)}^{\frac12}
\|\psi\|_{H^2(\Omega)}^{\frac12}\le
c\|\psi\|_{L^2(\Omega)}^{\frac14}
\|\psi\|_{H^2(\Omega)}^{\frac34}\le
c\|\varphi_h\|_{L^2(\Omega)}^{\frac14}
\|\Delta_h\varphi_h\|_{L^2(\Omega)}^{\frac34}.
\end{align*}
Meanwhile, the inverse inequality \cite[Lemma 2.2]{JinZhou:2023book} gives
\begin{equation*}
h^{\frac12}\|\Delta_h\varphi_h\|_{L^2(\Omega)}
\le
c\|\varphi_h\|_{L^2(\Omega)}^{\frac14}
\|\Delta_h\varphi_h\|_{L^2(\Omega)}^{\frac34}.
\end{equation*}
Consequently, we conclude 
\begin{equation}\label{eqn:discrete_multiplicative_trace}
\|\partial_{x_d}\varphi_h|_{x_d=\ell}\|_{L^2(\omega)}
\le
c\|\varphi_h\|_{L^2(\Omega)}^{\frac14}
\|\Delta_h\varphi_h\|_{L^2(\Omega)}^{\frac34}.
\end{equation}
By the Poincar\'e inequality, the smallest eigenvalue of $-\Delta_h$ is bounded below independently of $h$.
Thus, the semigroup smoothing estimates (see, e.g., \cite{Thomee:2006}) and
\eqref{eqn:discrete_multiplicative_trace} imply
\begin{equation*}
\|\partial_{x_d} e^{t\Delta_h} \chi_h\|_{L^2(\omega)}
\le ct^{-\frac34}e^{-ct}\|\chi_h\|_{L^2(\Omega)},
\quad
\forall\chi_h\in V_h^0,\quad t>0.
\end{equation*}
For $\alpha=\frac34+\epsilon$, the semigroup representation
of negative fractional powers gives
\begin{equation*}
(-\Delta_h)^{-\alpha}\chi_h
=
\frac{1}{\Gamma(\alpha)}
\int_0^\infty t^{\alpha-1}e^{t\Delta_h}\chi_h\,\mathrm dt.
\end{equation*}
Consequently, we derive that 
\begin{align*}
\|\partial_{x_d}(-\Delta_h)^{-\alpha}\chi_h\|_{L^2(\omega)}
\le
\frac{c}{\Gamma(\alpha)}
\int_0^\infty t^{\epsilon-1}e^{-ct}\,\mathrm dt\,
\|\chi_h\|_{L^2(\Omega)}\le \frac{c}{\epsilon}\|\chi_h\|_{L^2(\Omega)},
\end{align*}
where $c$ is independent of $h$ and $0<\epsilon\le\frac18$.
Taking $\chi_h=(-\Delta_h)^\alpha\varphi_h$ proves
\eqref{eqn:quantitative_discrete_trace}.
\end{proof}

\subsection{Error estimate for the direct problem}\label{subsec:forward_err} 
In this part, we analyze the error of the fully discrete solution $v_h^n$ defined by \eqref{eqn:fully-dis-v}. We separate the discretization error from the error caused by the noisy boundary data.
In the analysis, we use the discrete harmonic lifting $\widetilde G_h(t)\in V_h$ defined by
\begin{equation}\label{eqn:bdry_lift_dis0}
\left\{
\begin{aligned}
(\nabla\widetilde G_h(t),\nabla\varphi_h)
&=0,
&&\forall\varphi_h\in V_h^0,\\
\widetilde G_h(t)
&=\mathcal I_h^\partial b^\dagger(t),
&&\text{on }\partial\Omega.
\end{aligned}
\right.
\end{equation}
Since $G(t)$ is harmonic, we have
$\widetilde G_h(t)=R_hG(t)$.
Below we write $\widetilde G_h^n=\widetilde G_h(t_n)$.

Let
$\widetilde v_h^0=R_h(\partial_{x_d}u_0)$. Then for any coefficient sequence $q=\{q^n\}_{n=1}^N$,
for $n=1,\dots,N$, we define auxiliary functions $\widetilde v_h^n\in V_h$ by
$\widetilde v_h^n-\widetilde G_h^n\in V_h^0$ and
\begin{equation}\label{eqn:fully-dis-v0}
\bar\partial_\tau
  (\widetilde v_h^n-\widetilde G_h^n)
-\Delta_h(\widetilde v_h^n-\widetilde G_h^n)
=
-P_h^0\bigl(
  \bar\partial_\tau\widetilde G_h^n
  +q^n\widetilde v_h^n
\bigr).
\end{equation}

The next result gives several technical estimates for the lifted boundary data $G(t)$ in \eqref{eqn:bdry_lift}, $G_h(t)$ in \eqref{eqn:bdry_lift_dis} and $\widetilde{G}_h(t)$ in \eqref{eqn:bdry_lift_dis0}.

\begin{lem}\label{lem:err_lift}
Let Assumption \ref{assum:reg_err} hold. Let $G(t)$, $G_h(t)$ and $\widetilde{G}_h(t)$ be the lifted boundary data in  \eqref{eqn:bdry_lift}, \eqref{eqn:bdry_lift_dis} and \eqref{eqn:bdry_lift_dis0}, respectively. Then, for $t\in[0,T]$ and $n=1,\dots,N$, there hold
\begin{align*}
        \|\partial_{x_d}G(t)-\partial_{x_d}\widetilde{G}_h(t)\|_{L^2(\omega)}&\le c h^{\frac{1}{2}},\quad \| \partial_tG(t)- \partial_t\widetilde{G}_h(t)\|_{L^2(\Omega)}\le c  h ,\\
         \| \partial_t\widetilde{G}_h(t_n)- \bar{\partial}_\tau \widetilde{G}_h(t_n)\|_{L^2(\Omega)}&\le c  \tau,\quad
        \| \partial_{x_d}\widetilde{G}_h(t)-\partial_{x_d} G_h(t)\|_{L^2(\omega)}\le ch^{-1}\delta ,\\
        \| \bar{\partial}_\tau \widetilde{G}_h(t_n)-\bar{\partial}_\tau G_h(t_n)\|_{L^2(\Omega)}&\le c\tau^{-1}\delta,\quad
\|\widetilde G_h(t)-G_h(t)\|_{L^2(\Omega)}
\le c\delta.
\end{align*}
The constant $c$ is independent of $h$, $\tau$, $\delta$, $t$, and $n$.
\end{lem}

\begin{proof}
For $\eta_h\in V_h^\partial$, we denote by
$H_h\eta_h\in V_h$ its discrete harmonic extension. We first claim
\begin{equation}\label{eqn:discrete_harmonic_stability}
\|H_h\eta_h\|_{L^2(\Omega)}
+h^{\frac12}\|H_h\eta_h\|_{H^1(\Omega)}
\le c\|\eta_h\|_{L^2(\partial\Omega)}.
\end{equation}
Indeed, a bounded continuous lifting, the stability of
$\mathcal I_h$, the energy estimate, and the inverse estimate on the boundary mesh give
\[
\|H_h\eta_h\|_{H^1(\Omega)}
\le c\|\eta_h\|_{H^{\frac12}(\partial\Omega)}
\le ch^{-\frac12}\|\eta_h\|_{L^2(\partial\Omega)},
\]
To estimate the $L^2(\Omega)$ norm, let
$\psi\in H^2(\Omega)\cap H_0^1(\Omega)$ solve
$-\Delta\psi=H_h\eta_h$. By the elliptic regularity and discrete
harmonicity, we have
\begin{align*}
&\|H_h\eta_h\|_{L^2(\Omega)}^2 =
(\nabla H_h\eta_h,\nabla(\psi-\mathcal I_h\psi))
-\langle\eta_h,\partial_n \psi\rangle_{\partial\Omega}\\
\le&
c\bigl(
h\|H_h\eta_h\|_{H^1(\Omega)}
+\|\eta_h\|_{L^2(\partial\Omega)}
\bigr)\|\psi\|_{H^2(\Omega)}\le
c\|\eta_h\|_{L^2(\partial\Omega)}
 \|H_h\eta_h\|_{L^2(\Omega)}.
\end{align*}
This proves the estimate \eqref{eqn:discrete_harmonic_stability}.
The element-wise trace inequality also yields
\begin{equation}\label{eqn:discrete_harmonic_flux_stability}
\|\partial_{x_d} H_h\eta_h|_{x_d=\ell}\|_{L^2(\omega)}
\le ch^{-1}\|\eta_h\|_{L^2(\partial\Omega)}.
\end{equation}
Since $\widetilde G_h(t)=R_hG(t)$, the standard energy estimate gives
\[
\|G(t)-\widetilde G_h(t)\|_{H^1(\Omega)}
\le ch\|G(t)\|_{H^2(\Omega)}.
\]
By splitting the error $G(t)-\widetilde G_h(t)$ using $\mathcal I_hG(t)$ and 
the element-wise trace inequality, we obtain
\begin{align*}
\|\partial_{x_d}(G(t)-\widetilde G_h(t))
|_{x_d=\ell}\|_{L^2(\omega)}\le
ch^{\frac12}\|G(t)\|_{H^2(\Omega)}
+ch^{-\frac12}
  \|\mathcal I_hG(t)-\widetilde G_h(t)\|_{H^1(\Omega)}\le ch^{\frac12}.
\end{align*}
Next, let $e(t)=\partial_tG(t)-\partial_t\widetilde G_h(t)$. Differentiating the identity \eqref{eqn:bdry_lift_dis0} gives
\begin{equation*}
(\nabla e(t),\nabla\varphi_h)=0
\quad \forall\varphi_h\in V_h^0\quad \text{and}\quad e(t)|_{\partial\Omega}=
\bigl(\partial_tG(t)-\mathcal I_h\partial_tG(t)\bigr)
|_{\partial\Omega}.
\end{equation*}
Then we split $e(t)$ into $e(t)=\eta(t)+\xi_h(t)$, with 
\[
\eta(t)=\partial_tG(t)-\mathcal I_h\partial_tG(t)\quad\mbox{and}
\quad
\xi_h(t)=\mathcal I_h\partial_tG(t)-\partial_t\widetilde G_h(t) \in V_h^0.
\]
Taking $\varphi_h=\xi_h(t)$ in the orthogonality relation yields
\[
\|\nabla\xi_h(t)\|_{L^2(\Omega)}^2
=
-(\nabla\eta(t),\nabla\xi_h(t))
\le
\|\nabla\eta(t)\|_{L^2(\Omega)}
\|\nabla\xi_h(t)\|_{L^2(\Omega)}.
\]
Then the Poincar\'e inequality for $\xi_h(t)$, together with the interpolation estimate, implies
\[
\begin{aligned}
\|e(t)\|_{H^1(\Omega)}
&\le \|\eta(t)\|_{H^1(\Omega)}
+\|\xi_h(t)\|_{H^1(\Omega)} \\
&\le c\|\eta(t)\|_{H^1(\Omega)}
\le ch^{\frac12}\|\partial_tG(t)\|_{H^{\frac32}(\Omega)}.
\end{aligned}
\]
Moreover, since $e(t)=\eta(t)$ on $\partial\Omega$, the local trace inequality on each boundary element $K$, with $h_K=\operatorname{diam}(K)$, gives
\[
\|\eta(t)\|_{L^2(\partial K\cap\partial\Omega)}^2
\le c\left(
h_K^{-1}\|\eta(t)\|_{L^2(K)}^2
+h_K\|\nabla\eta(t)\|_{L^2(K)}^2
\right).
\]
By summing over the boundary elements and the local interpolation estimates, we obtain
\[
\|e(t)\|_{L^2(\partial\Omega)}^2
=
\|\eta(t)\|_{L^2(\partial\Omega)}^2
\le ch^2\|\partial_tG(t)\|_{H^{\frac32}(\Omega)}^2.
\]
To estimate the error in $L^2(\Omega)$, let $\varphi \in H^2(\Omega)\cap H_0^1(\Omega)$ solve the dual problem
$-\Delta \varphi=e(t)$.  We choose $\varphi_h\in V_h^0$ such that $\|\varphi-\varphi_h\|_{H^1(\Omega)}\le ch\|\varphi\|_{H^2(\Omega)}$. Integration by parts and the orthogonality relation then give
\[
\begin{aligned}
\|e(t)\|_{L^2(\Omega)}^2
&=(\nabla e(t),\nabla \varphi)
-\langle e(t),\partial_n \varphi\rangle_{\partial\Omega} =(\nabla e(t),\nabla(\varphi-\varphi_h))
-\langle e(t),\partial_n \varphi\rangle_{\partial\Omega} \\
&\le c\left(
h\|e(t)\|_{H^1(\Omega)}
+\|e(t)\|_{L^2(\partial\Omega)}
\right)\|\varphi\|_{H^2(\Omega)},
\end{aligned}
\]
Hence, combining the preceding estimates and the full elliptic regularity yields 
\[
\begin{aligned}
\|e(t)\|_{L^2(\Omega)}\le c\left(
h\|e(t)\|_{H^1(\Omega)}
+\|e(t)\|_{L^2(\partial\Omega)}
\right) \le c(h^{\frac32}+h)\|\partial_tG(t)\|_{H^{\frac32}(\Omega)}
\le ch\|\partial_tG(t)\|_{H^{\frac32}(\Omega)}.
\end{aligned}
\]
Note that $\|\partial_tG(t)\|_{H^{\frac32}(\Omega)}$ is uniformly bounded by Assumption \ref{assum:reg_err}.
Next, since  
$\widetilde G_h= H_h\mathcal I_h^\partial b^\dagger$,  
the stability estimates \eqref{eqn:discrete_harmonic_stability}
and \eqref{eqn:Ih_bdry_stab} imply
\[
\|\partial_{tt}\widetilde G_h(t)\|_{L^2(\Omega)}
\le c\|\partial_{tt}b^\dagger(t)\|_{L^2(\partial\Omega)}.
\]
Consequently, 
\begin{align*}
\|\partial_t\widetilde G_h(t_n)
-\bar\partial_\tau\widetilde G_h^n\|_{L^2(\Omega)} \le
\frac1\tau\int_{t_{n-1}}^{t_n}
\int_s^{t_n}
\|\partial_{tt}\widetilde G_h(r)\|_{L^2(\Omega)}
\,\mathrm dr\,\mathrm ds\le c\tau.
\end{align*}
Finally, for $\eta_h(t):=\mathcal I_h^\partial(b^\dagger(t)-b^\delta(t))$, by the noise assumption and \eqref{eqn:Ih_bdry_stab}, we have
\[
\|\eta_h(t)\|_{L^2(\partial\Omega)}\le c\delta.
\]
Since $\widetilde G_h(t)-G_h(t)=H_h\eta_h(t)$, the estimates \eqref{eqn:discrete_harmonic_stability} and \eqref{eqn:discrete_harmonic_flux_stability} imply
\[
\|\widetilde G_h(t)-G_h(t)\|_{L^2(\Omega)}
\le c\delta\quad\mbox{and}
\quad
\|\partial_{x_d}(\widetilde G_h(t)-G_h(t))|_{x_d=\ell}\|_{L^2(\omega)}
\le ch^{-1}\delta.
\]
Furthermore, with $\eta_h^n=\eta_h(t_n)$, we obtain
\[
\|\bar\partial_\tau\eta_h^n\|_{L^2(\partial\Omega)}
\le \tau^{-1}\left(
\|\eta_h^n\|_{L^2(\partial\Omega)}
+\|\eta_h^{n-1}\|_{L^2(\partial\Omega)}
\right)
\le c\tau^{-1}\delta.
\]
By the linearity of $H_h$ and the estimate \eqref{eqn:discrete_harmonic_stability}, we deduce
\[
\|\bar\partial_\tau\widetilde G_h(t_n)-\bar\partial_\tau G_h(t_n)\|_{L^2(\Omega)}
=
\|H_h\bar\partial_\tau\eta_h^n\|_{L^2(\Omega)}
\le c\tau^{-1}\delta.
\]
This completes the proof of the lemma.
\end{proof}

Now we derive the  error estimates of the fully discrete scheme \eqref{eqn:fully-dis-v}, which plays an important role in the error analysis. First, we estimate the error between the exact solution $v^\dagger=\partial_{x_d}u(q^\dagger)$ to problem \eqref{eqn:dudx} and the  discrete solution $\widetilde{v}_h^n$ to problem \eqref{eqn:fully-dis-v0}.

\begin{lem}\label{lem:wnh0-wdag}
Let Assumption \ref{assum:reg_err} hold.
Let $v^\dagger=\partial_{x_d}u(q^\dagger)$ and
$\widetilde v_h^n=\widetilde v_h^n(q^\dagger)$ be the
solutions to \eqref{eqn:dudx} and \eqref{eqn:fully-dis-v0},
respectively.
Then there exists $c>0$, independent of $h$, $\tau$, and $n$, such that
\begin{equation*}
\max_{1\le n\le N}
\| \partial_{x_d}(\widetilde v_h^n-v^\dagger(t_n)) \|_{L^2(\omega)}\le c\left(
(1+|\log h|)h^{\frac12}
+(1+|\log\tau|) \tau^{\frac14}
\right).
\end{equation*}
\end{lem}
\begin{proof}
Let $w_0=\partial_{x_d}u_0-G(0)$, $E(t)=e^{t\Delta}$, and $E_h(t)=e^{t\Delta_h}$. 
Note that $p(t):=q^\dagger(t)v^\dagger(t)$ satisfies
\begin{equation}\label{eqn:product_time_regularity}
\|p\|_{W^{1,\infty}(0,T;L^2(\Omega))}\le c.
\end{equation}
Indeed, $w^\dagger=v^\dagger-G$ satisfies
\[
\partial_tw^\dagger-\Delta w^\dagger
=-G'-q^\dagger(w^\dagger+G),
\quad w^\dagger(0)=w_0.
\]
By Assumption \ref{assum:reg_err},
$w_0\in D(-\Delta)$,
$G''\in C([0,T];L^2(\Omega))$, and
$q^\dagger\in W^{1,\infty}(0,T;L^\infty(\omega))$.
Then the standard energy argument and the
time-differentiated equation of $w^\dag$, 
justified by time difference
quotients, together with the initial value 
$\Delta w_0-G'(0)-q^\dagger(0)(w_0+G(0))
\in L^2(\Omega)$ for $\partial_tw^\dagger$, yield
$v^\dagger\in W^{1,\infty}(0,T;L^2(\Omega))$.
The product rule then gives
\eqref{eqn:product_time_regularity}.
Now we define auxiliary solutions
$$\bar v_h(t)=\widetilde G_h(t)+\bar w_h(t)\quad\mbox{and}\quad
\bar v_h^n=\widetilde G_h^n+\bar w_h^n,$$ where
the functions $\bar w_h(t),\bar w_h^n\in V_h^0$ satisfy
\begin{align*}
\partial_t\bar w_h-\Delta_h\bar w_h
&=-P_h^0\bigl(\widetilde G_h'+p\bigr), \quad \text{with}
\quad \bar w_h(0)=R_hw_0,\\
\bar\partial_\tau\bar w_h^n-\Delta_h\bar w_h^n
&=-P_h^0\bigl(
\bar\partial_\tau\widetilde G_h^n+p(t_n)\bigr),
 \quad \text{with}
\quad\bar w_h^0=R_hw_0.
\end{align*}
Then we can split the error as
$$v^\dagger(t_n)-\widetilde v_h^n =(v^\dagger(t_n)-\bar v_h(t_n)) 
+ (\bar v_h(t_n)-\bar v_h^n)
+ (\bar v_h^n-\widetilde v_h^n)
=: \mathrm{I}^n + \mathrm{II}^n + \theta_h^n.$$
Below we bound the three terms separately. The lengthy proof is divided into three steps.

\smallskip
\noindent
\textbf{Step 1. Spatial discretization error $\mathrm{I}^n$.}
The solution representations
give
\begin{align*}
\mathrm I^n
={}&(G(t_n)-\widetilde G_h(t_n))
+(E(t_n)w_0-E_h(t_n)R_hw_0)\\
&+\int_0^{t_n}(E_h(t_n-s)P_h^0-E(t_n-s))
\bigl(G'(s)+p(s)\bigr)\,\mathrm ds \\
&+\int_0^{t_n}E_h(t_n-s)P_h^0
\bigl(\widetilde G_h'(s)-G'(s)\bigr)\,\mathrm ds,
\end{align*}
The first term can be directly bounded using Lemma \ref{lem:err_lift}. Now we bound the second term. 
For any $f\in L^2(\Omega)$ and a fixed $\theta\in(\pi/2,\pi)$, let $u=(z-\Delta)^{-1}f$ and
$u_h=(z-\Delta_h)^{-1}P_h^0f$. Then the standard finite element resolvent estimates \cite[Lemma 2.6]{JinZhou:2023book} give
\begin{equation*}
\|u\|_{H^2(\Omega)}\le c\|f\|_{L^2(\Omega)}\quad\text{and}
\quad
\|u-u_h\|_{H^1(\Omega)}
\le ch\|f\|_{L^2(\Omega)}.
\end{equation*}
By splitting $u-u_h$ through $\mathcal I_hu$ and using the estimates
\eqref{eqn:dis_trace}, \eqref{eqn:cts_trace_ineq}, and
\eqref{eqn:error_I_h}, we obtain
\begin{align*}
\|\partial_{x_d}(u-u_h)|_{x_d=\ell}\|_{L^2(\omega)}\le
ch^{\frac12}\|u\|_{H^2(\Omega)}
+ch^{-\frac12}\|\mathcal I_hu-u_h\|_{H^1(\Omega)}\le ch^{\frac12}\|f\|_{L^2(\Omega)}.
\end{align*}
This proves the boundary resolvent error estimate \begin{equation}\label{eqn:flux_resolvent_error}
\|
\partial_{x_d}
[((z-\Delta_h)^{-1}P_h^0
-(z-\Delta)^{-1})f]|_{x_d=\ell}\|_{L^2(\omega)}
\le ch^{\frac12}\|f\|_{L^2(\Omega)},
\quad \forall z\in\Sigma_\theta,~~ f\in L^2(\Omega).
\end{equation}
Meanwhile, using the relation $(-\Delta_h)R_h=P_h^0(-\Delta)$,
we have 
\begin{align*}
E(t)w_0-E_h(t)R_hw_0
=w_0-R_hw_0-\frac{1}{2\pi\mathrm i}
\int_{\Gamma_{\theta,1/t}}
e^{zt}z^{-1}
\left((z-\Delta)^{-1}
-(z-\Delta_h)^{-1}P_h^0\right)
(-\Delta)w_0\,\mathrm dz.
\end{align*}
Thus, the element-wise trace
inequalities  
\eqref{eqn:dis_trace} and \eqref{eqn:cts_trace_ineq}, the estimates  \eqref{eqn:error_R_h}, \eqref{eqn:error_I_h} and \eqref{eqn:flux_resolvent_error} imply
\begin{align*}
\|\partial_{x_d}
\bigl(E(t_n)w_0-E_h(t_n)R_hw_0\bigr)
|_{x_d=\ell}\|_{L^2(\omega)}
\le ch^{\frac12}\|(-\Delta)w_0\|_{L^2(\Omega)}
\le ch^{\frac12}.
\end{align*}
Similarly, the contour representation and
\eqref{eqn:flux_resolvent_error} also yield
\[
\left\| 
\partial_{x_d}\bigl(E_h(t)P_h^0-E(t)\bigr)f |_{x_d=\ell}
\right\|_{L^2(\omega)}
\le ch^{\frac12}t^{-1}\|f\|_{L^2(\Omega)}.
\]
Meanwhile, the trace estimate, semigroup smoothing, and
\eqref{eqn:quantitative_discrete_trace} with
$\epsilon=\frac18$ give
\[
\left\| 
\partial_{x_d}\bigl(E_h(t)P_h^0-E(t)\bigr)f |_{x_d=\ell}
\right\|_{L^2(\omega)}
\le ct^{-\frac78}\|f\|_{L^2(\Omega)},
\quad 0<t\le T.
\]
Consequently, we conclude 
\begin{equation}\label{eqn:flux_semigroup_error}
\left\| 
\partial_{x_d}\bigl(E_h(t)P_h^0-E(t)\bigr)f |_{x_d=\ell}
\right\|_{L^2(\omega)}
\le
c\min\{h^{\frac12}t^{-1},t^{-\frac78}\}
\|f\|_{L^2(\Omega)}.
\end{equation}
Since the function $G'+p$ is bounded in
$L^\infty(0,T;L^2(\Omega))$, splitting the time integral
at $h^4$ gives
\begin{align*}
&\left\|
\partial_{x_d}
\int_0^{t_n}(E_h(t_n-s)P_h^0-E(t_n-s))
\bigl(G'(s)+p(s)\bigr)\,\mathrm ds |_{x_d=\ell}
\right\|_{L^2(\omega)}\\
\le&
c\int_0^{t_n}
\min\{h^{\frac12}r^{-1},r^{-\frac78}\}\,\mathrm dr
\le ch^{\frac12}(1+|\log h|).
\end{align*}
For the remaining convolution, Lemma \ref{lem:err_lift}
and discrete semigroup smoothing imply
\begin{align*}
&\left\|
\left.
\partial_{x_d}
\left[
\int_0^{t_n}E_h(t_n-s)P_h^0
\bigl(\widetilde G_h'(s)-G'(s)\bigr)\,\mathrm ds
\right]
\right|_{x_d=\ell}
\right\|_{L^2(\omega)}\\
\le&
c\int_0^{t_n}(t_n-s)^{-\frac78}
\|\widetilde G_h'(s)-G'(s)\|_{L^2(\Omega)}
\,\mathrm ds
\le ch.
\end{align*}
Combining the preceding estimates yields
\begin{equation}\label{eqn:spatial_flux_error_exact}
\left\|
\left.\partial_{x_d}\mathrm I^n\right|_{x_d=\ell}
\right\|_{L^2(\omega)}
\le ch^{\frac12}(1+|\log h|).
\end{equation}
Repeating the argument also provides the estimate in $L^2(\Omega)$ for the term $\mathrm I^n$:
\begin{equation}\label{eqn:spatial_bulk_error_exact}
\|\mathrm I^n\|_{L^2(\Omega)}\le ch.
\end{equation}

\smallskip
\noindent
\textbf{Step 2. Temporal discretization error $\mathrm{II}^n$.}
The solution representation \eqref{eqn:dis-solrep-w},
on the auxiliary problem, give
\begin{align*}
\mathrm{II}^n
={}&
\bigl(E_h(t_n)-E_{h,\tau}^{n-1}\bigr)R_hw_0
+\tau\sum_{j=1}^n E_{h,\tau}^{n-j}P_h^0
\bigl(\bar\partial_\tau\widetilde G_h^j
-\widetilde G_h'(t_j)\bigr)\\
&+\tau\sum_{j=1}^n E_{h,\tau}^{n-j}P_h^0\bigl(\widetilde G_h'+p\bigr)(t_j)
-\int_0^{t_n}E_h(t_n-s)P_h^0\bigl(\widetilde G_h'(s)+p(s)\bigr)\,\mathrm ds.
\end{align*}
By the identity
$\widetilde G_h=H_h\mathcal I_h^\partial b^\dagger$,
and the estimates \eqref{eqn:discrete_harmonic_stability},
\eqref{eqn:Ih_bdry_stab}, and
\eqref{eqn:product_time_regularity}, we obtain
\begin{align*}
\|P_h^0\bigl(\widetilde G_h'(0)+p(0)\bigr)\|_{L^2(\Omega)}
+\|P_h^0\bigl(\widetilde G_h''(t)+p'(t)\bigr)\|_{L^\infty(0,T;L^2(\Omega))}
\le c.
\end{align*}
Let $\alpha=\frac34+\epsilon$, where
$0<\epsilon\le\frac18$.
Since $(-\Delta_h)R_hw_0=P_h^0(-\Delta)w_0$,
the standard backward Euler error estimate (see \cite[Chapter~7]{Thomee:2006}) and the stability of $P_h^0$
give
\begin{align*}
\left\|
(-\Delta_h)^\alpha
\bigl(E_h(t_n)-E_{h,\tau}^{n-1}\bigr)R_hw_0
\right\|_{L^2(\Omega)}\le
c\tau t_n^{-\alpha}
\|(-\Delta_h)R_hw_0\|_{L^2(\Omega)}\le
c\tau t_n^{-\alpha}\|\Delta w_0\|_{L^2(\Omega)}
\le c\tau t_n^{-\alpha}.
\end{align*}
By the discrete smoothing property in  Lemma \ref{lem:dis-sol-op} and Lemma \ref{lem:err_lift}, we derive
\[
\big\|
(-\Delta_h)^\alpha
\tau\sum_{j=1}^n E_{h,\tau}^{n-j}P_h^0
\bigl(\bar\partial_\tau\widetilde G_h^j
-\widetilde G_h'(t_j)\bigr)
\big\|_{L^2(\Omega)}
\le
c\tau\sum_{j=1}^n
t_{n-j+1}^{-\alpha}\|\bigl(\bar\partial_\tau\widetilde G_h^j
-\widetilde G_h'(t_j)\bigr)\|_{L^2(\Omega)}
\le c\tau.
\]
Further, the backward Euler convolution error estimate (by the standard resolvent argument \cite[Theorem 3.4]{JinZhou:2023book}) yields
\begin{align*}
&\big\|
(-\Delta_h)^\alpha
\big[
\tau\sum_{j=1}^n E_{h,\tau}^{n-j}P_h^0\bigl(\widetilde G_h'+p\bigr)(t_j)
-\int_0^{t_n}E_h(t_n-s)P_h^0\bigl(\widetilde G_h'(s)+p(s)\bigr)\,\mathrm ds
\big]
\big\|_{L^2(\Omega)}\\
\le&
c\tau t_n^{-\alpha}\|(\widetilde G_h'+p\bigr)(0)\|_{L^2(\Omega)}
+c\tau\int_0^{t_n}(t_n-s)^{-\alpha}
\|(\widetilde G_h'+p\bigr)'(s)\|_{L^2(\Omega)}\,\mathrm ds
\le c\tau(1+t_n^{-\alpha}).
\end{align*} 
Since $\mathrm{II}^n\in V_h^0$, combining these estimates
with \eqref{eqn:quantitative_discrete_trace} gives
\begin{align*}
\left\|
\left.\partial_{x_d}\mathrm{II}^n\right|_{x_d=\ell}
\right\|_{L^2(\omega)}\le
\frac{c}{\epsilon}
\|(-\Delta_h)^{\frac34+\epsilon}\mathrm{II}^n\|_{L^2(\Omega)}\le
\frac{c}{\epsilon}\tau
\bigl(1+t_n^{-\frac34-\epsilon}\bigr)
\le
\frac{c}{\epsilon}\tau^{\frac14-\epsilon},
\end{align*}
where we used $t_n\ge\tau$ and $0<\tau\le1$.  By choosing
$\epsilon =\min\left\{\frac18,\frac1{1+|\log\tau|}\right\}$,
we have
$\epsilon^{-1}\le c(1+|\log\tau|)$ and
$\tau^{-\epsilon}\le e$.
Consequently,
\begin{equation}\label{eqn:temporal_flux_error_exact}
\max_{1\le n\le N}
\left\|
\left.\partial_{x_d}\mathrm{II}^n\right|_{x_d=\ell}
\right\|_{L^2(\omega)}
\le c\tau^{\frac14}(1+|\log\tau|).
\end{equation}
The backward Euler estimate also gives
\begin{equation}\label{eqn:temporal_bulk_error_exact}
\max_{1\le n\le N}
\|\mathrm{II}^n\|_{L^2(\Omega)}
\le c\tau.
\end{equation}

\smallskip
\noindent
\textbf{Step 3. Coefficient perturbation.}
Set
$\eta^n=v^\dagger(t_n)-\bar v_h^n
=\mathrm I^n+\mathrm{II}^n$.
Steps 1 and 2 provide the estimate
\begin{equation}\label{eqn:auxiliary_direct_error}
\begin{aligned}
\|v^\dagger(t_n)-\bar v_h^n\|_{L^2(\Omega)}
&\le c(h+\tau).
\end{aligned}
\end{equation}
Subtracting \eqref{eqn:fully-dis-v0} from the auxiliary
discrete equation yields
\[
\bar\partial_\tau\theta_h^n-\Delta_h\theta_h^n
+P_h^0\bigl(q^\dagger(t_n)\theta_h^n\bigr)
=P_h^0\bigl(q^\dagger(t_n)(\bar v_h^n-v^\dagger(t_n))\bigr),
\quad \theta_h^0=0.
\]
Testing with $\theta_h^n\in V_h^0$, and applying the estimate \eqref{eqn:auxiliary_direct_error} and 
the discrete Gronwall inequality gives
\[
\max_{1\le j\le n}\|\theta_h^j\|_{L^2(\Omega)}
\le
c\max_{1\le j\le n}\|\eta^j\|_{L^2(\Omega)}
\le c(h+\tau).
\]
The discrete solution representation further yields
\[
\theta_h^n
=
-\tau\sum_{j=1}^n E_{h,\tau}^{n-j}P_h^0
\left[
q^\dagger(t_j)(\eta^j+\theta_h^j)
\right].
\]
By the estimate \eqref{eqn:quantitative_discrete_trace}
with $\epsilon=\frac18$ and
Lemma \ref{lem:dis-sol-op}, we obtain
\begin{align*}
\left\|
\left.\partial_{x_d}\theta_h^n\right|_{x_d=\ell}\right\|_{L^2(\omega)}\le
c\tau\sum_{j=1}^n t_{n-j+1}^{-\frac78}
\left(
\|\eta^j\|_{L^2(\Omega)}
+\|\theta_h^j\|_{L^2(\Omega)}
\right)\le c(h+\tau).
\end{align*}
Finally,  the triangle inequality yields the desired result.  
\end{proof}

Under additional initial compatibility conditions, the temporal
error estimate in Lemma \ref{lem:wnh0-wdag} can be improved to first order,
as stated in the following corollary.
\begin{cor} 
\label{cor:wnh0-wdag-first-order}
Let the assumptions of Lemma \ref{lem:wnh0-wdag} hold.
Assume additionally that
\begin{equation}\label{eqn:initial_compatibility_first_order}
w_0=\partial_{x_d}u_0-G(0)=0,
\quad
\partial_t b^\dagger(0)=0,
\quad
q^\dagger(0)v^\dagger(0)=0.
\end{equation}
Then there exists a constant $c>0$, independent of $h$,
$\tau$, and $n$, such that
\[
\max_{1\le n\le N}
\left\|
\left.
\partial_{x_d}
\bigl(\widetilde v_h^n-v^\dagger(t_n)\bigr)
\right|_{x_d=\ell}
\right\|_{L^2(\omega)}
\le
c\left((1+|\log h|)h^{\frac12}+\tau\right).
\]
\end{cor}
\begin{proof}
We use the notation and error decomposition in the proof
of Lemma \ref{lem:wnh0-wdag}. Only the estimate for the
temporal error $\mathrm{II}^n$ needs to be improved. Fix $\alpha=7/8$.
Since $\widetilde G_h=H_h\mathcal I_h^\partial b^\dagger$,
the compatibility conditions
\eqref{eqn:initial_compatibility_first_order} imply $P_h^0\bigl(\widetilde G_h'(0)+p(0)\bigr)=0$. Then the backward Euler convolution error estimate (by the standard resolvent argument
\cite[Theorem 3.4]{JinZhou:2023book}) yields
\begin{align*}
&\big\|
(-\Delta_h)^\alpha
\big[
\tau\sum_{j=1}^n E_{h,\tau}^{n-j}P_h^0\bigl(\widetilde G_h'+p\bigr)(t_j)
-\int_0^{t_n}E_h(t_n-s)P_h^0\bigl(\widetilde G_h'(s)+p(s)\bigr)\,\mathrm ds
\big]
\big\|_{L^2(\Omega)}\\
\le& c\tau\int_0^{t_n}(t_n-s)^{-\alpha}
\|(\widetilde G_h'+p\bigr)'(s)\|_{L^2(\Omega)}\,\mathrm ds \le c\tau.
\end{align*} 
Consequently, we have
$\max_{1\le n\le N}
\|(-\Delta_h)^\alpha\mathrm{II}^n\|_{L^2(\Omega)}
\le c\tau$.
By \eqref{eqn:quantitative_discrete_trace}
with the fixed choice $\epsilon=1/8$, we obtain
\begin{align*}
\max_{1\le n\le N}
\left\|
\left.
\partial_{x_d}\mathrm{II}^n
\right|_{x_d=\ell}
\right\|_{L^2(\omega)}
\le c\tau.
\end{align*} 
The spatial error estimate from Step 1 and the coefficient
perturbation estimate from Step 3 remain unchanged. Then the proof is complete. 
\end{proof}

Next, we estimate the error due to data noise, namely, the
difference between the discrete solutions $\widetilde v_h^n$
and $v_h^n$ to problems \eqref{eqn:fully-dis-v0} and
\eqref{eqn:fully-dis-v}, respectively.

\begin{lem}\label{lem:wnh-wnh0}
Let Assumption \ref{assum:reg_err} hold. Let $q^\dagger$ be the exact coefficient, $v_h^n=v_h^n(q^\dagger)$ be the discrete solution to problem \eqref{eqn:fully-dis-v} and  $\widetilde{v}_h^n=\widetilde{v}_h^n(q^\dagger)$ be the discrete solution to problem \eqref{eqn:fully-dis-v0}.   Then there exists $c>0$   independent of $\delta,\, h,\,\tau$ and $t_n$, such that
    \begin{equation*}
        \|\partial_{x_d} v_h^n-\partial_{x_d} \widetilde{v}_h^n \|_{L^2(\omega)}\le c(h^{ -1}\delta+ \tau^{-1}\delta ).
    \end{equation*}  
\end{lem}
\begin{proof}
Since the noisy data agree with the exact initial data at
$t=0$, we have $G_h^0=\widetilde G_h^0$.
Subtracting the discrete solution representations yields
\begin{align*}
v_h^n-\widetilde v_h^n=G_h^n-\widetilde G_h^n
-\tau\sum_{j=1}^{n}E_{h,\tau}^{n-j}P_h^0
\bigl(\bar\partial_\tau G_h^j
-\bar\partial_\tau\widetilde G_h^j\bigr)
-\tau\sum_{j=1}^{n}E_{h,\tau}^{n-j}P_h^0
\bigl((q^\dagger)^j(v_h^j-\widetilde v_h^j)\bigr).
\end{align*}
By Lemmas \ref{lem:dis-sol-op} and \ref{lem:err_lift},
the $L^2$ stability of $P_h^0$, and the boundedness
of $q^\dagger$, we obtain
\begin{align*}
\|v_h^n-\widetilde v_h^n\|_{L^2(\Omega)}
&\le
c\delta
+c\tau\sum_{j=1}^{n}\tau^{-1}\delta
+c\tau\sum_{j=1}^{n}
\|v_h^j-\widetilde v_h^j\|_{L^2(\Omega)}
\\
&\le
c\tau^{-1}\delta
+c\tau\sum_{j=1}^{n}
\|v_h^j-\widetilde v_h^j\|_{L^2(\Omega)},
\end{align*}
The discrete Gronwall inequality then gives
\begin{equation*}
\max_{1\le n\le N}
\|v_h^n-\widetilde v_h^n\|_{L^2(\Omega)}
\le c\tau^{-1}\delta.
\end{equation*}
Next, by \eqref{eqn:quantitative_discrete_trace}
with $\epsilon=1/8$, together with the discrete smoothing
estimate from Lemma \ref{lem:dis-sol-op} and the lifting
estimates from Lemma \ref{lem:err_lift}, we derive
\begin{align*}
\left\|
\partial_{x_d}(v_h^n-\widetilde v_h^n)
|_{x_d=\ell}
\right\|_{L^2(\omega)}\le ch^{-1}\delta
+c\tau\sum_{j=1}^{n}(t_{n+1}-t_j)^{-\frac78}\tau^{-1}\delta \le
ch^{-1}\delta
+c\tau^{-1}\delta.
\end{align*}
This completes the proof of the lemma.
\end{proof}

\subsection{Error estimate for (ICP2)}\label{subsec:inverse_err}

We first introduce weighted discrete norms. For $\lambda>0$ and $p\in[1,\infty]$, we define
\begin{equation}\label{eqn:dis_weighted_norm}
    \|v\|_{\ell_\lambda^p}
    =
    \begin{cases}
        \displaystyle \left(\tau\sum_{n=1}^{\infty}|e^{-\lambda t_n}v^n|^p\right)^{1/p}, & 1\le p<\infty,\\[2mm]
        \displaystyle \sup_{n\ge1}|e^{-\lambda t_n}v^n|, & p=\infty.
    \end{cases}
\end{equation}
For a finite sequence $(v^n)_{n=1}^{N}$, the norm is defined by extending the sequence by zero for $n>N$. We denote by $\|\cdot\|_{\ell^p}$ the unweighted norm (i.e., $\lambda=0$). For finite sequences with $t_N=T$, we have
\[
    e^{-\lambda T}\|v\|_{\ell^p}\le \|v\|_{\ell_\lambda^p}\le \|v\|_{\ell^p}, \quad p\in[1,\infty].
\]
Thus, the weighted and unweighted norms are equivalent, with equivalence constants independent of $\tau$ and $N$.
For a Banach space $X$, we define $\|v\|_{\ell_\lambda^p(X)}$ analogously by replacing $|v^n|$ with $\|v^n\|_X$ in \eqref{eqn:dis_weighted_norm}. The  zero-extension convention and norm equivalence apply to finite $X$-valued sequences.

We first establish a uniform $L^\infty$ bound for the fully discrete solutions.

\begin{lem}\label{lem:bound_w}
Let Assumption \ref{assum:reg_err} hold. Let $q=(q^n)_{n=1}^{N}$ satisfy $|q^n|\le b$ a.e. in $\omega$, and let $v_h^n=v_h^n(q)$ be the solutions to \eqref{eqn:fully-dis-v}. Then there exist constants $\widetilde{\tau}>0$ and $c>0$, independent of $\delta$, $h$, $\tau$, $N$, and $q$, such that, for all $0<\tau\le\widetilde{\tau}$,
\begin{equation*}
\|(v_h^n)_{n=1}^{N} \|_{\ell^\infty(L^\infty(\Omega))} \le c.
\end{equation*} 
\end{lem}

\begin{proof} 
By the solution representation \eqref{eqn:dis-solrep-w}, we have
\begin{align*}
v_h^n =G_h^n+E_{h,\tau}^{n-1}R_h(\partial_{x_d}u_0-G(0))-\tau\sum_{j=1}^{n}E_{h,\tau}^{n-j}P_h^0(\bar{\partial}_\tau G_h^j)
    -\tau\sum_{j=1}^{n}E_{h,\tau}^{n-j}P_h^0(q^jv_h^j)
    =\sum_{i=1}^{4}\mathrm{I}_i^n.
\end{align*}
The discrete elliptic stability \cite[Theorem 1.1]{LeykekhmanLi:2021} and the boundary interpolation stability \eqref{eqn:Ih_bdry_stab} yield
\begin{align*}
\|\mathrm I_1^n\|_{L^\infty(\Omega)}
&\le c\|\mathcal I_h^\partial b^\delta(t_n)\|_{L^\infty(\partial\Omega)}
\le c\|b^\delta(t_n)\|_{L^\infty(\partial\Omega)}\\
&\le c\left(\|z^\delta(t_n)\|_{L^\infty(\omega)}+\|\partial_{x_d}g(t_n)\|_{L^\infty(\partial\omega\times(0,\ell))}\right).
\end{align*}
Next, the inverse inequality \eqref{eqn:inverse_ineq}, the $L^\infty$ stability of $P_h$ \cite[Lemma 6.1]{Thomee:2006}, and the approximation estimates \eqref{eqn:error_P_h}--\eqref{eqn:error_R_h} imply
\begin{align*}
    \|R_h(\partial_{x_d}u_0-G(0))\|_{L^\infty(\Omega)}
    &\le \|(R_h-P_h)(\partial_{x_d}u_0-G(0))\|_{L^\infty(\Omega)}
    +\|P_h(\partial_{x_d}u_0-G(0))\|_{L^\infty(\Omega)}\\
    &\le ch^{-d/2}\|(R_h-P_h)(\partial_{x_d}u_0-G(0))\|_{L^2(\Omega)}
    +c\|\partial_{x_d}u_0-G(0)\|_{L^\infty(\Omega)}\\
    &\le ch^{2-d/2}\|\partial_{x_d}u_0-G(0)\|_{H^2(\Omega)}
    +c\|\partial_{x_d}u_0-G(0)\|_{L^\infty(\Omega)}.
\end{align*}
Hence, Lemma \ref{lem:dis-sol-op} gives
\begin{align*}
\|\mathrm{I}_2^n\|_{L^\infty(\Omega)}
\le c\|R_h(\partial_{x_d}u_0-G(0))\|_{L^\infty(\Omega)}\le c\left(
\|\partial_{x_d}u_0-G(0)\|_{H^2(\Omega)} +\|\partial_{x_d}u_0-G(0)\|_{L^\infty(\Omega)}
\right).
\end{align*}
Since the discrete lifting operator is linear and independent of time, $G_h'(s)$ is the discrete lifting of the time derivatives of the boundary data for almost every $s\in(0,T)$. the relation $\bar{\partial}_\tau G_h^j=\frac{1}{\tau}\int_{t_{j-1}}^{t_j}G_h'(s)\mathrm{d} s $, Lemma \ref{lem:dis-sol-op}, the $L^\infty$ stability of $P_h^0$, the discrete elliptic stability \cite[Theorem 1.1]{LeykekhmanLi:2021}, and \eqref{eqn:Ih_bdry_stab}, we obtain
\begin{align*}
    \|\mathrm{I}_3^n\|_{L^\infty(\Omega)}
    &\le c\tau\sum_{j=1}^{n}\|P_h^0(\bar{\partial}_\tau G_h^j)\|_{L^\infty(\Omega)}
    \le c\tau\sum_{j=1}^{n}\|\bar{\partial}_\tau G_h^j\|_{L^\infty(\Omega)}\\
    &\le c\int_0^{t_n}\|G_h'(s)\|_{L^\infty(\Omega)}\,\mathrm{d}s\le c\int_0^{t_n}\left(
\|\partial_tz^\delta(s)\|_{L^\infty(\omega)}  +\|\partial_t\partial_{x_d}g(s)\|_{L^\infty(\partial\omega\times(0,\ell))}
    \right)\,\mathrm{d}s\\
    &\le c\left(
    \|z^\delta\|_{W^{1,1}(0,T;L^\infty(\omega))}
    +\|\partial_{x_d}g\|_{W^{1,1}(0,T;L^\infty(\partial\omega\times(0,\ell)))}
    \right).
\end{align*}
Moreover, Lemma \ref{lem:dis-sol-op}, the $L^\infty$ stability of $P_h^0$, and the box constraint imply
\begin{align*}
\|\mathrm{I}_4^n\|_{L^\infty(\Omega)}\le c\tau\sum_{j=1}^{n}\|P_h^0(qv_h^j)\|_{L^\infty(\Omega)}\le c\tau\sum_{j=1}^{n}\|qv_h^j\|_{L^\infty(\Omega)}\le c\tau\sum_{j=1}^{n}\|v_h^j\|_{L^\infty(\Omega)}.
\end{align*}
Combining these estimates and using Assumption \ref{assum:reg_err}, together with the continuous embedding of $W^{1,1}(0,T;X)$ into $C([0,T];X)$ for the boundary data spaces, we obtain
\begin{align*}
\|v_h^n\|_{L^\infty(\Omega)}    \le c+c\tau\sum_{j=1}^{n}\|v_h^j\|_{L^\infty(\Omega)},\quad 1\le n\le N,
\end{align*}
where $c$ is independent of $h$, $\tau$, $N$, $n$, and $q$. For sufficiently small $\tau>0$, the discrete Gr\"onwall inequality therefore gives the desired estimate and completes the proof.
\end{proof}

The next lemma establishes the stability of the boundary derivative of the discrete solution in the weighted $\ell^\infty_\lambda$ norm.

\begin{lem}\label{lem:dis_forward_stab}
Let Assumption \ref{assum:reg_err} hold. For $i=1,2$, let $q_i=(q_i^n)_{n=1}^{N}$ satisfy $|q_i^n|\le b$ a.e. in $\omega$, and let $v_{i,h}^n=v_h^n(q_i)$ be the solutions to \eqref{eqn:fully-dis-v}. For any fixed $\epsilon\in(0,\frac18]$, there exist constants $\lambda_0>0$, $\widetilde{\tau}>0$, and $c>0$, independent of $h$, $\tau$, $N$, $\lambda$, $q_1$, and $q_2$, such that
\begin{equation*}
    \left\|
    \left(\left.\partial_{x_d}(v_{1,h}^n-v_{2,h}^n)\right|_{x_d=\ell}\right)_{n=1}^{N}
\right\|_{\ell^\infty_\lambda(L^2(\omega))}
    \le c\bigl(\tau^{\frac14-\epsilon}+\lambda^{-\frac14+\epsilon}\bigr)
\|(q_1^n-q_2^n)_{n=1}^{N} \|_{\ell^\infty_\lambda(L^2(\omega))}
\end{equation*}
for all $\lambda\ge\lambda_0$ and $0<\tau\le\widetilde{\tau}$. The constants may depend on $\epsilon$.
\end{lem}

\begin{proof}
Let $w_h^n=v_{1,h}^n-v_{2,h}^n$. Since the initial and boundary data are the same for both solutions, the solution representation \eqref{eqn:dis-solrep-w} gives
\begin{equation*}
    w_h^n
    =\tau\sum_{j=1}^{n}E_{h,\tau}^{n-j}P_h^0
    \bigl(-q_2^jw_h^j+(q_2^j-q_1^j)v_{1,h}^j\bigr).
\end{equation*}
By Lemma \ref{lem:dis-sol-op}, the $L^2$ stability of $P_h^0$, the box constraint on $q_2$, and Lemma \ref{lem:bound_w}, we obtain
\begin{align*}
\|w_h^n\|_{L^2(\Omega)}\le c\tau\sum_{j=1}^{n}
    \|-q_2^jw_h^j+(q_2^j-q_1^j)v_{1,h}^j\|_{L^2(\Omega)}\le c\tau\sum_{j=1}^{n}
    \left(
        \|w_h^j\|_{L^2(\Omega)}
        +\|q_2^j-q_1^j\|_{L^2(\omega)}
    \right).
\end{align*}
Multiplying this estimate by $e^{-\lambda t_n}$ and taking the maximum over $1\le n\le N$ yields
\begin{align*}
\|(w_h^n)_{n=1}^{N}\| _{\ell^\infty_\lambda(L^2(\Omega))} 
    &\le c\left(\tau\sum_{k=1}^{N}e^{-\lambda t_{k-1}}\right)
    \left(
        \|(w_h^n)_{n=1}^{N}\|_{\ell^\infty_\lambda(L^2(\Omega))}
        +\|(q_1^n-q_2^n)_{n=1}^{N}\|_{\ell^\infty_\lambda(L^2(\omega))}
    \right)\\
    &\le c(\tau+\lambda^{-1})
    \left(
        \|(w_h^n)_{n=1}^{N}\|_{\ell^\infty_\lambda(L^2(\Omega))}
        +\|(q_1^n-q_2^n)_{n=1}^{N}\|_{\ell^\infty_\lambda(L^2(\omega))}
    \right),
\end{align*}
where we used
\begin{equation*}
    \tau\sum_{k=1}^{N}e^{-\lambda t_{k-1}}
    \le \tau+\int_0^{t_{N-1}}e^{-\lambda s}\,\mathrm{d}s
    \le \tau+\lambda^{-1}.
\end{equation*}
Now we choose $\lambda_0\ge1$ sufficiently large and $\widetilde{\tau}\le1$ sufficiently small such that $c(\widetilde{\tau}+\lambda_0^{-1})\le1/2$ and Lemma \ref{lem:bound_w} applies for $0<\tau\le\widetilde{\tau}$. Then, for all $\lambda\ge\lambda_0$ and $0<\tau\le\widetilde{\tau}$, absorption gives
\begin{equation*}
    \|(w_h^n)_{n=1}^{N}\|_{\ell^\infty_\lambda(L^2(\Omega))}
    \le c(\tau+\lambda^{-1})
    \|(q_1^n-q_2^n)_{n=1}^{N}\|_{\ell^\infty_\lambda(L^2(\omega))}.
\end{equation*}
Next, set $\beta=\frac34+\epsilon\in(0,1)$. Lemmas \ref{lem:dis-sol-op} and \ref{lem:dis_embedding}, together with the solution representation and the preceding bounds, imply
\begin{equation*}
    \left\|\left.\partial_{x_d}w_h^n\right|_{x_d=\ell}\right\|_{L^2(\omega)}
    \le c\tau\sum_{j=1}^{n}t_{n-j+1}^{-\beta}
    \left(
        \|w_h^j\|_{L^2(\Omega)}
        +\|q_2^j-q_1^j\|_{L^2(\omega)}
    \right).
\end{equation*}
By multiplying by $e^{-\lambda t_n}$ and taking the maximum over $1\le n\le N$, we obtain
\begin{align*}
&\left\|
\left(\left.\partial_{x_d}w_h^n\right|_{x_d=\ell}\right)_{n=1}^{N}
    \right\|_{\ell^\infty_\lambda(L^2(\omega))}\\
\le& c\left(
        \tau\sum_{k=1}^{N}e^{-\lambda t_{k-1}}t_k^{-\beta}
    \right)
    \left(
        \|(w_h^n)_{n=1}^{N}\|_{\ell^\infty_\lambda(L^2(\Omega))}
        +\|(q_1^n-q_2^n)_{n=1}^{N}\|_{\ell^\infty_\lambda(L^2(\omega))}
    \right)\\
\le& c\left(
        \tau\sum_{k=1}^{N}e^{-\lambda t_{k-1}}t_k^{-\beta}
    \right)(1+\tau+\lambda^{-1})
    \|(q_1^n-q_2^n)_{n=1}^{N}\|_{\ell^\infty_\lambda(L^2(\omega))}.
\end{align*}
Since $0<\beta<1$, we have
\begin{align*}
    \tau\sum_{k=1}^{N}e^{-\lambda t_{k-1}}t_k^{-\beta}
    \le \tau^{1-\beta}
    +\int_0^{t_{N-1}}e^{-\lambda s}s^{-\beta}\,\mathrm{d}s\le \tau^{1-\beta}+\Gamma(1-\beta)\lambda^{\beta-1}\le c\bigl(\tau^{\frac14-\epsilon}+\lambda^{-\frac14+\epsilon}\bigr).
\end{align*}
Combining these estimates and using $\tau\le1$ and $\lambda\ge1$ gives
\begin{equation*}
    \left\|
    \left(\left.\partial_{x_d}w_h^n\right|_{x_d=\ell}\right)_{n=1}^{N}
    \right\|_{\ell^\infty_\lambda(L^2(\omega))}
    \le c\bigl(\tau^{\frac14-\epsilon}+\lambda^{-\frac14+\epsilon}\bigr)
    \|(q_1^n-q_2^n)_{n=1}^{N}\|_{\ell^\infty_\lambda(L^2(\omega))},
\end{equation*}
where $c$ is independent of $h$, $\tau$, $N$, $\lambda$, $q_1$, and $q_2$, but may depend on $\epsilon$.
\end{proof}

The next result establishes the convergence of the scheme \eqref{eqn:recon_dis} and an error estimate for its limit.

\begin{thm}\label{thm:error}
Let Assumption \ref{assum:reg_err} hold. Fix $\epsilon\in(0,\frac18]$. There exist $\lambda_0\ge1$, $\widetilde{\tau}\in(0,1]$, and $c,C>0$, independent of $\delta$, $h$, $\tau$, $N$, and $\lambda$, such that, for all $\lambda\ge\lambda_0$ and $0<\tau\le\widetilde{\tau}$, the following statements hold. For any initial guess $q^0$ satisfying $|q^{0,n}|\le b$ a.e. in $\omega$, $n=0,\dots,N$, the scheme \eqref{eqn:recon_dis} converges to a unique fixed point $q^*$ satisfying the box constraint. Moreover,
\begin{align*}
\|(q^{k,n}-q^{*,n})_{n=1}^{N}\|_{\ell^\infty_\lambda(L^2(\omega))}\le
\left[c\bigl(\tau^{\frac14-\epsilon}+\lambda^{-\frac14+\epsilon}\bigr)\right]^k
\|(q^{0,n}-q^{*,n})_{n=1}^{N}\|_{\ell^\infty_\lambda(L^2(\omega))},
    \quad k\ge0,
\end{align*}
where $c(\tau^{\frac14-\epsilon}+\lambda^{-\frac14+\epsilon})\le\frac12$.
The fixed point satisfies
\begin{equation}\label{eqn:error}
\begin{aligned}
\|(q^{*,n}-q^\dagger(t_n))_{n=1}^{N}\|_{\ell^\infty(L^2(\omega))}\le C\left(
   h^{-1}\delta+\tau^{-1}\delta
+(1+|\log h|)h^{\frac12}
+(1+|\log\tau|)\tau^{\frac14}
    \right).
\end{aligned}
\end{equation}
If, in addition, the initial compatibility conditions \eqref{eqn:initial_compatibility_first_order} hold, then
\begin{equation}\label{eqn:error_first_order}
\begin{aligned}
\|(q^{*,n}-q^\dagger(t_n))_{n=1}^{N}\|_{\ell^\infty(L^2(\omega))}\le C\left(
h^{-1}\delta+\tau^{-1}\delta
 +(1+|\log h|)h^{\frac12}+\tau
\right).
    \end{aligned}
\end{equation}
The constants may depend on $\epsilon$ and the fixed problem data.
\end{thm}
\begin{proof}
For an admissible sequence $q=(q^n)_{n=1}^{N}$, we define
\begin{equation*}
\mathcal{G}_{h,\tau}(q)=\left(P_{\mathcal{Q}_b}\left[\left.\partial_{x_d}v_h^n(q)\right|_{x_d=\ell}\right]\right)_{n=1}^{N}.
\end{equation*}
The element-wise derivative $\left.\partial_{x_d}v_h^n(q)\right|_{x_d=\ell}$ belongs to $L^2(\omega)$, and the cut-off operator ensures that each component of $\mathcal{G}_{h,\tau}(q)$ satisfies the box constraint. Hence $\mathcal{G}_{h,\tau}$ is a self-map on the set of admissible sequences. This set is closed in $(L^2(\omega))^N$ and therefore complete in the weighted norm. For any two admissible sequences $q_i=(q_i^n)_{n=1}^{N}$, $i=1,2$, Lemmas \ref{l4} and \ref{lem:dis_forward_stab} give
\begin{align*}
\|\mathcal{G}_{h,\tau}(q_1)-\mathcal{G}_{h,\tau}(q_2)\|_{\ell^\infty_\lambda(L^2(\omega))}
&\le \left\|
(\partial_{x_d}
\bigl(v_h^n(q_1)-v_h^n(q_2)\bigr)
|_{x_d=\ell})_{n=1}^{N}
\right\|_{\ell^\infty_\lambda(L^2(\omega))}\\
    &\le
    c\bigl(\tau^{\frac14-\epsilon}+\lambda^{-\frac14+\epsilon}\bigr)
    \|(q_1^n-q_2^n)_{n=1}^{N}\|_{\ell^\infty_\lambda(L^2(\omega))}.
\end{align*}
The constant $c$ is uniform with respect to the noisy data satisfying Assumption \ref{assum:reg_err}, since the bounds on those data are independent of $\delta$. Now we choose $\lambda_0\ge1$ sufficiently large and $\widetilde{\tau}\in(0,1]$ sufficiently small, within the range of validity of the preceding lemmas, such that $\widetilde{\tau}b<1$ and $c\bigl(\widetilde{\tau}^{\frac14-\epsilon}+\lambda_0^{-\frac14+\epsilon}\bigr)\le\frac12$.
Then $\mathcal{G}_{h,\tau}$ is a contraction for every $\lambda\ge\lambda_0$ and $0<\tau\le\widetilde{\tau}$. The Banach fixed-point theorem therefore yields a unique fixed point $(q^{*,n})_{n=1}^{N}$ and
\begin{equation*}
\|(q^{k,n}-q^{*,n})_{n=1}^{N}\|_{\ell^\infty_\lambda(L^2(\omega))}\le\left[c\bigl(\tau^{\frac14-\epsilon}+\lambda^{-\frac14+\epsilon}\bigr)\right]^k\|(q^{0,n}-q^{*,n})_{n=1}^{N}\|_{\ell^\infty_\lambda(L^2(\omega))},\quad k\ge0.
\end{equation*}
At the initial time level, $v_h^0=R_h(\partial_{x_d}u_0)$ is independent of the coefficient. Consequently,
\begin{equation*}
    q^{k,0}=q^{*,0}
    =
    P_{\mathcal{Q}_b}\left[
        \left.\partial_{x_d}
        R_h(\partial_{x_d}u_0)\right|_{x_d=\ell}
    \right],
    \quad k\ge1.
\end{equation*}
This proves convergence and uniqueness for the full sequence $(q^{*,n})_{n=0}^{N}$. We next estimate the reconstruction error. Since $\mathcal{G}_{h,\tau}$ does not depend on $\lambda$, neither does its unique fixed point. We therefore fix $\lambda=\lambda_0$ for the remainder of the proof. The exact coefficient $q^\dag$ satisfies
\begin{equation*}
    q^\dagger(t_n)
    =\left.\partial_{x_d}v^\dagger(t_n)\right|_{x_d=\ell},
    \quad 1\le n\le N.
\end{equation*}
Since $q^\dagger(t_n)$ satisfies the box constraint, Lemma \ref{l4} and the fixed-point identity imply
\begin{align*}
    & \|q^{*,n}-q^\dagger(t_n)\|_{L^2(\omega)}
    \le
    \left\|
        \left.\partial_{x_d}
        \bigl(v_h^n(q^*)-v^\dagger(t_n)\bigr)
        \right|_{x_d=\ell}
    \right\|_{L^2(\omega)}\\
    \le&
    \left\|
        \left.\partial_{x_d}
        \bigl(v_h^n(q^*)-v_h^n(q^\dagger)\bigr)
        \right|_{x_d=\ell}
    \right\|_{L^2(\omega)}+
    \left\|
        \left.\partial_{x_d}
        \bigl(v_h^n(q^\dagger)-\widetilde v_h^n(q^\dagger)\bigr)
        \right|_{x_d=\ell}
    \right\|_{L^2(\omega)}\\
    &+
    \left\|
        \left.\partial_{x_d}
        \bigl(\widetilde v_h^n(q^\dagger)-v^\dagger(t_n)\bigr)
        \right|_{x_d=\ell}
    \right\|_{L^2(\omega)}.
\end{align*}
By multiplying by $e^{-\lambda_0t_n}$, taking the maximum over $1\le n\le N$, and applying Lemmas \ref{lem:dis_forward_stab}, \ref{lem:wnh-wnh0}, and \ref{lem:wnh0-wdag}, together with the choice of $\lambda_0$ and $\widetilde{\tau}$, we obtain
\begin{align*}
\|(q^{*,n}-q^\dagger(t_n))_{n=1}^{N}\|_{\ell^\infty_{\lambda_0}(L^2(\omega))}
    &\le
    \frac12
    \|(q^{*,n}-q^\dagger(t_n))_{n=1}^{N}\|_{\ell^\infty_{\lambda_0}(L^2(\omega))}\\
    &\quad+
    C\left(
        h^{-1}\delta+\tau^{-1}\delta
        +(1+|\log h|)h^{\frac12}
        +(1+|\log\tau|)\tau^{\frac14}
    \right).
\end{align*}
This directly yields
\begin{equation*}
\|(q^{*,n}-q^\dagger(t_n))_{n=1}^{N}\|_{\ell^\infty_{\lambda_0}(L^2(\omega))}\le C\left(h^{-1}\delta+\tau^{-1}\delta+(1+|\log h|)h^{\frac12}+(1+|\log\tau|)\tau^{\frac14}\right).
\end{equation*}
Finally, by the equivalency between norms, we conclude that 
\begin{equation*}
    \|(q^{*,n}-q^\dagger(t_n))_{n=1}^{N}\|_{\ell^\infty(L^2(\omega))}
    \le e^{\lambda_0T}
    \|(q^{*,n}-q^\dagger(t_n))_{n=1}^{N}\|_{\ell^\infty_{\lambda_0}(L^2(\omega))}.
\end{equation*}
Since $\lambda_0$ depends only on $\epsilon$ and the fixed problem data, the factor $e^{\lambda_0T}$ can be absorbed into $C$ without introducing any dependence on $\delta$, $h$, $\tau$, $N$, or $\lambda$. This proves \eqref{eqn:error}.
Under the additional conditions \eqref{eqn:initial_compatibility_first_order}, Corollary \ref{cor:wnh0-wdag-first-order} replaces the temporal discretization term $(1+|\log\tau|)\tau^{1/4}$ by $\tau$. The same argument then proves \eqref{eqn:error_first_order}.
\end{proof}

\begin{rmk}\label{rmk:err_highorder}
The spatial error estimate in Theorem \ref{thm:error} can be improved under additional regularity assumptions. Specifically, if $G$ in \eqref{eqn:bdry_lift} satisfies $G\in C^1([0,T];W^{2,\infty}(\Omega))$, the $L^\infty(\Omega)$ error estimate in \cite[Theorem 6.8]{Thomee:2006} yields an improved spatial discretization term $(1+|\log h|)h$. If, in addition, the initial compatibility conditions \eqref{eqn:initial_compatibility_first_order} hold, then
\begin{equation*}
\|(q^{*,n}-q^\dagger(t_n))_{n=1}^{N}\|_{\ell^\infty(L^2(\omega))}\le C\left(h^{-1}\delta+\tau^{-1}\delta+(1+|\log h|)h+\tau\right).
\end{equation*}
For sufficiently small $\delta>0$, choosing $\tau\asymp\delta^{\frac12}$ and $h\asymp\delta^{\frac12}[\log(\delta^{-1})]^{-\frac12}$ gives
\begin{equation*}
\|(q^{*,n}-q^\dagger(t_n))_{n=1}^{N}\|_{\ell^\infty(L^2(\omega))}\le C\delta^{\frac12}[\log(\delta^{-1})]^{\frac12}.
\end{equation*}
\end{rmk}

Combining the convergence and reconstruction error estimates in Theorem \ref{thm:error} gives the following estimate for the iterates.

\begin{cor}
Under the assumptions of Theorem \ref{thm:error}, for $0<\tau\le\widetilde{\tau}$ and $k\ge0$, the iterate $q^k$ satisfies
\begin{align*}
\|(q^{k,n}-q^\dagger(t_n))_{n=1}^{N}\|_{\ell^\infty(L^2(\omega))}
&\le e^{\lambda_0T}\left[c\bigl(\tau^{\frac14-\epsilon}+\lambda_0^{-\frac14+\epsilon}\bigr)\right]^k\|(q^{0,n}-q^{*,n})_{n=1}^{N}\|_{\ell^\infty_{\lambda_0}(L^2(\omega))}\\
&\quad+C\left(h^{-1}\delta+\tau^{-1}\delta+(1+|\log h|)h^{\frac12}+(1+|\log\tau|)\tau^{\frac14}\right).
\end{align*}
If, in addition, the initial compatibility conditions \eqref{eqn:initial_compatibility_first_order} hold, then the temporal discretization term $(1+|\log\tau|)\tau^{\frac14}$ can be replaced by $\tau$. The constants are as in Theorem \ref{thm:error}.
\end{cor}

\section{Numerical results and discussions}\label{sec:num}
Now we present numerical experiments to illustrate the reconstruction algorithm. Fix a square domain $\Omega = (0, 1)^2$ and $T=1$. The initial and boundary data are fixed at $u_0\equiv 1$ and $g\equiv 1$. The noisy data $z^\delta$ is measured on the boundary segment $(0, 1) \times \{1\}$ over the time interval $(0, T)$:
\begin{equation*}
    z^\delta(x_1,t)=\partial_{x_2}u^\dagger(x_1,1,t)+\delta\| \partial_{x_2}u^\dagger(\cdot,1,t)\|_{L^\infty(\omega)}\xi(x_1,t),
\end{equation*}
where $\delta$ denotes the relative noise level and $\xi$ follows the standard Gaussian distribution. The exact state $u^\dagger=u(q^\dagger)$  is computed on  a highly refined mesh  with $h=1/200$ and $\tau =1/1000$. The accuracy of a reconstruction $q^*$
relative to the exact one $q^\dagger$ is measured by the relative error  $e_q={\|q^*-q^\dagger\|_{\ell^\infty(L^2(\omega))}}/{\|q^\dagger\|_{L^\infty(0,T;L^2(\omega))}} $. 

The first example is about a smooth coefficient.
\begin{exa}\label{ex:smooth}
       $q^\dagger(x_1,t)=\sin(2\pi x_1) t\sin(2\pi t)$.
\end{exa}

\begin{figure}[htbp]
	\centering
	\begin{tabular}{ccc}
		\includegraphics[width=0.30\textwidth]{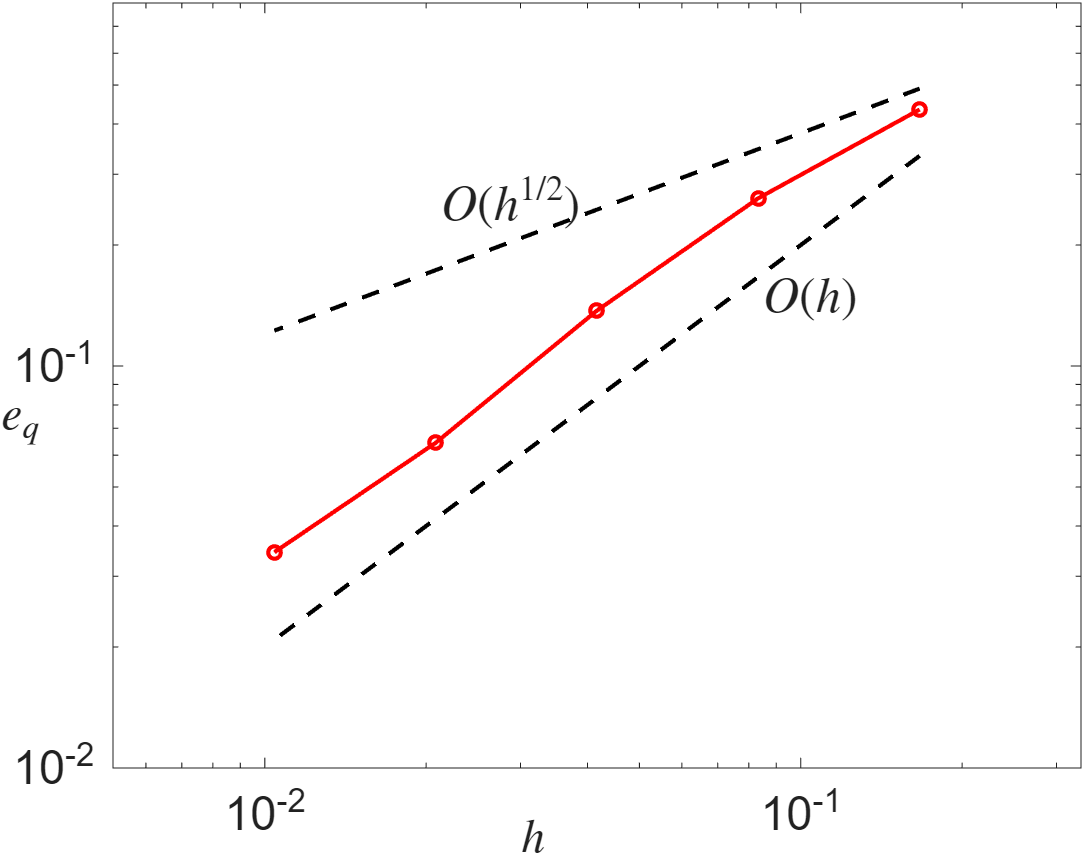}&
		\includegraphics[width=0.30\textwidth]{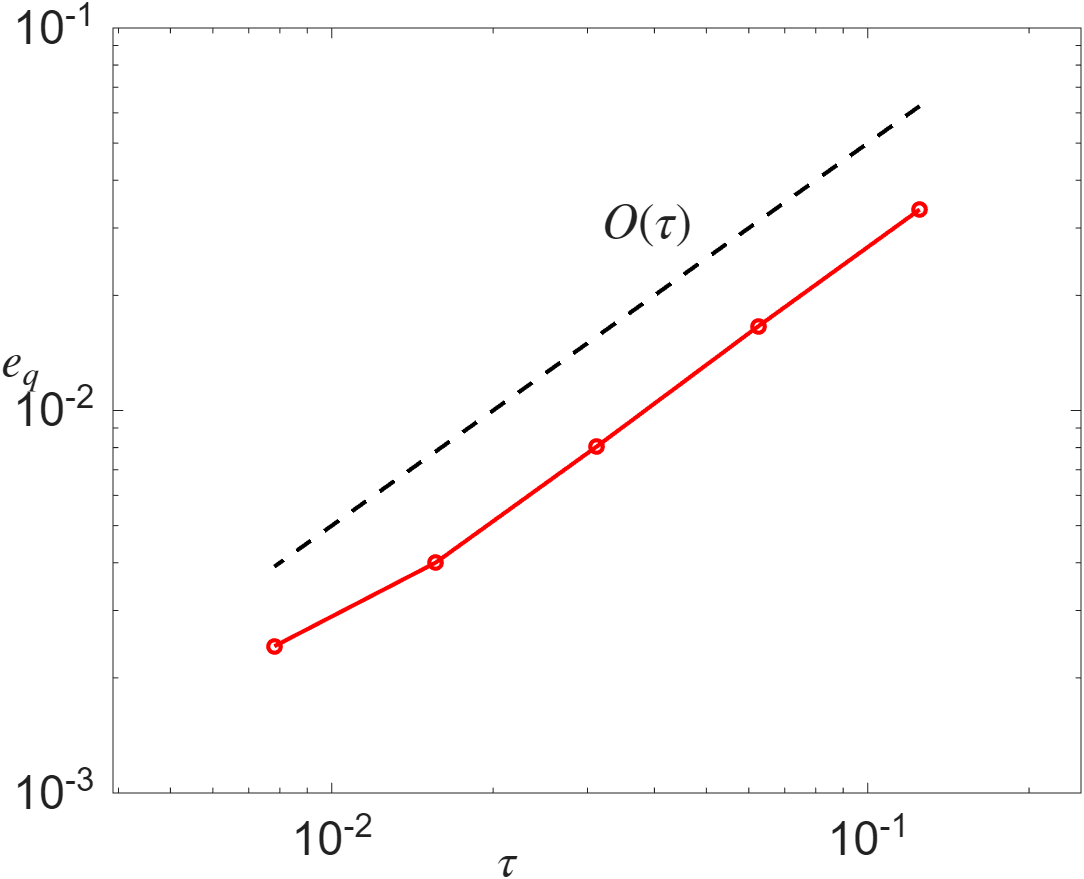}&
		\includegraphics[width=0.30\textwidth]{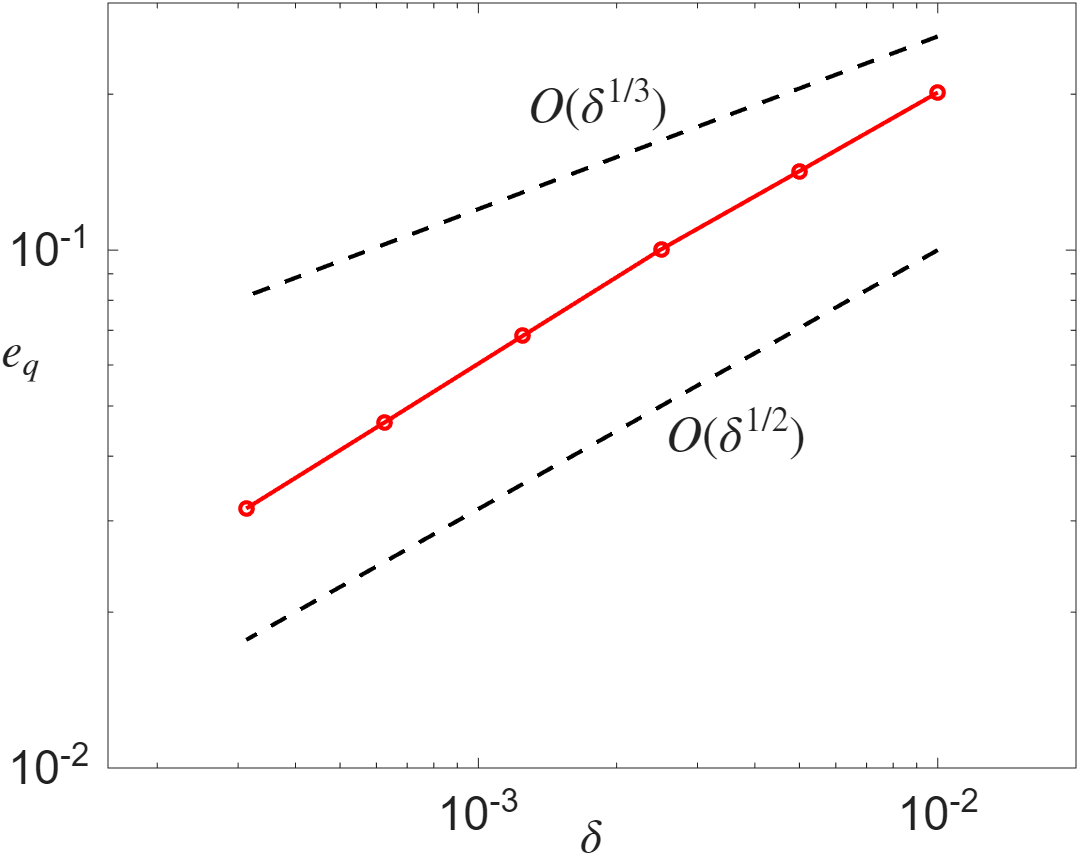}\\
		(a) $e_q$ versus $h$ & (b) $e_q$ versus $\tau$  & (c) $e_q$ versus $\delta$  
	\end{tabular}
	\caption{The convergence of $e_q$ with respect to $h$, $\tau$ and $\delta$ for Example \ref{ex:smooth}.}
	\label{Fig:conv_smooth}
\end{figure}

\begin{figure}[htbp]
	\centering
	\begin{tabular}{ccc}
		\includegraphics[width=0.25\textwidth]{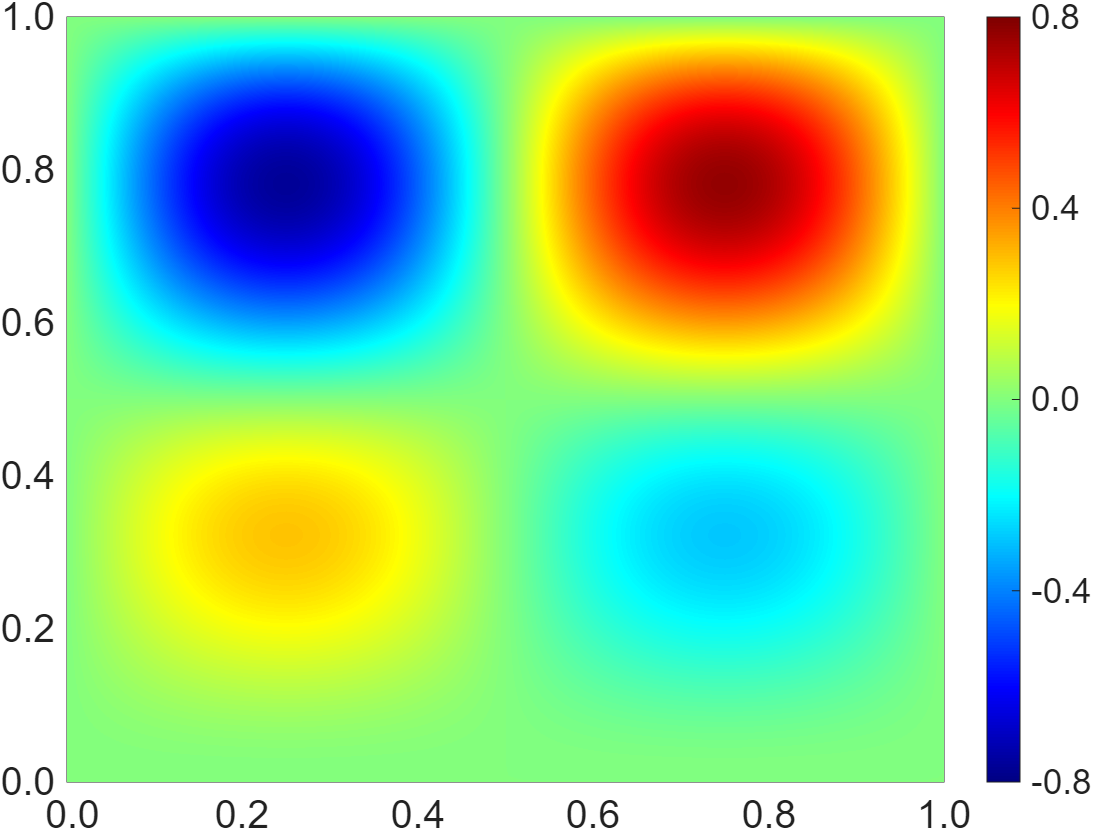}&
		\includegraphics[width=0.25\textwidth]{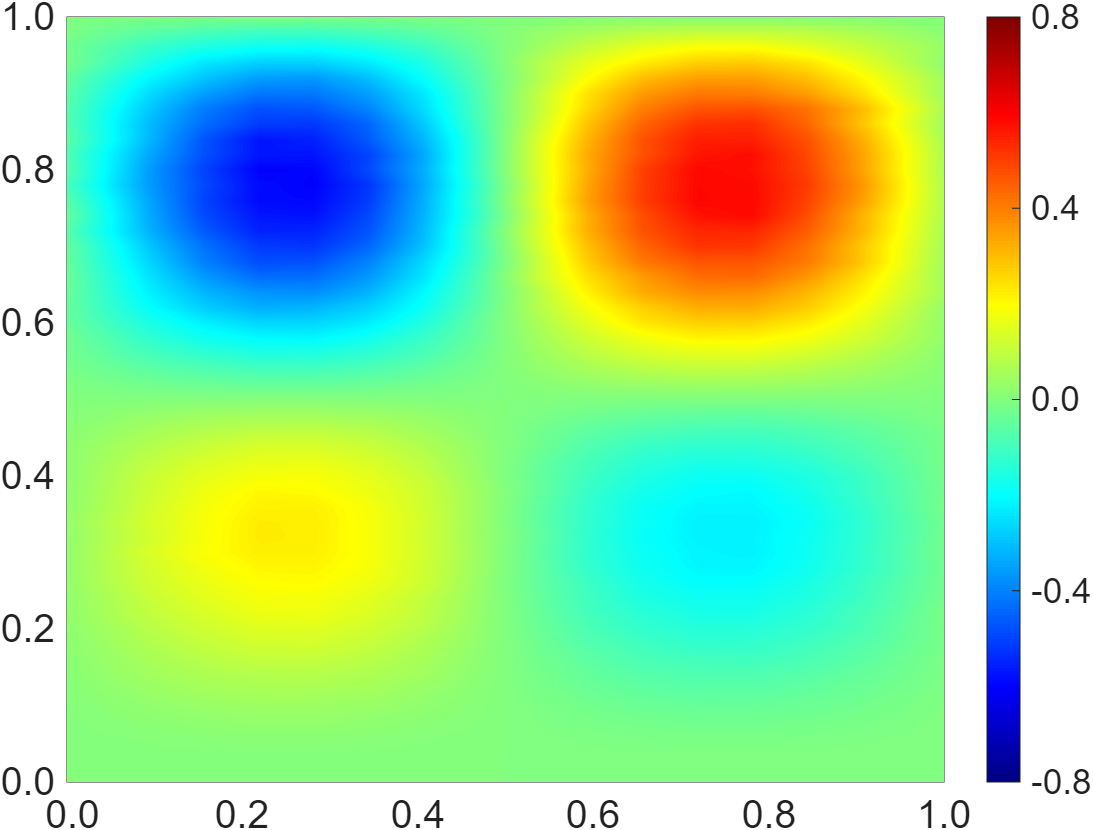}&
		\includegraphics[width=0.25\textwidth]{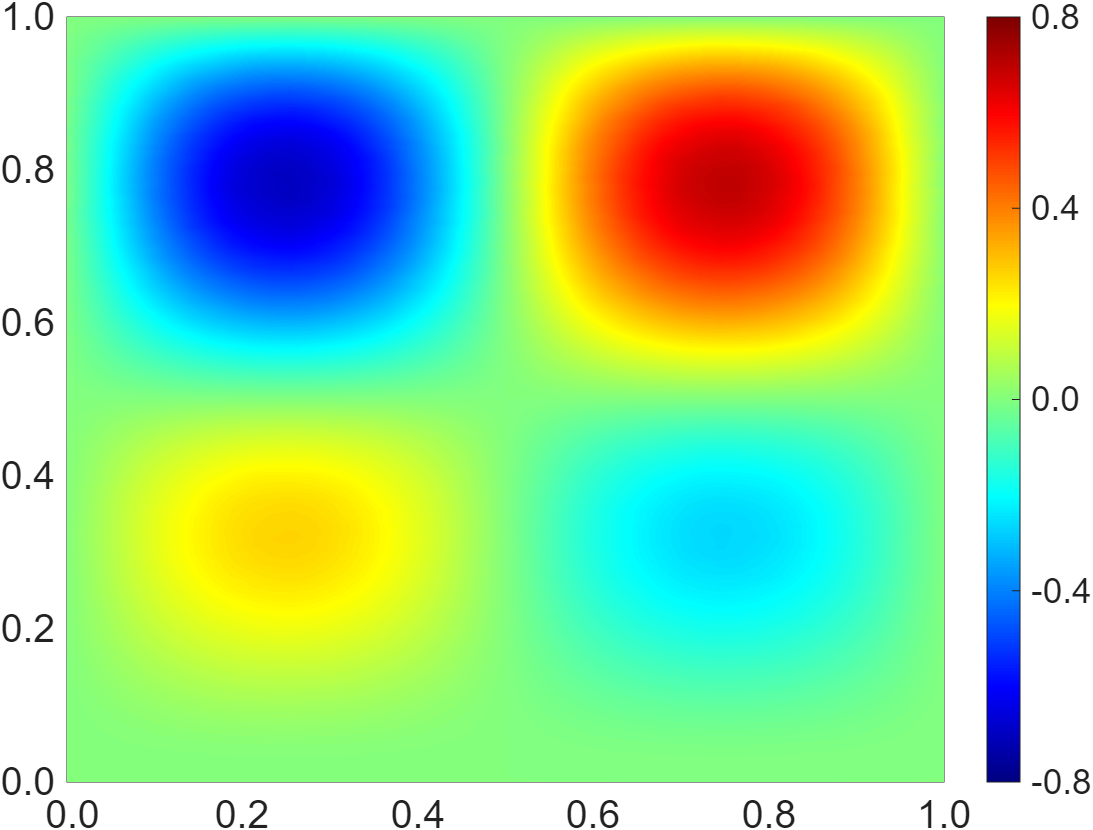} \\
		(a) exact & (b) $\delta=\text{1e-2}$  & (c) $\delta=\text{1.25e-3}$  
	\end{tabular}
	\caption{The reconstruction $q^*$  for Example \ref{ex:smooth}   at different noise levels.}
	\label{Fig:recon_smooth}
\end{figure}

First we verify the following error bound in Theorem \ref{thm:error} and Remark \ref{rmk:err_highorder}:
\begin{equation*} 
        \| ((q^*)^n-(q^\dagger)^n)_{n=0}^{N} \|_{ \ell^\infty(L^2(\omega))}\le   c(h^{-1}\delta+\tau^{-1}\delta+|\log h|  h    + \tau  ). 
\end{equation*}
 For $\delta = 0$, Theorem \ref{thm:error} and Remark \ref{rmk:err_highorder} predict a spatial convergence rate   $O(|\log h|h)$ and a temporal rate   $O(\tau)$.  Fig.~\ref{Fig:conv_smooth} (a)-(b) indicates that the empirical convergence rates are $O(h^{0.93})$ in space and $O(\tau^{0.96})$ in time, which align well with these theoretical predictions. To compute the convergence rate with respect to the noise level $\delta$, we initially fix $h=\frac{1}{16}$ and $\tau=\frac{1}{50}$ for $\delta=10^{-2}$ and   refine the discretization parameters according to the scaling relations $h^{-1}\delta\sim h$ and $\tau^{-1}\delta\sim \tau$. The   empirical rate $O(\delta^{0.53})$  in Fig.~\ref{Fig:conv_smooth} (c)  matches the theoretical prediction   $O(\delta^{\frac{1}{2}}|\log\delta|^{\frac{1}{2}})$ closely. The recovered coefficients $q^*$ at different noise levels are presented in Fig.~\ref{Fig:recon_smooth}. The numerical results show that the reconstructions maintain high quality and accuracy.


\begin{figure}[htbp]
	\centering
	\begin{tabular}{ccc}
		\includegraphics[width=0.30\textwidth]{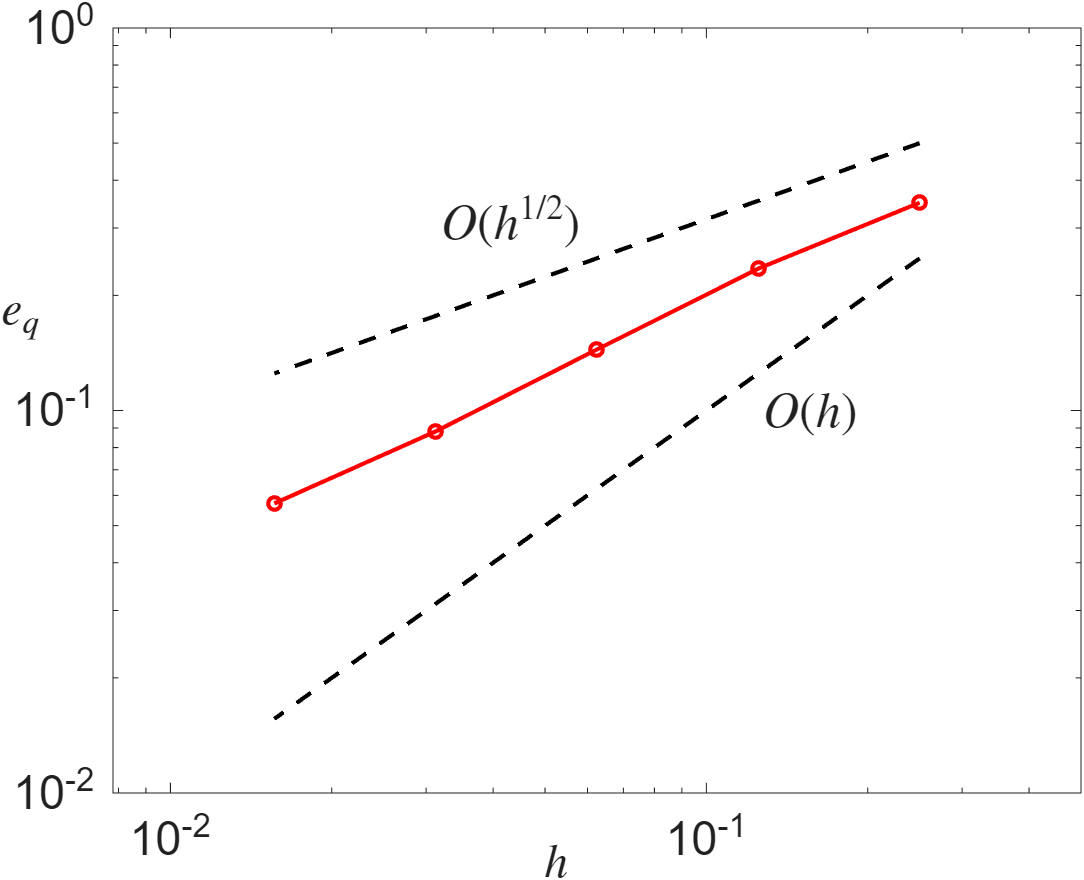}&
		\includegraphics[width=0.30\textwidth]{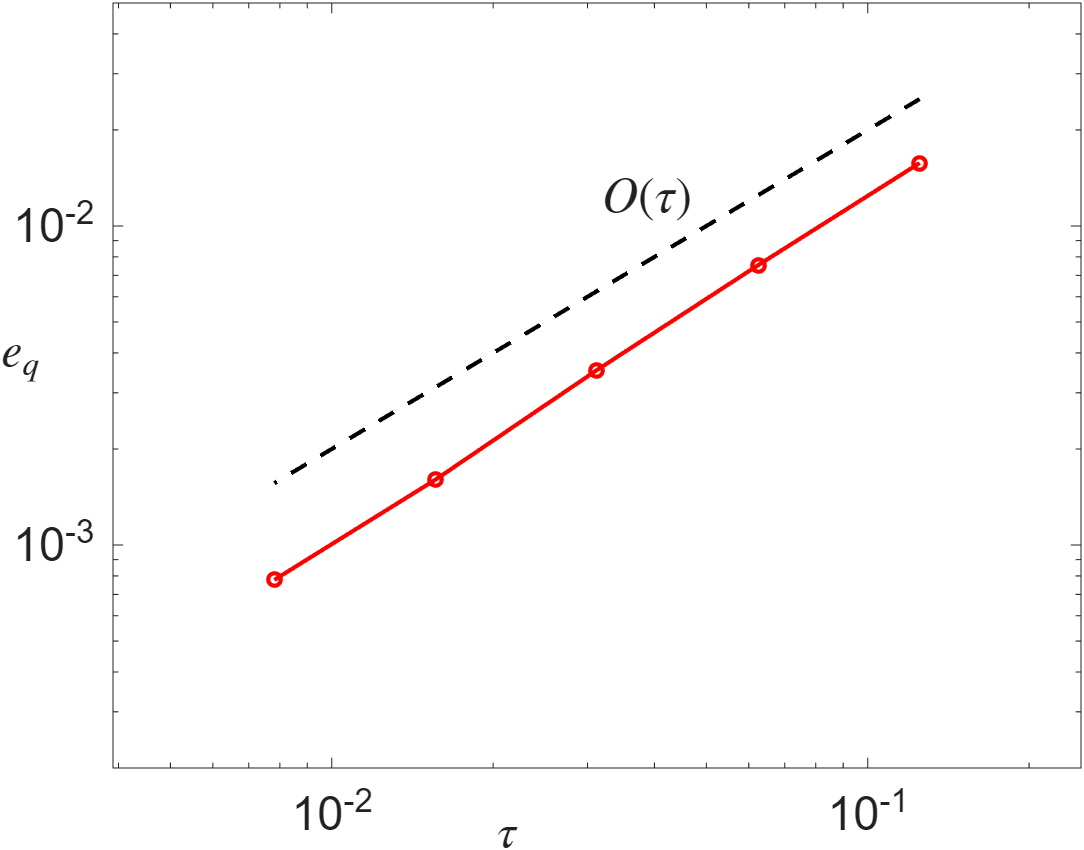}&
		\includegraphics[width=0.30\textwidth]{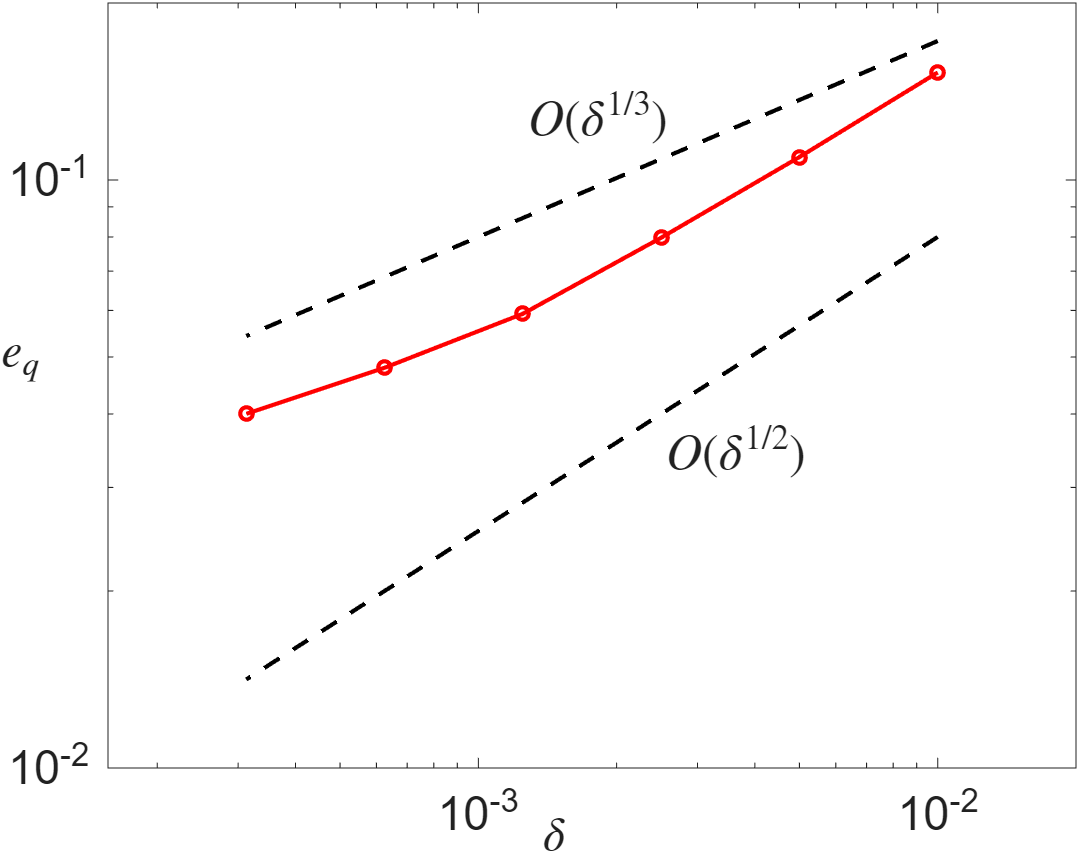}\\
		(a) $e_q$ versus $h$ & (b) $e_q$ versus $\tau$  & (c) $e_q$ versus $\delta$  
	\end{tabular}
	\caption{The convergence of $e_q$ with respect to $h$, $\tau$ and $\delta$ for Example \ref{ex:nonsmooth}.}
	\label{Fig:conv_nonsmooth}
\end{figure}

Next, we investigate  nonsmooth coefficient reconstruction to verify the sharpness of the error estimate in Theorem \ref{thm:error}.
\begin{exa}\label{ex:nonsmooth}
       $q^\dagger(x_1,t)=x_1^{0.4}(1-x_1)^{0.4} t(e^{t}-1)$.
\end{exa}
The coefficient $q^\dag$ exhibits weak singularities at $x_1=0$ and $x_1=1$, and thus belongs to $H^{0.9-\epsilon}(\omega)$ in space for an arbitrarily small $\epsilon>0$. We investigate the sharpness of the error bound in Theorem \ref{thm:error}:
\begin{equation*} 
        \| ((q^*)^n-(q^\dagger)^n)_{n=0}^{N} \|_{ \ell^\infty(L^2(\omega))}\le   c(h^{-1}\delta+\tau^{-1}\delta+|\log h|  h^{\frac{1}{2}}    + \tau  ). 
\end{equation*}
When $\delta = 0$, Fig.~\ref{Fig:conv_nonsmooth} (a)-(b) show an empirical convergence rate $O(h^{0.66})$ in space, and $O(\tau^{1.09})$ in time. The spatial rate is slightly higher than the theoretical prediction $O(|\log h|h^{\frac{1}{2}})$. Note that the theoretical analysis relies on the regularity conditions of the measurement data rather than the coefficient itself, and the implicit mapping between the coefficient $q^\dagger$ and the data $\partial_{x_d}u(q^\dagger)(t,x',\ell)$ makes directly verifying the sharpness of Theorem \ref{thm:error} challenging. Nevertheless, as the regularity of $q^\dagger$ decreases, the empirical spatial rate approaches the theoretical prediction. The experiment is consistent with  slower spatial convergence for less regular coefficients, as indicated by Theorem \ref{thm:error}. 
To evaluate the convergence behavior with respect  to the noise level $\delta$, we take the parameters $h = \frac{1}{16}$ and $\tau= \frac{1}{50}$ for $\delta=10^{-2}$ and then choose $h $ and $ \tau$ based on the balancing relations in Theorem \ref{thm:error}, namely $h^{-1}\delta\sim h^{\frac{1}{2}}$ and $\tau^{-1}\delta\sim \tau$. 
The empirical convergence rate is about $O(\delta^{0.39})$, which is close to the theoretical prediction   $O(\delta^{\frac{1}{3}}|\log\delta|^{\frac{2}{3}})$.

\begin{figure}[htbp]
	\centering
	\begin{tabular}{ccc}
		\includegraphics[width=0.25\textwidth]{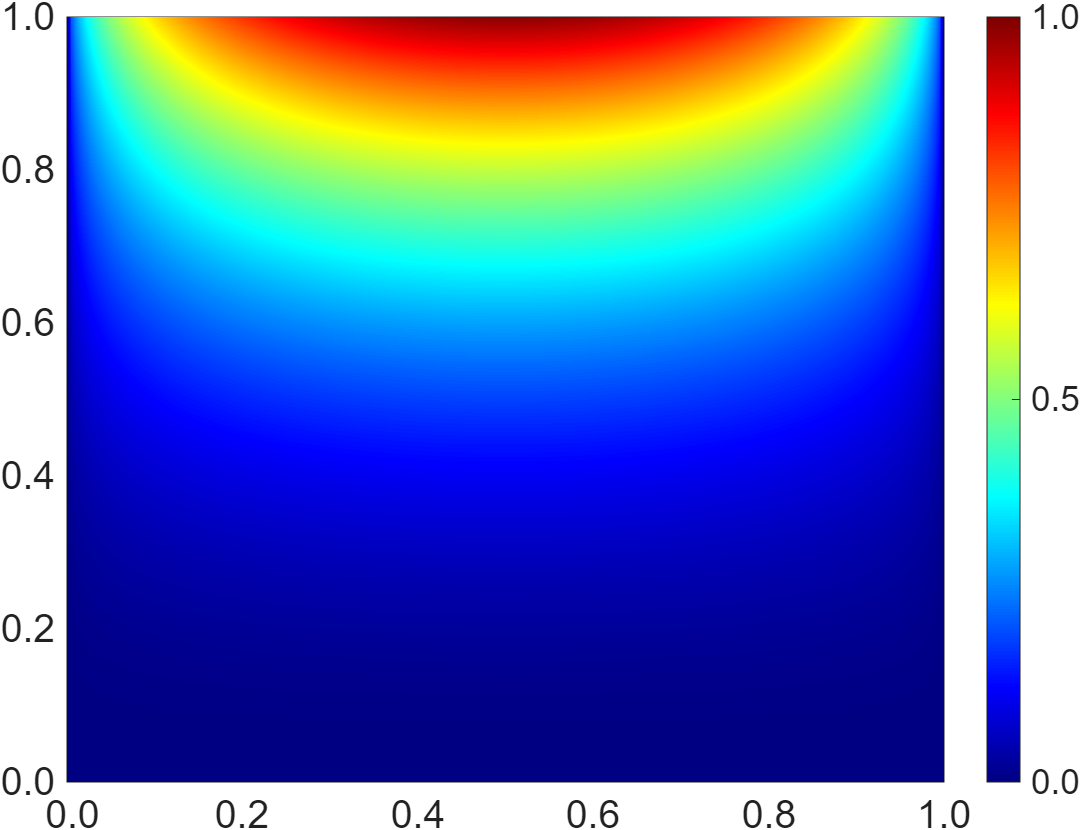}&
		\includegraphics[width=0.25\textwidth]{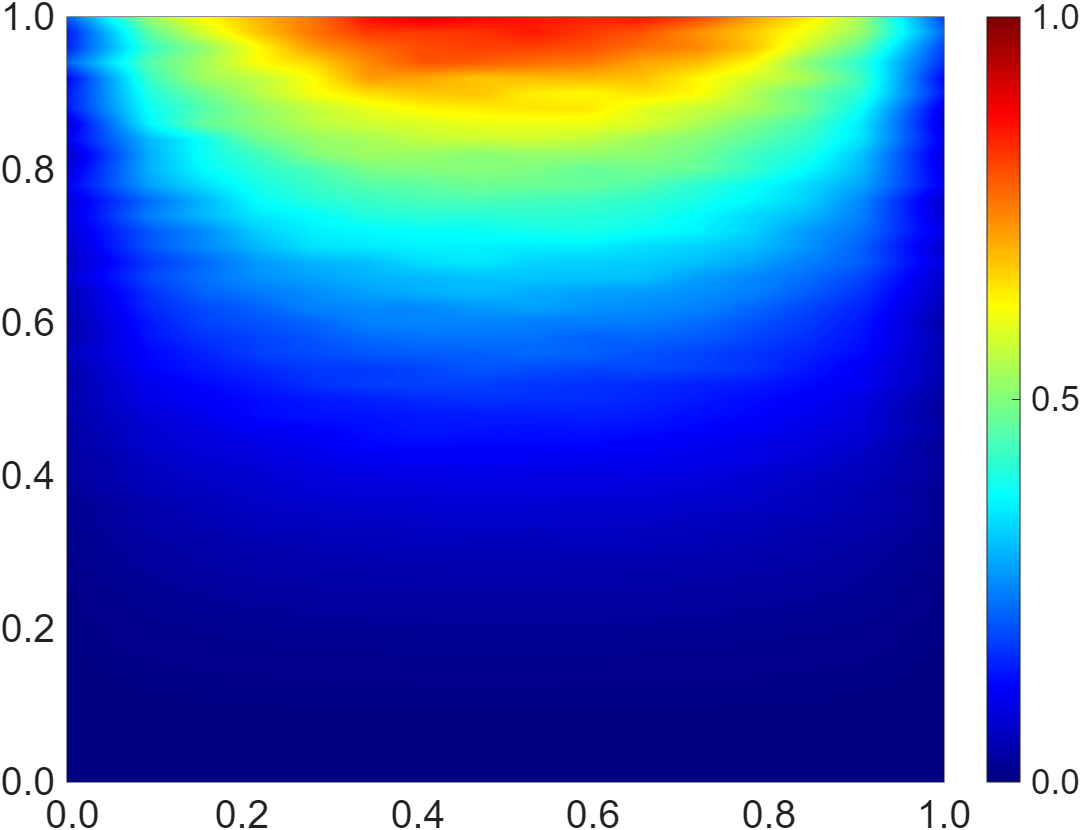}&
		\includegraphics[width=0.25\textwidth]{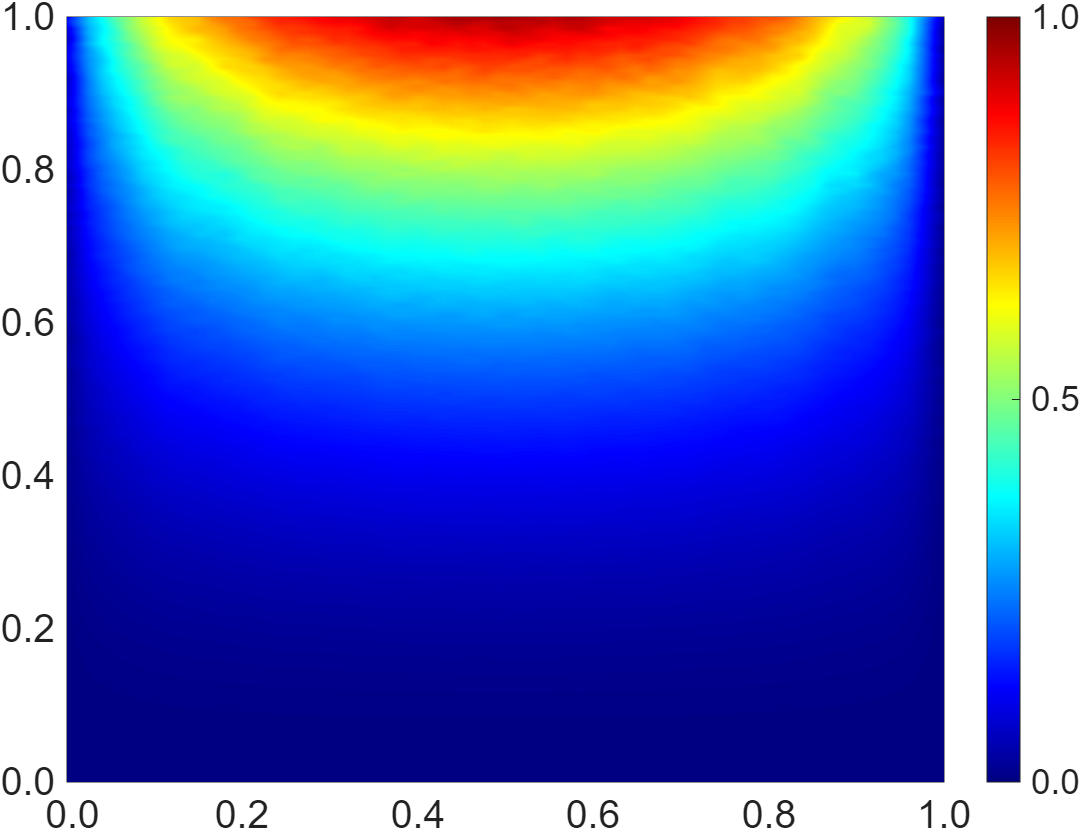} \\
		(a) exact & (b) $\delta=\text{1e-2}$  & (c) $\delta=\text{1.25e-3}$  
	\end{tabular}
	\caption{The reconstruction $q^*$  for Example \ref{ex:nonsmooth}   at different noise levels.}
	\label{Fig:recon_nonsmooth}
\end{figure}

The next example investigates the convergence behavior of the iteration scheme \eqref{eqn:recon_dis}.
\begin{exa}\label{ex:compare}
       $q^\dagger(x_1,t)=7( \cos(2\pi x_1)-1) t (e^t-1)$.
\end{exa}

\begin{figure}[htbp]
	\centering
	\begin{tabular}{cccc}
    		
		\includegraphics[width=0.20\textwidth]{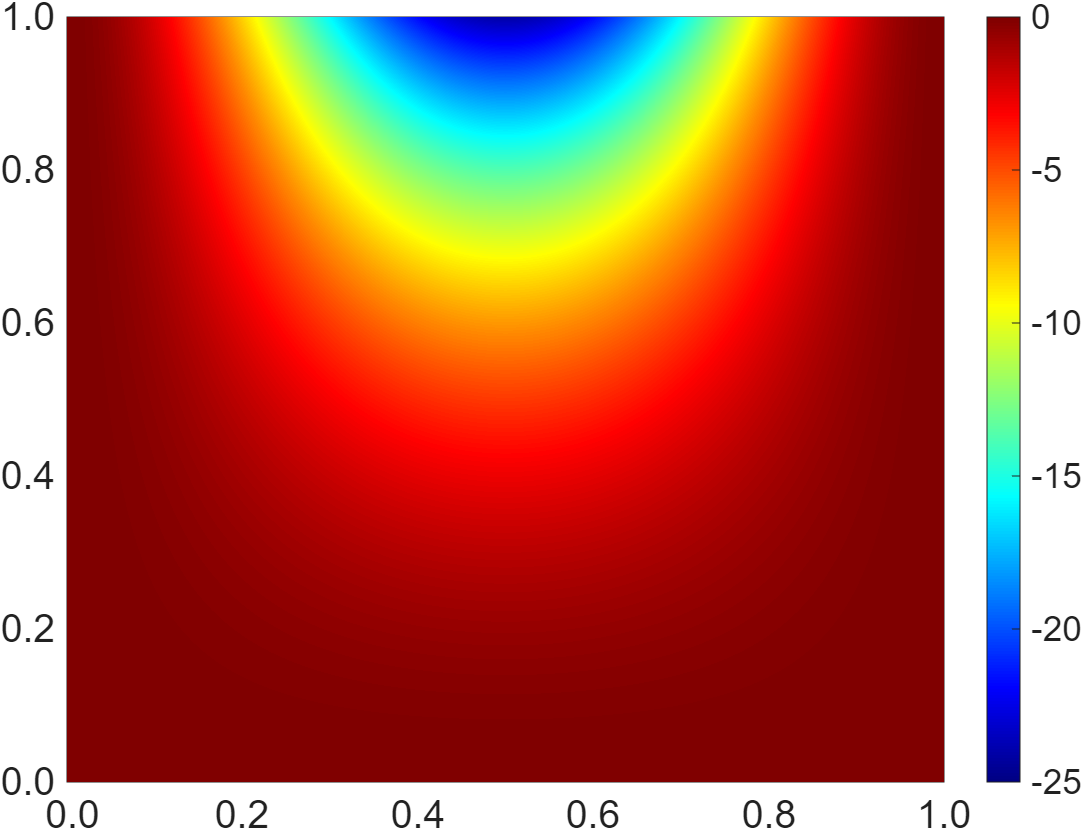}&
		\includegraphics[width=0.20\textwidth]{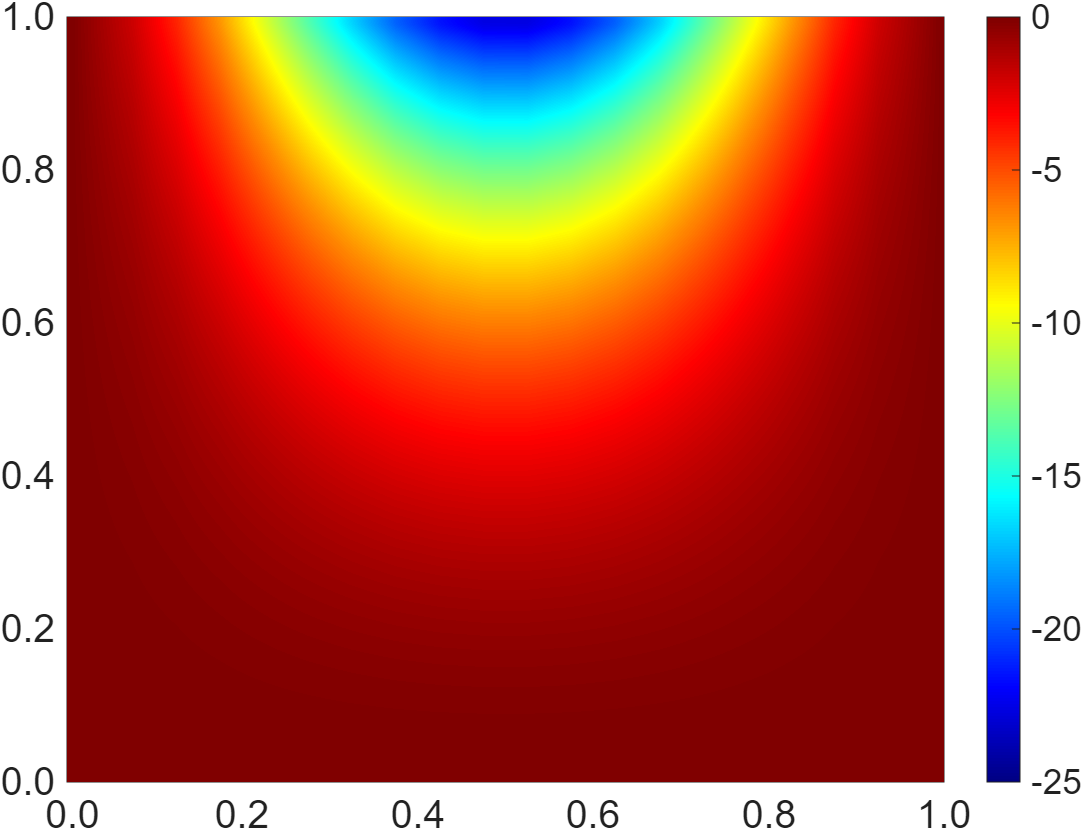}&
        \includegraphics[width=0.20\textwidth]{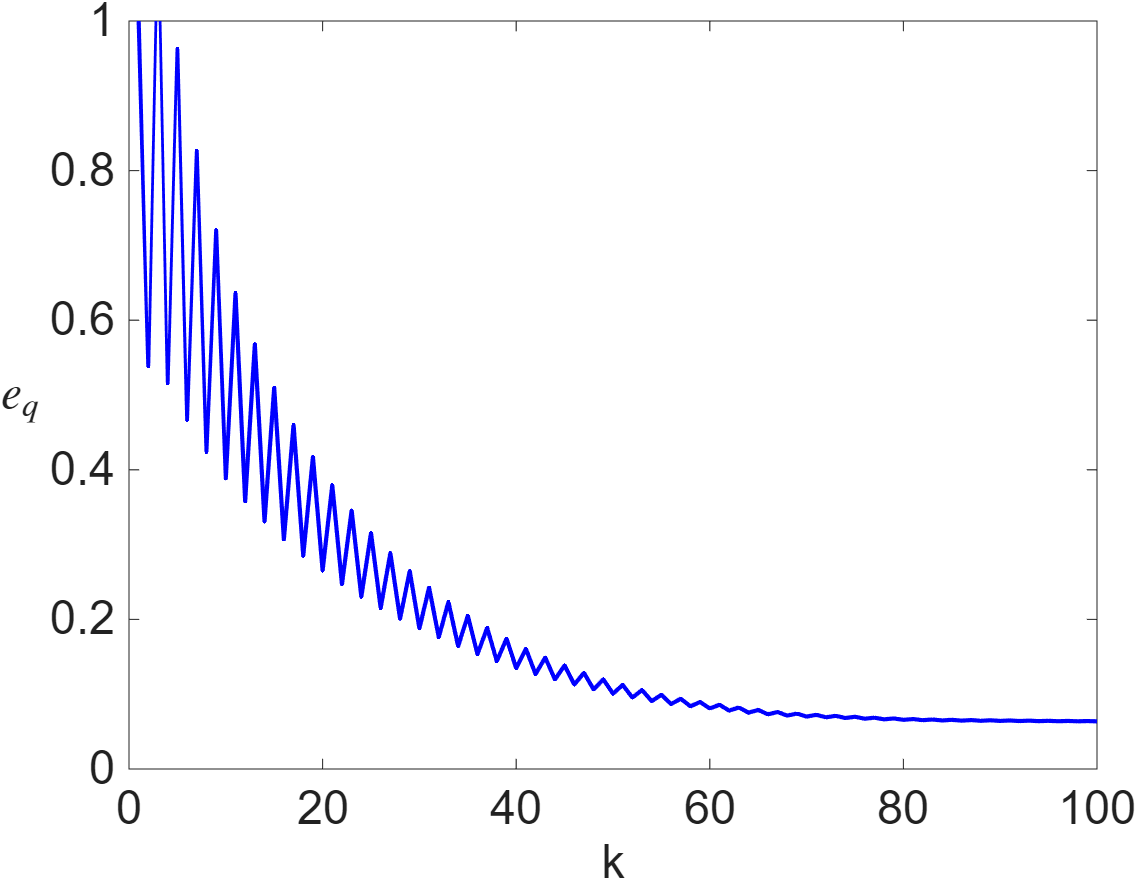}&
		\includegraphics[width=0.20\textwidth]{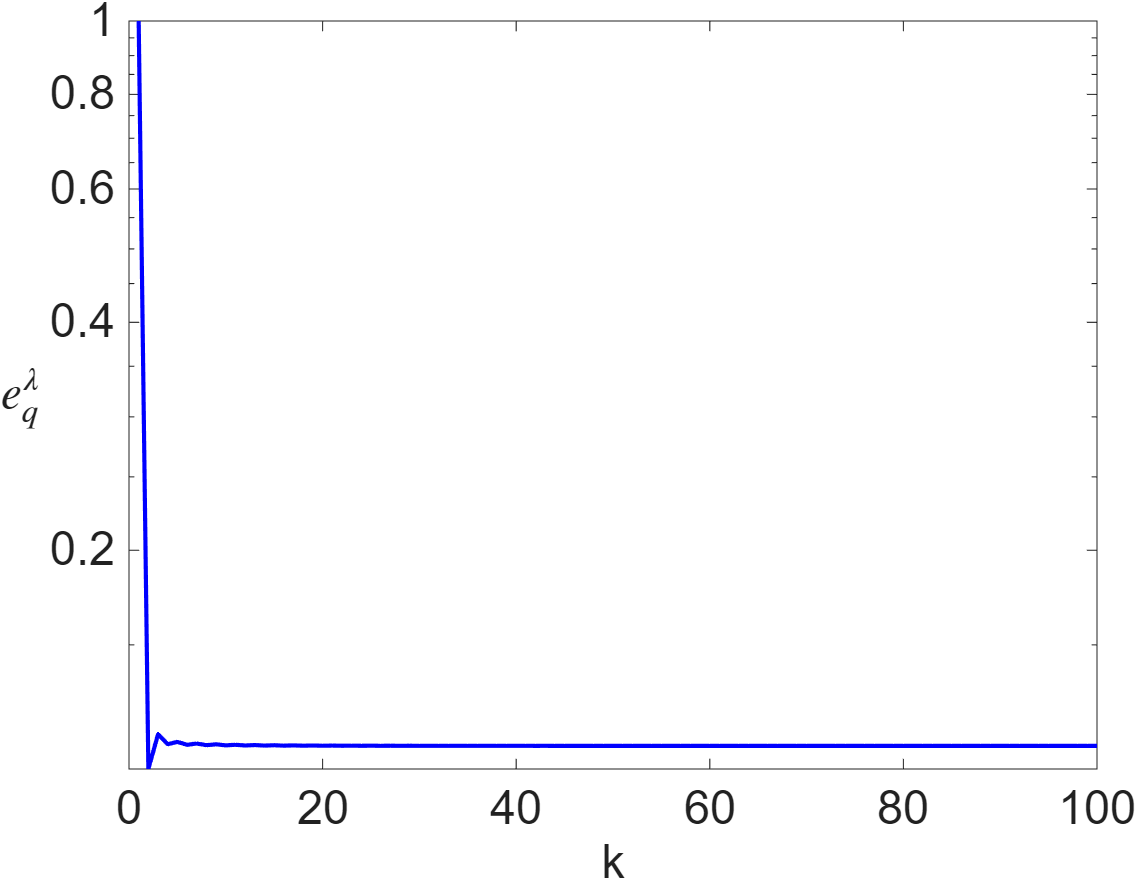}   \\
		  (a) exact $q^\dagger$ & (b) $q^*$ with $\delta=\text{1e-2}$  & (c) $e_q$ with $\lambda=0$ & (d) $e_q^\lambda$ with $\lambda=10$  
	\end{tabular}
	\caption{The reconstruction $q^*$ and the evolution of errors with respect to iteration number $k$.}
	\label{Fig:compare}
\end{figure}

The example is implemented with a mesh size $h=\frac{1}{20}$, time step size $\tau=\frac{1}{160}$, noise level $\delta=10^{-2}$ and $\mathcal{Q}_b=\{q\in L^\infty(0,T;L^\infty(\omega)):\|q\|_{L^\infty(0,T;L^\infty(\omega))}\le 25\}$. 
Fig.~\ref{Fig:compare} (a)-(b) present the reconstructed profile $q^*$. In Fig.~\ref{Fig:compare} (c)-(d), we plot the evolution of both the standard error $e_q$ and the relative weighted error $e_q^\lambda:={\|q^*-q^\dagger\|_{\ell^\infty_\lambda(L^2(\omega))}}/{\|q^\dagger\|_{L^\infty_\lambda(0,T;L^2(\omega))}}$ during the iteration. The scheme is initialized with $q^0 = 0$. While the trajectory of the unweighted error $e_q$ exhibits initial oscillations before stabilizing within 100 iterations, the convergence in the weighted norm is remarkably stable and rapid. This behavior matches our theoretical linear convergence result in Theorem \ref{thm:error}:  the weighted norm is essential for guaranteeing the contraction property of the fixed point iteration. Nonetheless, it is worth pointing out that the weight is  merely an analytic device for establishing the contraction of the iteration.

\begin{figure}[htbp]
	\centering
	\begin{tabular}{cccc}
		\includegraphics[width=0.25\textwidth]{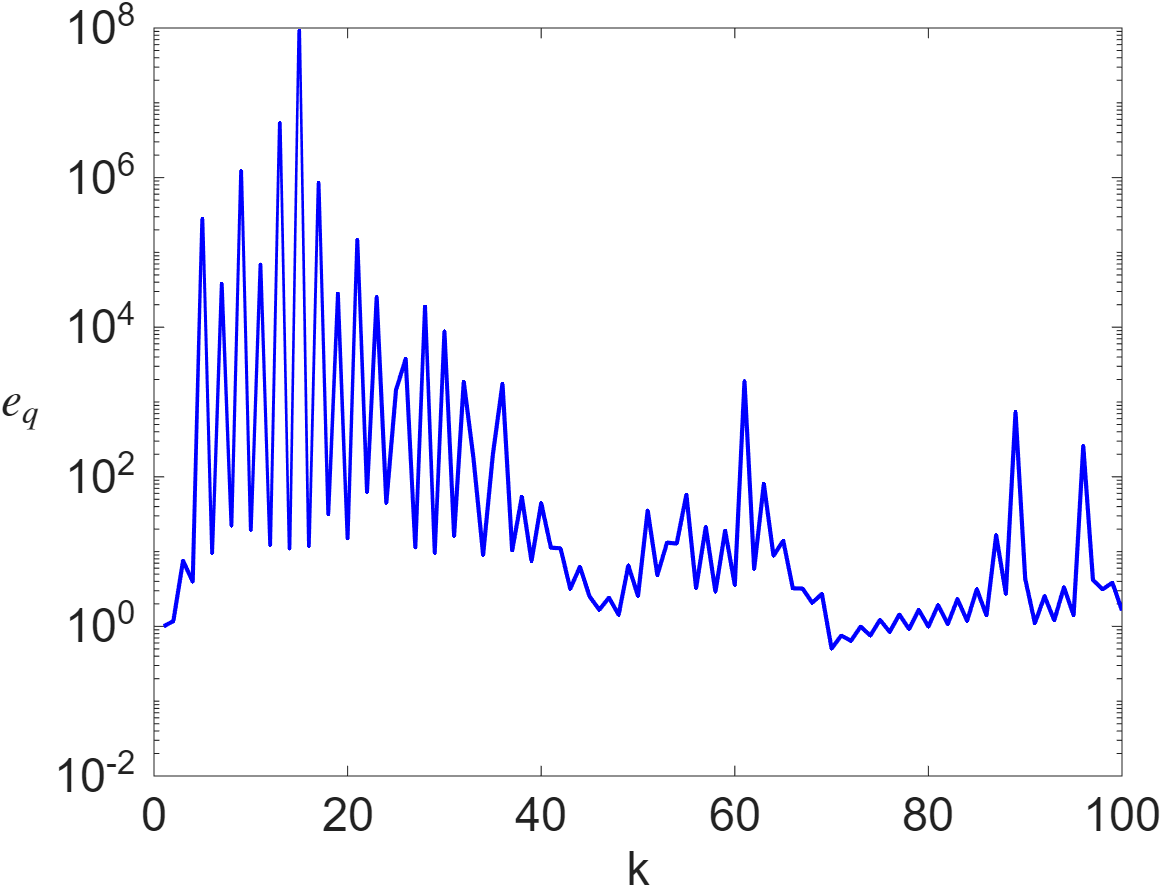}&
		\includegraphics[width=0.25\textwidth]{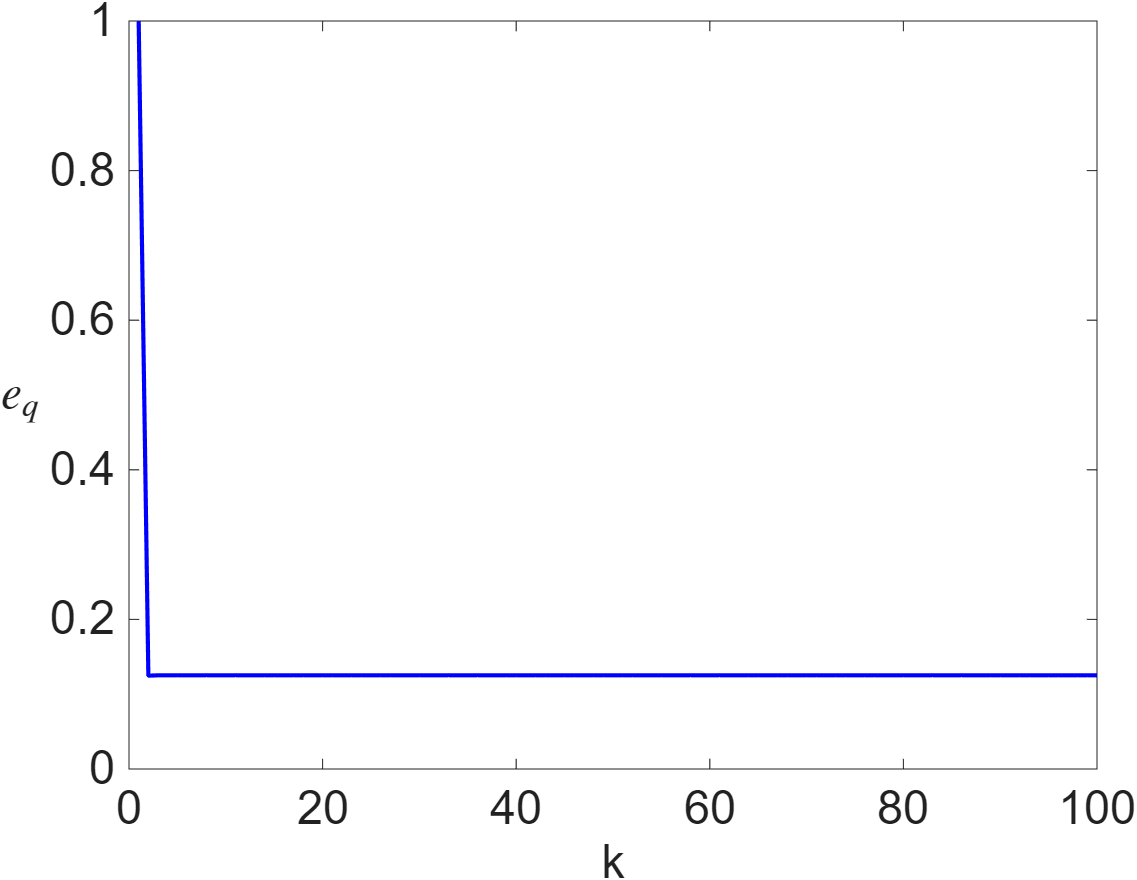}  &
		\includegraphics[width=0.25\textwidth]{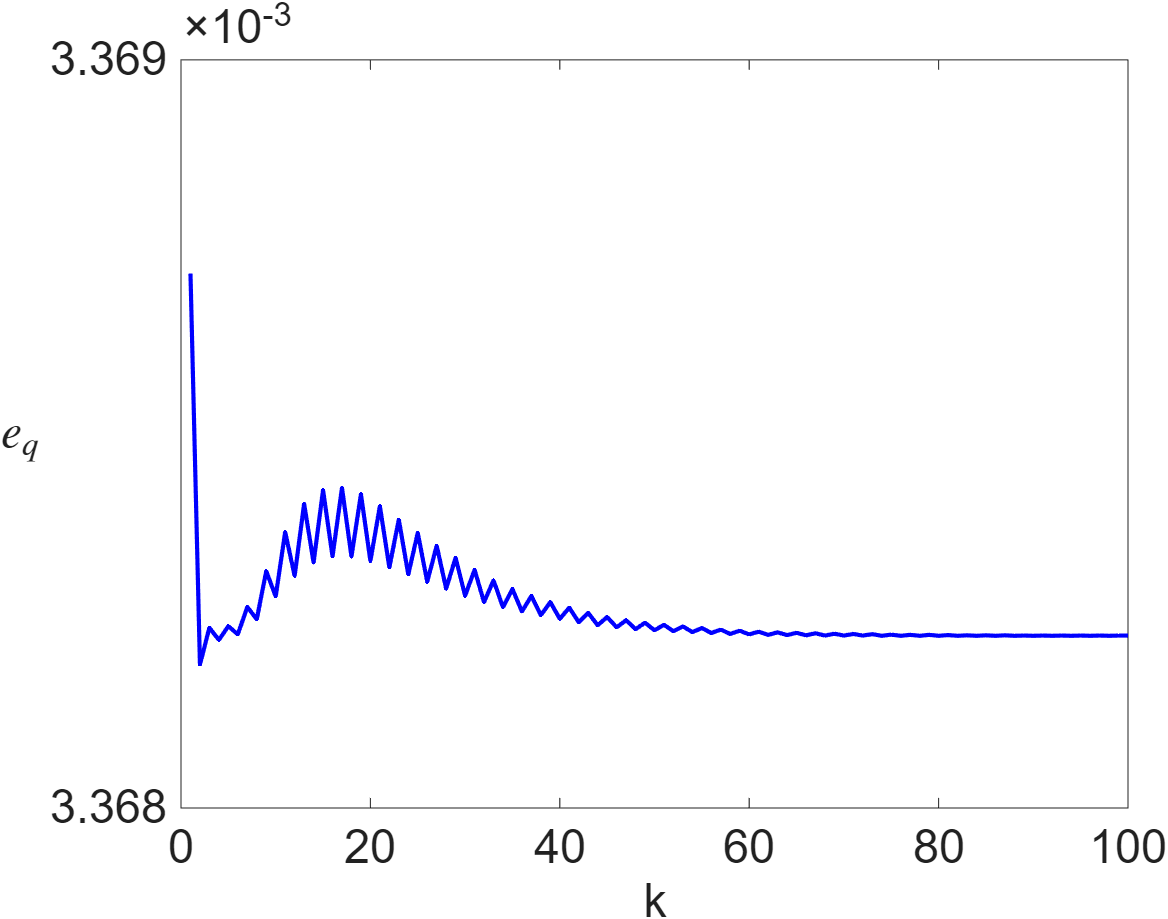} \\
		(a) without $P_{\mathcal{Q}_b}$ & (b) small $T$  & (c)  small $P$
	\end{tabular}
	\caption{The fixed point iteration without projection, and with small $T$ or small $P$.}
	\label{Fig:compare_loss}
\end{figure}

With $\delta=0$, we test the necessity of the \textit{a priori} condition $q^\dagger\in  {\mathcal{Q}_b}$.  Fig.~\ref{Fig:compare_loss} (a) shows the error evolution of scheme \eqref{eqn:recon_dis} without the cut-off operator $P_{\mathcal{Q}_b}$. The iteration diverges rapidly, which indicates the failure of the unprojected fixed point iteration in the specific scenario. Next we   investigate a local reconstruction scenario with small $T$ or $P$, following the framework of Theorem \ref{t2}. Specifically, we take $T=0.1$ and choose $q$ such that $\|q-q^\dagger\|_{L^\infty(0,T;L^2(\omega))}\le 10^{-2}\| q^\dagger\|_{L^\infty(0,T;L^2(\omega))}$ so as to generate $P=\|\mathcal{M}_{q^\dagger}(m)\|_{C([0,T];H^{s_1}(\omega))}$. Fig.~\ref{Fig:compare_loss} (b)-(c) show that the fixed-point iteration remains locally well-posed even without the \textit{a priori} bound $q^\dagger \in \mathcal{Q}_b$, which aligns well with the theoretical analysis in Theorem \ref{t2}


\bibliographystyle{abbrv}
\bibliography{ref}

@article{YBZFK38,
author = {Zeldovich, Y. B. and Frank-Kamenetsky, D. A.},
title = {A theory of thermal propagation of flame},
journal = {Acta Physicochim. URSS}, 
volume = {9},
pages = {341-350},
year = {1938}
}

@article{NW,
author = {Newell, A. C. and Whitehead, J. A.},  title = {Finite bandwidth, finite amplitude convection},
journal = {J. Fluid Mech.}, 
volume = {38},
pages = {3118-139},
year = {1987}
}

@article{FIS,
author = {Fisher, R. A.},
title = {The wave of advance of advantageous genes},
journal = {Ann. Eugen.},
volume = {7},
number = {4},
pages = {355-369},
doi = {https://doi.org/10.1111/j.1469-1809.1937.tb02153.x},
url = {https://onlinelibrary.wiley.com/doi/abs/10.1111/j.1469-1809.1937.tb02153.x},
eprint = {https://onlinelibrary.wiley.com/doi/pdf/10.1111/j.1469-1809.1937.tb02153.x},
year = {1937}
}

@article {KianYamamoto:2019,
    AUTHOR = {Kian, Yavar and Yamamoto, Masahiro},
     TITLE = {Reconstruction and stable recovery of source terms and
              coefficients appearing in diffusion equations},
   JOURNAL = {Inverse Problems},
  FJOURNAL = {Inverse Problems. An International Journal on the Theory and
              Practice of Inverse Problems, Inverse Methods and Computerized
              Inversion of Data},
    VOLUME = {35},
      YEAR = {2019},
    NUMBER = {11},
     PAGES = {115006, 24 pp.},
      ISSN = {0266-5611},
   MRCLASS = {35R30 (35B30 35B35 35R11)},
  MRNUMBER = {4036388},
       DOI = {10.1088/1361-6420/ab2d42},
       URL = {https://doi.org/10.1088/1361-6420/ab2d42},
}

@article {ZhangZhangZhou:2022,
    AUTHOR = {Zhang, Zhengqi and Zhang, Zhidong and Zhou, Zhi},
     TITLE = {Identification of potential in diffusion equations from
              terminal observation: analysis and discrete approximation},
   JOURNAL = {SIAM J. Numer. Anal.},
  FJOURNAL = {SIAM Journal on Numerical Analysis},
    VOLUME = {60},
      YEAR = {2022},
    NUMBER = {5},
     PAGES = {2834--2865},
      ISSN = {0036-1429},
   MRCLASS = {65M30 (65M15 65M32 65M60 65R20)},
  MRNUMBER = {4492543},
MRREVIEWER = {Zhu Wang},
       DOI = {10.1137/21M1446708},
       URL = {https://doi.org/10.1137/21M1446708},
}

@article {Isakov:1991,
    AUTHOR = {Isakov, Victor},
     TITLE = {Inverse parabolic problems with the final overdetermination},
   JOURNAL = {Comm. Pure Appl. Math.},
  FJOURNAL = {Communications on Pure and Applied Mathematics},
    VOLUME = {44},
      YEAR = {1991},
    NUMBER = {2},
     PAGES = {185--209},
      ISSN = {0010-3640},
   MRCLASS = {35R30 (35K10 35N10)},
  MRNUMBER = {1085828},
MRREVIEWER = {Dang Dinh Ang},
       DOI = {10.1002/cpa.3160440203},
       URL = {https://doi.org/10.1002/cpa.3160440203},
}

@article {LiSun:2017,
    AUTHOR = {Li, Buyang and Sun, Weiwei},
     TITLE = {Maximal {$L^p$} analysis of finite element solutions for
              parabolic equations with nonsmooth coefficients in convex
              polyhedra},
   JOURNAL = {Math. Comp.},
  FJOURNAL = {Mathematics of Computation},
    VOLUME = {86},
      YEAR = {2017},
    NUMBER = {305},
     PAGES = {1071--1102},
      ISSN = {0025-5718},
   MRCLASS = {65M60 (35K20 65M12)},
  MRNUMBER = {3614012},
MRREVIEWER = {B\"{u}lent Karas\"{o}zen},
       DOI = {10.1090/mcom/3133},
       URL = {https://doi.org/10.1090/mcom/3133},
}

@book {Thomee:2006,
    AUTHOR = {Thom\'{e}e, Vidar},
     TITLE = {Galerkin {F}inite {E}lement {M}ethods for {P}arabolic {P}roblems},
   EDITION = {Second},
 PUBLISHER = {Springer-Verlag, Berlin},
      YEAR = {2006},
     PAGES = {xii+370},
      ISBN = {978-3-540-33121-6; 3-540-33121-2},
   MRCLASS = {65-02 (65M15 65M60)},
  MRNUMBER = {2249024},
}

@article {JinShinZhou:2023,
    AUTHOR = {Jin, Bangti and Shin, Kwancheol and Zhou, Zhi},
     TITLE = {Numerical recovery of a time-dependent potential in
              subdiffusion},
   JOURNAL = {Inverse Problems},
  FJOURNAL = {Inverse Problems. An International Journal on the Theory and
              Practice of Inverse Problems, Inverse Methods and Computerized
              Inversion of Data},
    VOLUME = {40},
      YEAR = {2024},
    NUMBER = {2},
     PAGES = {025008, 34 pp.},
      ISSN = {0266-5611,1361-6420},
   MRCLASS = {65M32 (65M60)},
  MRNUMBER = {4687005},
       DOI = {10.1088/1361-6420/ad14a0},
       URL = {https://doi.org/10.1088/1361-6420/ad14a0},
}

@book {JinZhou:2023book,
    AUTHOR = {Jin, Bangti and Zhou, Zhi},
     TITLE = {Numerical {T}reatment and {A}nalysis of {T}ime-{F}ractional {E}volution
              {E}quations},  
 PUBLISHER = {Springer, Cham},
      YEAR = {2023},
     PAGES = {xiii+427},
      ISBN = {978-3-031-21049-5; 978-3-031-21050-1},
   MRCLASS = {65-02 (26A33 35R11 65M99)},
  MRNUMBER = {4572861},
       DOI = {10.1007/978-3-031-21050-1},
       URL = {https://doi.org/10.1007/978-3-031-21050-1},
}

@book {Brenner:2008,
    AUTHOR = {Brenner, Susanne C. and Scott, L. Ridgway},
     TITLE = {{The {M}athematical {T}heory of {F}inite {E}lement {M}ethods}},
   EDITION = {Third},
 PUBLISHER = {Springer, New York},
      YEAR = {2008},
     PAGES = {xviii+397},
      ISBN = {978-0-387-75933-3},
   MRCLASS = {65-01 (65-02)},
  MRNUMBER = {2373954},
       DOI = {10.1007/978-0-387-75934-0},
       URL = {https://doi.org/10.1007/978-0-387-75934-0},
}

@article{BH,
  author  = {Behzadan, A. and Holst, M.},
  title   = {Multiplication in {Sobolev} spaces, revisited},
  journal = {Ark. Mat.},
  volume  = {59},
  number  = {2},
  pages   = {275--306},
  year    = {2021}
}

@book{CH,
  author    = {Cazenave, T. and Haraux, A.},
  title     = {{An Introduction to Semilinear Evolution Equations}},  
  publisher = {Clarendon Press},
  address   = {Oxford},
  year      = {1998}
}

@article{CGR,
  author  = {Caro, P. and Garc{\'i}a-Ferrero, M. A. and Rogers, K. M.},
  title   = {Reconstruction for the {Calder{\'o}n} problem with {Lipschitz} conductivities},
  journal = {Anal. PDE},
  volume  = {18},
  number  = {8},
  pages   = {2033--2060},
  year    = {2025}
}

@article{CK2,
  author  = {Choulli, M. and Kian, Y.},
  title   = {Logarithmic stability in determining the time-dependent zero order coefficient in a parabolic equation from a partial {Dirichlet-to-Neumann} map. {Application} to the determination of a nonlinear term},
  journal = {J. Math. Pures Appl.},
  volume  = {114},
  pages   = {235--261},
  year    = {2018}
}

@misc{CJKZ:2026,
  author        = {Cen, S. and Jin, B. and Kian, Y. and Zhou, Z.},
  title         = {{Numerical Analysis of Space-Time Dependent Source Identification in Subdiffusion Equations}},
   howpublished = {Preprint, arXiv:2405.10708},
  year          = {2026} 
}

@article{Fei,
  author  = {Feizmohammadi, A.},
  title   = {An inverse boundary value problem for isotropic nonautonomous heat flows},
  journal = {Math. Ann.},
  volume  = {388},
  pages   = {1569--1607},
  year    = {2024}
}

@article{FKOU,
  author  = {Feizmohammadi, A. and Krupchyk, K. and Oksanen, L. and Uhlmann, G.},
  title   = {Reconstruction in the {Calder{\'o}n} problem on conformally transversally anisotropic manifolds},
  journal = {J. Funct. Anal.},
  volume  = {281},
  number  = {9},
  pages   = {109191, 25 pp.},
  year    = {2021}
}

@article{GK,
  author  = {Gaitan, P. and Kian, Y.},
  title   = {A stability result for a time-dependent potential in a cylindrical domain},
  journal = {Inverse Problems},
  volume  = {29},
  number  = {6},
  pages   = {065006, 18 pp.},
  year    = {2013}
}

@article{Ga,
  author  = {Garde, H.},
  title   = {Reconstruction of piecewise constant layered conductivities in electrical impedance tomography},
  journal = {Commun. Partial Differ. Equations},
  volume  = {45},
  pages   = {1118--1133},
  year    = {2020}
}

@book{H,
  author    = {Henry, D.},
  title     = {{Geometric Theory of Semilinear Parabolic Equations}},
  publisher = {Springer-Verlag},
  address   = {Berlin},
  year      = {1981}
}

@article{I1,
  author  = {Isakov, V.},
  title   = {Completeness of products of solutions and some inverse problems for {PDE}},
  journal = {J. Diff. Equat.},
  volume  = {92},
  pages   = {305--316},
  year    = {1991}
}

@article{I2,
  author  = {Isakov, V.},
  title   = {On uniqueness in inverse problems for semilinear parabolic equations},
  journal = {Arch. Rat. Mech. Anal.},
  volume  = {124},
  pages   = {1--12},
  year    = {1993}
}

@article{JKZ,
  author  = {Jin, B. and Kian, Y. and Zhou, Z.},
  title   = {Reconstruction of a space-time-dependent source in subdiffusion models via a perturbation approach},
  journal = {SIAM J. Math. Anal.},
  volume  = {53},
  pages   = {4445--4473},
  year    = {2021}
}

@article{KLLY,
  author  = {Kian, Y. and Li, Z. and Liu, Y. and Yamamoto, M.},
  title   = {The uniqueness of inverse problems for a fractional equation with a single measurement},
  journal = {Math. Ann.},
  volume  = {380},
  number  = {3-4},
  pages   = {1465--1495},
  year    = {2021}
}

@article{KiUh,
 author = {Kian, Yavar and Uhlmann, Gunther},
 title = {Recovery of nonlinear terms for reaction diffusion equations from boundary measurements},
 fjournal = {Archive for Rational Mechanics and Analysis},
 journal = {Arch. Ration. Mech. Anal.},
 issn = {0003-9527},
 volume = {247},
 number = {1},
 pages = {6, 20 pp.}, 
 year = {2023}, 
}

@article{KRS,
 author = {Koch, Herbert and R{\"u}land, Angkana and Salo, Mikko},
 title = {On instability mechanisms for inverse problems},
 fjournal = {Ars Inveniendi Analytica},
 journal = {Ars Inven. Anal.},
 issn = {2769-8505},
 volume = {2025},
 pages = {1, 94 pp.}, 
 year = {2025} 
}

@article{LeykekhmanLi:2021,
  author  = {Leykekhman, D. and Li, B.},
  title   = {Weak discrete maximum principle of finite element methods in convex polyhedra},
  journal = {Math. Comput.},
  volume  = {90},
  pages   = {1--18},
  year    = {2021}
}

@book{LM1,
  author    = {Lions, J.-L. and Magenes, E.},
  title     = {{Non-Homogeneous Boundary Value Problems and Applications}},
  volume    = {I},
  publisher = {Springer-Verlag},
  address   = {Berlin},
  year      = {1972}
}

@book{LM2,
 author = {Lions, J. L. and Magenes, E.},
 title = {Non-homogeneous boundary value problems and applications. {Vol}. {II}. {Translated} from the {French} by {P}. {Kenneth}},
 fseries = {Grundlehren der Mathematischen Wissenschaften},
 series = {Grundlehren Math. Wiss.},
 issn = {0072-7830},
 volume = {182},
 year = {1972},
 publisher = {Springer, Cham},
 language = {English},
 zbMATH = {3360568},
 Zbl = {0227.35001}
}

@article{Na,
 author = {Nachman, Adrian I.},
 title = {Reconstructions from boundary measurements},
 fjournal = {Annals of Mathematics. Second Series},
 journal = {Ann. Math. (2)},
 issn = {0003-486X},
 volume = {128},
 number = {3},
 pages = {531--576},
 year = {1988} 
}

@article{Pennes:1948,
  author  = {Pennes, H. H.},
  title   = {Analysis of tissue and arterial blood temperatures in the resting human forearm},
  journal = {J. Appl. Physiol.},
  volume  = {1},
  pages   = {93--122},
  year    = {1948}
}

@book{Pietro:2012book,
  author    = {Di Pietro, D. A. and Ern, A.},
  title     = {{Mathematical Aspects of Discontinuous Galerkin Methods}},
  publisher = {Springer},
  address   = {Berlin},
  year      = {2012}
}

@article{ScottZhang:1990,
  author  = {Scott, L. R. and Zhang, S.},
  title   = {Finite element interpolation of nonsmooth functions satisfying boundary conditions},
  journal = {Math. Comput.},
  volume  = {54},
  pages   = {483--493},
  year    = {1990}
}

@article{SU,
  author  = {Sylvester, J. and Uhlmann, G.},
  title   = {A global uniqueness theorem for an inverse boundary value problem},
  journal = {Ann. of Math.},
  volume  = {125},
  pages   = {153--169},
  year    = {1987}
}

@article{Trucu:2010,
  author  = {Trucu, D. and Ingham, D. B. and Lesnic, D.},
  title   = {Space-dependent perfusion coefficient identification in the transient bio-heat equation},
  journal = {J. Engrg. Math.},
  volume  = {67},
  pages   = {307--315},
  year    = {2010}
}

@article{Yamamoto:2001,
  author  = {Yamamoto, M. and Zou, J.},
  title   = {Simultaneous reconstruction of the initial temperature and heat radiative coefficient},
  journal = {Inverse Problems},
  volume  = {17},
  pages   = {1181--1202},
  year    = {2001}
}

@article{YGD,
  author  = {Ye, H. and Gao, J. and Ding, Y.},
  title   = {A generalized {Gronwall} inequality and its application to a fractional differential equation},
  journal = {J. Math. Anal. Appl.},
  volume  = {328},
  pages   = {1075--1081},
  year    = {2007}
}

@article {Zlamal:1973,
    AUTHOR = {Zl{\'a}mal, Milo{\v s}},
     TITLE = {Curved elements in the finite element method. {I}},
   JOURNAL = {SIAM J. Numer. Anal.},
  FJOURNAL = {SIAM Journal on Numerical Analysis},
    VOLUME = {10},
      YEAR = {1973},
     PAGES = {229--240},
      ISSN = {0036-1429},
   MRCLASS = {65N10},
  MRNUMBER = {395263},
MRREVIEWER = {Gilbert\ Strang},
       DOI = {10.1137/0710022},
       URL = {https://doi.org/10.1137/0710022},
}

@article {Zlamal:1974,
    AUTHOR = {Zl{\'a}mal, Milo{\v s}},
     TITLE = {Curved elements in the finite element method. {II}},
   JOURNAL = {SIAM J. Numer. Anal.},
  FJOURNAL = {SIAM Journal on Numerical Analysis},
    VOLUME = {11},
      YEAR = {1974},
     PAGES = {347--362},
      ISSN = {0036-1429},
   MRCLASS = {65N30},
  MRNUMBER = {343660},
MRREVIEWER = {Gilbert\ Strang},
       DOI = {10.1137/0711031},
       URL = {https://doi.org/10.1137/0711031},
}

@article {Lenoir:1986,
    AUTHOR = {Lenoir, M.},
     TITLE = {Optimal isoparametric finite elements and error estimates for
              domains involving curved boundaries},
   JOURNAL = {SIAM J. Numer. Anal.},
  FJOURNAL = {SIAM Journal on Numerical Analysis},
    VOLUME = {23},
      YEAR = {1986},
    NUMBER = {3},
     PAGES = {562--580},
      ISSN = {0036-1429},
   MRCLASS = {65N15 (65N30)},
  MRNUMBER = {842644},
MRREVIEWER = {Stephen\ W.\ Brady},
       DOI = {10.1137/0723036},
       URL = {https://doi.org/10.1137/0723036},
}

@article {Ciarlet:1972,
    AUTHOR = {Ciarlet, P. G. and Raviart, P.-A.},
     TITLE = {Interpolation theory over curved elements, with applications
              to finite element methods},
   JOURNAL = {Comput. Methods Appl. Mech. Engrg.},
  FJOURNAL = {Computer Methods in Applied Mechanics and Engineering},
    VOLUME = {1},
      YEAR = {1972},
     PAGES = {217--249},
      ISSN = {0045-7825,1879-2138},
   MRCLASS = {65N15 (65D05)},
  MRNUMBER = {375801},
MRREVIEWER = {Yu.\ K.\ Dem\cprime yanovich},
       DOI = {10.1016/0045-7825(72)90006-0},
       URL = {https://doi.org/10.1016/0045-7825(72)90006-0},
}

\end{document}